\documentclass[12pt]{amsart}
\usepackage{amsmath}
\usepackage{amssymb}
\usepackage{mathrsfs}
\usepackage{amsfonts}
\usepackage{amssymb}
\usepackage{amsfonts, amscd, amsmath, mathrsfs, amssymb, amsthm, amsxtra, bbding, epsfig, graphicx, latexsym, url, mathbbol, bbold}
\usepackage[papersize={8.7in,12in},textwidth=18cm,textheight=24.4cm,centering]{geometry}
\usepackage{enumerate}
\usepackage{caption} 
\usepackage{bm}

\usepackage{graphicx}
\usepackage{tikz}
\usepackage{float}

\usepackage{xcolor}
\definecolor{cite}{rgb}{0.00,0.60,1.00}
\definecolor{url}{rgb}{1.00,0.10,0.80}
\definecolor{link}{rgb}{0.00,0.00,1.00}
\usepackage[colorlinks,linkcolor=link,urlcolor=url,citecolor=cite,pagebackref,breaklinks]{hyperref}

\hypersetup{
	pdfstartpage=1,
	pdfstartview=FitH}

\DeclareFontFamily{U}{mathx}{\hyphenchar\font45}
\DeclareFontShape{U}{mathx}{m}{n}{
	<5> <6> <7> <8> <9> <10>
	<10.95> <12> <14.4> <17.28> <20.74> <24.88>
	mathx10
}{}
\DeclareSymbolFont{mathx}{U}{mathx}{m}{n}
\DeclareMathAccent{\widecheck}{\mathalpha}{mathx}{"71}

\numberwithin{equation}{section}

\allowdisplaybreaks

\newtheorem{theorem}{Theorem}[section]
\newtheorem{lemma}{Lemma}[section]

\newtheorem{definition}{Definition}[section]
\newtheorem{proposition}{Proposition}[section]

\newtheorem{corollary}{Corollary}[section]

\newtheorem{condition}{Condition}[section]

\makeatletter
\newcounter{roem}
\renewcommand{\theroem}{\Roman{roem}}

\newcommand{\c@org@eq}{}
\let\c@org@eq\c@equation
\newcommand{\org@theeq}{}
\let\org@theeq\theequation

\newcommand{\setroem}{
	\let\c@equation\c@roem
	\let\theequation\theroem}

\newcommand{\setarab}{
	\let\c@equation\c@org@eq
	\let\theequation\org@theeq}
\makeatother

\newtheorem*{claim*}{Claim}

\theoremstyle{theorem}
\newtheorem{remark}{\bf Remark}

\newcommand{\ud}{\mathrm{d}}
\newcommand{\ue}{\mathrm{e}}

\newcommand{\kl}{\mathrm{Kl}}

\DeclareMathOperator{\Mod}{mod}

\renewcommand{\bmod}[1]{\,(\Mod{ #1})}

\newcommand{\bh}{\mathbf{h}}

\newcommand{\balpha}{\boldsymbol{\alpha}}
\newcommand{\bbeta}{\boldsymbol{\beta}}
\newcommand{\bgamma}{\boldsymbol{\gamma}}
\newcommand{\bdelta}{\boldsymbol{\delta}}
\newcommand{\blambda}{\boldsymbol{\lambda}}
\newcommand{\brho}{\boldsymbol{\rho}}

\newcommand{\bC}{\mathbf{C}}

\newcommand{\bR}{\mathbf{R}}
\newcommand{\bZ}{\mathbf{Z}}

\newcommand{\sA}{\mathscr{A}}

\newcommand{\sC}{\mathscr{C}}

\newcommand{\sG}{\mathscr{G}}
\newcommand{\sK}{\mathscr{K}}
\newcommand{\sL}{\mathscr{L}}
\newcommand{\sM}{\mathscr{M}}
\newcommand{\sN}{\mathscr{N}}
\newcommand{\sP}{\mathscr{P}}
\newcommand{\sR}{\mathscr{R}}

\newcommand{\sT}{\mathscr{T}}

\newcommand{\cB}{\mathcal{B}}
\newcommand{\cC}{\mathcal{C}}

\newcommand{\cE}{\mathcal{E}}

\newcommand{\cI}{\mathcal{I}}
\newcommand{\cJ}{\mathcal{J}}
\newcommand{\cK}{\mathcal{K}}
\newcommand{\cL}{\mathcal{L}}

\newcommand{\cP}{\mathcal{P}}

\newcommand{\cT}{\mathcal{T}}
\newcommand{\cU}{\mathcal{U}}

\newcommand{\fa}{\mathfrak{a}}

\newcommand{\fS}{\mathfrak{S}}
\newcommand{\fR}{\mathfrak{R}}

\def\leq{\leqslant}

\newcommand{\lqs}{\leqslant}
\newcommand{\gqs}{\geqslant}
\usepackage{graphicx}
\usepackage{tikz}

\begin{document}
	
	\title{Primes in simultaneous arithmetic progressions}
	
	\author{Zongkun Zheng}
	
	\address{School of Mathematics and Statistics, Xi'an Jiaotong University, Xi'an 710049, P. R. China}

	\email{zongkunzheng@stu.xjtu.edu.cn}
	
	\subjclass[2020]{11N05, 11N13, 11N36, 11T23, 11L05, 11L07, 11N75}
	
	\keywords{primes in arithmetic progressions, dispersion method, sums of Kloosterman sums, van der Corput method, shifted Chen primes, linear sieve}

	\begin{abstract} 
		We prove a new mean value theorem on the distribution of primes in two simultaneous arithmetic progressions. Our approach builds on previous arguments of Bombieri, Fouvry, Friedlander, and Iwaniec appealing to spectral theory of Kloosterman sums, as well as the $q$-analogue of van der Corput method. In particular, we need estimates for exponential sums coming from the spectral theory of automorphic forms (sums of Kloosterman sums) and from algebraic geometry (Weil--Deligne bound for algebraic exponential sums). As an application, we show that the greatest prime factor of $p + 6$ for Chen prime $p$ is infinitely often greater than $p^{0.217}$.
	\end{abstract}
	
	\maketitle
	
	\setcounter{tocdepth}{1}
	\section{Introduction}\label{Introduction}
	
	\subsection{Primes in arithmetic progressions}\label{primesin}
	For $(q,a)=1,$ let
	\begin{align*}
		\pi(x;q,a):=\#\{p\lqs x:p\equiv a\bmod q\}
	\end{align*}
	count the number of primes up to $x$ which are congruent to $a$ mod $q$, and write $\pi(x;1,1)=\pi(x)$. The Siegel--Walfisz theorem \cite{Wa36} provides the pointwise asymptotic 
	$$\pi(x;q,a)\sim \frac{\pi(x)}{\varphi(q)}$$
	in the range $q \lqs \log^B x$ for any fixed $B>0$.
	Assuming Generalized Riemann Hypothesis (GRH), this range can be extended to $q \lqs x^{1/2-\varepsilon}$.
	However, in many applications, it suffices to establish equidistributions for almost all $q\lqs Q$. Specifically, for some large $A>0$, we expect bounds of the form
	\begin{align}\label{eq:bvtype}
		\sum_{\substack{q\lqs Q\\(q,a)=1}}\gamma_q\Big(\pi(x;q,a)-\frac{\pi(x)}{\varphi(q)}\Big)\ll \frac{x}{\log^Ax}
	\end{align}
	with some weight $\bgamma=(\gamma_q).$
	The celebrated Bombieri--Vinogradov theorem \cite{Bo65, Vi65} establishes this for $Q\lqs x^{1/2-\varepsilon}$, with an arbitrary weight $\bgamma$ and a residue class $a$ that may depend on $q$. On average over $q$, this result matches what would be expected under GRH. Elliott and Halberstam \cite{EH68} conjectured that (\ref{eq:bvtype}) holds for all $Q \lqs x^{1-\varepsilon}$.

	Remarkably, through the work of Fouvry--Iwaniec \cite{FI83}, Fouvry \cite{Fo84a}, and Bombieri--Friedlander-Iwaniec \cite{BFI86}, the $x^{1/2}$-barrier has been broken under additional restrictions: a fixed residue class $a$ (or of slow growth with $x$) and a well-factorable weight $\bgamma$. In his breakthrough on bounded gaps between primes, Zhang \cite{Zh14} broke the $x^{1/2}$-barrier with $\bgamma$ supported on smooth numbers. Recent work by Maynard \cite{Ma25b} provides equidistribution results for moduli up to $x^{3/5-\varepsilon}$ and $x^{7/12-\varepsilon}$, with triply-well-factorable weights and upper-bound linear sieve weights, respectively.
	Then for triply-well-factorable weights, Lichtman \cite{Li23} used optimized Deshouillers--Iwaniec estimates via Kim--Sarnak’s bound towards Ramanujan conjecture to improve the exponent to $66/107-\varepsilon$. More recently, Pascadi \cite{Pa24,Pa25a} proved large sieve inequalities for exceptional Maass forms with special sequences, and thus was able to improve Maynard's exponents to $5/8-\varepsilon$ and $3/5-\varepsilon$, respectively. For the residue $a$ that might be very large compared to $x$, Assing--Blomer--Li \cite{ABL21} extended the range to $x^{1/2+\delta}$ for some $\delta>0.$
	
	In this paper, we consider the distribution of primes $p$ satisfying
	$p\equiv a_1\bmod d$ and $p\equiv a_2\bmod q$ simultaneously. If $(d,q)=1,$
	the two congruence conditions are equivalent to $p\equiv a_1q\overline{q}+a_2d\overline{d}\bmod{dq}$ by the Chinese remainder theorem. Furthermore, if $a_1=a_2,$ the residue class becomes $a_1\bmod {dq}$, and one is in a good position to apply the above existing equidistributions. While for $a_1\neq a_2$, the residue class depends heavily on $d$ and $q$, which case is our main concern in this paper.
	
	To state the main results, we first recall the standard definition that a sequence $(a_n)$ is said to be {\it divisor-bounded} if $|a_n|\leqslant \tau(n)^C$ for some absolute constant $C.$ We also recall the convention of {\it well-factorable} originally introduced by Fouvry and Iwaniec \cite{FI83}.
	
	\begin{definition}[Well-factorable] Let $D\geqslant1.$ We say that a sequence $(\lambda_d)$ is well-factorable of level $D$ if$,$ for any factorization $D=D_1D_2$ with $D_1,D_2\gqs 1,$ there exist two sequences $\lambda'(d_1),\lambda''(d_2)$ such that
		\begin{itemize}
			\item $|\lambda'(d_1)|,|\lambda''(d_2)|\lqs 1$ for all $d_1,d_2;$
			\item $\lambda'(d_1),\lambda''(d_2)$ are supported on $1\lqs d_1\lqs D_1,1\lqs d_2\lqs D_2,$ respectively$;$
			\item For any $d\lqs D$ we have
			\begin{align*}
				\lambda_d=\sum_{d_1d_2=d}\lambda'(d_1)\lambda''(d_2).
			\end{align*}
		\end{itemize}
	\end{definition}
	
	We are now ready to formulate our main theorem.
	
	\begin{theorem}\label{thm:meanvalue-primes}
		Let $\varepsilon>0$ and $0\lqs\theta\leqslant 7/36.$ Suppose $|a_1|\lqs \log^B x$ for some $B>0,$ $|a_2|\leqslant x$ and $a_1\neq a_2$. 
		Let $\bgamma=(\gamma_q)$ and $\blambda=(\lambda_d)$ be divisor-bounded sequences supported on positive integers with
		\begin{align*}
			q\sim x^\theta,~d\lqs x^{\cL(\theta)-\varepsilon},~(d,a_1)=(q,a_2)=\mu^2(d)=\mu^2(q)=1.
		\end{align*}
		Assume further that $\blambda$ is well-factorable of level $x^{\cL(\theta)-\varepsilon}.$ Then for all $\cL(\theta)$ satisfying
		$$\cL(\theta)=\begin{cases}
			7/13-3\theta, & \theta\in [0,1/78], \\
			1/2, & \theta\in [1/78,1/40], \\
			7/13-20\theta/13, & \theta\in [1/40,1/35], \\
			19/36-20\theta/17, & \theta\in [1/35,17/192], \\
			295/576-\theta, & \theta\in [17/192,7/72],\\
			151/288-9\theta/8, & \theta\in [7/72,7/36],
		\end{cases}$$
		we have
		\begin{align}\label{eq:dq}
			\mathop{\sum_{q}\sum_{d}}_{(q,d)=1}\gamma_q\lambda_d\Big(\sum_{\substack{p\lqs x\\p\equiv a_1\bmod d\\p\equiv a_2\bmod q}}1-\frac{\pi(x)}{\varphi(qd)}\Big)\ll \frac{x}{\log^Ax} 
		\end{align}        
		for any $A>0,$ where the implied constant depends at most on $(\varepsilon,A,B).$
	\end{theorem}

	The arguments in proving Theorem \ref{thm:meanvalue-primes} also allow us to bound the following quadrilinear form.
	
	\begin{theorem}\label{thm:meanvalue-multilinearform}
		Let $\varepsilon>0,$ $0\lqs\theta\leqslant2/23$ and $0<\nu<1.$ Suppose $|a_1|\lqs \log^B x$ for some $B>0,$ $|a_2|\leqslant x$ and $a_1\neq a_2$. 
		Let $\balpha=(\alpha_m),\bbeta=(\beta_n),\bgamma=(\gamma_q)$ and $\blambda=(\lambda_d)$ be divisor-bounded sequences supported on positive integers with
		\begin{align*}
			&m\sim M=x^{1-\nu},~n\sim N=x^\nu,~q\sim x^\theta,~d\lqs x^{\cL(\theta,\nu)-\varepsilon},\\
			&(d,a_1)=(q,a_2)=\mu^2(d)=\mu^2(q)=1,~P^-(n)> \log^Cx,
		\end{align*}
		where $P^-(n)$ denotes the least prime factor of $n.$ 
		Assume further that $\blambda$ is well-factorable of level $x^{\cL(\theta)-\varepsilon}$ and  $\bbeta$ satisfies the Siegel--Walfisz condition $($Definition $\ref{Def:SW})$. Then for any $A>0,$ there exists some constant $C=C(A)>0$ such that
		\begin{align*}
			\mathop{\sum_{q}\sum_{d}}_{(q,d)=1}\gamma_q\lambda_d\bigg(\sum_{\substack{mn\equiv a_1\bmod d\\mn\equiv a_2\bmod q}}\alpha_m\beta_n-\frac{1}{\varphi(qd)}\sum_{\substack{(mn,qd)=1}}\alpha_m\beta_n\bigg)\ll \frac{x}{\log^Ax} 
		\end{align*}
		holds with an implied constant depending at most on $(\varepsilon,A,B),$ provided that one of the following conditions is valid$:$
		\begin{enumerate}[$(1)$]
			\item $\theta\lqs 1/60,$ $1/10 \lqs \nu\lqs 4/15$ and 
			$$\cL(\theta,\nu)=\begin{cases}
				(1+\nu)/2-3\theta,& \ \nu\in [1/10,(1+\theta)/5], \\
				(2-\nu)/3-17\theta/6,& \ \nu\in [(1+\theta)/5,\min\{(2-13\theta)/7,4/15\}], \\
				1-3\nu/2-5\theta,& \ \nu\in [(2-13\theta)/7,\max\{(2-13\theta)/7,4/15\}];
			\end{cases}$$
			
			\item $\theta\lqs 1/30,$ $31/90 \lqs \nu\lqs 83/220$ and 
			$$\cL(\theta,\nu)= (3+6\nu-12\theta)/10;$$
			
			\item $1/60\lqs \theta\lqs 2/23,$ $8/23 \lqs \nu\lqs 56/117$ and 
			$$\cL(\theta,\nu)= (16+5\nu-40\theta)/34.$$
		\end{enumerate}
	\end{theorem}

	\begin{remark}
		Note that if $\theta$ tends to zero $($or even $q=1$ for simplicity$)$, the problem in $(\ref{eq:dq})$ reduces to the type considered in $(\ref{eq:bvtype})$ with a single residue class, which has been studied extensively in \cite{FI83, Fo84a, BFI86, Li23, Ma25a, Ma25b, Pa25a}. However, the aim of this paper is not to substantially overcome the $x^{1/2}$-barrier, but to handle the new difficulty coming from the extra congruence condition. Therefore, for tiny $\theta$, we apply Heath-Brown identity as in \cite{FI83} to simplify the proof. It is worthwhile to mention that Pascadi \cite{Pa24} improved the result in Lemma $\ref{lm:kloostermania}$ below for dispersion coefficients by establishing a new large sieve inequality for the Fourier coefficient of exceptional Maass forms. Unfortunately, in our current approach we have some loss in order to make separations of variables arising from the extra congruence conditions, and Pascadi's refinement does not impact the key terms in our estimates for exponential sums. 
	\end{remark}

	As one may see from the subsequent arguments, we follow previous approaches of Bombieri, Fouvry, Friedlander, and Iwaniec appealing to spectral theory of Kloosterman sums. To do so, certain separations of variables become necessary. However, the loss might be quite substantial when $\theta$ in Theorem \ref{thm:meanvalue-primes} is suitably large. Therefore, as a complementary approach, we also employ the factorizations of moduli in algebraic exponential sums, building on the $q$-analogue of van der Corput method. See Sections \ref{firstmethod}, \ref{secondmethod} and \ref{thirdmethod} for details.
	
	It is highly desirable to explore more efficient approaches to make separations of variables, so that we have more opportunities to take advantage of the new ideas introduced in 
	\cite{Ma25a, Ma25b, Ma25c, Li23, Pa25a}.

	\subsection{An application to shifted Chen primes}\label{sec:shiftedChenprimes}
	The twin prime conjecture, a central problem in prime number theory, asserts the existence of infinitely many primes $p$ such that $p+2$ is also prime. A landmark result is due to Chen \cite{Ch73}, who proved that there are infinitely many primes $p$ such that $p+2\in\sP_2$,
	where we denote by $\sP_k$ the set of almost primes with at most $k$ prime factors. Such primes $p$ are called {\it Chen primes}.
	
	Quantitatively, Chen's theorem establishes that for all sufficiently large $x$, the number of Chen primes satisfies
	\begin{align}\label{eq:Chen}
		\#\{p \leq x : p + 2 \in \sP_2\} > 0.335 \sC_0 \frac{x}{(\log x)^2},
	\end{align}
	where $\sC_0$ is the twin prime constant defined in Proposition \ref{prop:Liu1}. The constant $0.335$ in the lower bound has been progressively increased through refinements of sieve methods and improvements on levels of distributions of primes in arithmetic progressions$:$
	
	\begin{itemize}
		\item Halberstam \cite{Ha75}: $0.3445$ 
		\item Chen \cite{Ch78a,Ch78b}: $0.3772$ and $0.405$ (more elaborated weighted sieve)
		\item Fouvry--Grupp \cite{FG86}: $0.71$ (leveraging the well-factorable structure of $\bgamma$ in $(\ref{eq:bvtype})$)
		\item Liu \cite{Li89}: $1.015$ 
		\item Wu \cite{Wu90}: $1.05$ 
		\item Cai \cite{Ca02}: $1.0974$ 
		\item Wu \cite{Wu04}: $1.104$ 
		\item Cai \cite{Ca08}: $1.13$
	\end{itemize}

	From another perspective to approximate the twin prime conjecture, one may consider the greatest prime factors of shifted primes. Goldfeld \cite{Go69} first observed, via the Brun--Titchmarsh theorem, that $P^+(p+a)>p^\delta$ holds infinitely often for some $\delta>1/2$, where $P^+(n)$ denotes the greatest prime factor of $n$. Subsequent progresses include:
	\begin{itemize}
		\item Motohashi \cite{Mo70}: $\delta=0.6105$
		\item Hooley \cite{Ho72}: $\delta=0.6199$
		\item Hooley \cite{Ho73}: $\delta=0.625-o(1)$ (via the Brun--Titchmarsh theorem on average)
		\item Deshouillers--Iwaniec \cite{DI81}: $\delta=0.6563$
		\item Fouvry \cite{Fo84b}: $\delta=0.6578$
		\item Fouvry \cite{Fo85}: $\delta=0.6687$ ($\delta>2/3$ relates to the first case of Fermat's last theorem)
		\item Baker--Harman \cite{BH96}: $\delta=0.676$
	\end{itemize}

	An extension of the twin prime conjecture predicts the infinitude of primes $p$ such that $p+2,p+6$ are both primes, as another special case of the Hardy--Littlewood conjecture for prime tuples in general (\cite{HL23}). Porter \cite{Po72} first established that the values $n(n+2)(n+6)$ are infinitely often in $\sP_8$. Maynard \cite{Ma13} improved Porter's bound from $8$ to $7$. Building on Chen's work, Heath-Brown and Li \cite{HL16} proved the existence of infinitely many Chen primes $p$ with $p+6\in\sP_{76}$, which was subsequently improved by Cai \cite{Ca17} to $\sP_{14}$.
	
	Inspired by these results, we study the greatest prime factor of $p+6$ for infinitely many Chen primes $p$, with the aid of the new type of mean value theorems established in Theorems \ref{thm:meanvalue-primes} and \ref{thm:meanvalue-multilinearform}.
	
	\begin{theorem}\label{thm:p+6}
		There are infinitely many Chen primes $p$ such that
		\begin{align*}
			P^+(p+6)>p^{0.217}.
		\end{align*}
	\end{theorem}
	
	The proof of Theorem \ref{thm:p+6} will be given in Section $\ref{greatest}$. A direct consequence of Cai \cite{Ca17} yields a weaker exponent $1/14$ in place of $0.217$ in Theorem $\ref{thm:p+6}$. As another application, Theorems $\ref{thm:meanvalue-primes}$ and $\ref{thm:meanvalue-multilinearform}$ might lead to an improvement over \cite{Ca17} by careful applications of weighted sieves.

	\subsection*{Acknowledgments}
	I am deeply grateful to Professor Ping Xi for his initial idea, his careful revisions of the manuscript and his consistently helpful comments. This work is supported in part by Shaanxi NSF (No. 2025JC-QYCX-002) and NSFC (No. 12025106).

	\section{Notation and preliminaries}
	\subsection{Notation and conventions}
	We write $\ue(t)=\exp({2\pi it})$ and $\sL=\log x$. The letter $p$ (with or without subscripts) is reserved for primes. 
	For $(n, q)$ = 1, we use $\overline{n}$ to denote the inverse of $n$ modulo $q$, and the modulus will be clear from the context. 
	To ease the presentation, we may sometimes leave out the coprime constraint $(n,q)=1$ in the summation of exponent like $\ue(\overline{n}/q)$. We will use $(a, b)$ and $[a,b]$ to denote, respectively, the greatest common factor and the least common multiple of $a$ and $b$ (or sometimes $(a,b,c)$ denotes the greatest common factor of $a,b$ and $c$) when it does not conflict with notation for ordered pairs. As most usual notation, we denote $\Lambda,\mu,\varphi,\tau$ as von Mangoldt function, M\"{o}bius function, Euler function and divisor function, respectively. Let $\omega(n)$ be the number of distinct prime factors of $n$, and $P^-(n)$ the least prime factor of $n$.
	
	We use the standard asymptotic notation $\ll$, $\gg$, $O(\cdot)$ and $o(\cdot)$ from analytic number theory. Let $x\sim X$ denote $X<x\lqs 2X$. We use $\varepsilon$ and $A$ to denote, respectively, a sufficiently small positive number and  a sufficiently large number, which might be different at each occurrence (sometimes we may write $x^{2\varepsilon}\sL^A \ll x^\varepsilon,$ for example).

	\begin{definition}[Siegel--Walfisz condition]\label{Def:SW}
		We say that a complex sequence $(\beta_n)$ satisfies the Siegel--Walfisz condition$,$ if for any $r\gqs 1,$ $q\gqs 1,$ $(a,q)=1,$ and $A>1,$ the inequality
		\begin{align*}
			\Big|\sum_{\substack{n\sim N, (n,r)=1\\n\equiv a\bmod q}}\beta_n-\frac{1}{\varphi(q)}\sum_{\substack{n\sim N\\(n,rq)=1}}\beta_n\Big|\ll_A \frac{N\tau^B(r)}{\log^AN}
		\end{align*}
		holds for some constant $B$.
	\end{definition}
	\begin{definition}[Fourier transform]
		Define the Fourier transform of a function $f:\bR\rightarrow\bC$ by
		\begin{align*} 
			\widehat{f}(v)=\int_{\bR} f(u)\ue(-vu)\ud u,
		\end{align*}
		and the Fourier transform of a function $g:\bZ/q\bZ\rightarrow\bC$ by
		\begin{align}\label{eq:finiteFouriertransform}
			\widehat{g}(v)=\frac{1}{\sqrt{q}}\sum_{u\bmod q}g(u)\ue\Big(\frac{-vu}{q}\Big).
		\end{align}
		We say $\widehat{g}$ is the Fourier transform of $g$ modulo $q.$
	\end{definition}
	
	\subsection{Combinatorial decompositions of primes} Recall that $\Lambda(n)=\log p$ if $n=p^k$ with $k\in\bZ^+$ and a prime $p$, and vanishes otherwise. The following identity of Heath-Brown \cite{HB82} allows one to decompose primes in a combinatorial manner, so that multilinear forms can be utilized to understand sums over primes.
	
	\begin{lemma}[Heath-Brown identity]\label{lm:HB} For any $n<2x,$ we have
		\begin{align}\label{eq:HB} 
			\Lambda(n)=\sum_{1\lqs j\lqs 4}(-1)^{j+1}\binom{4}{j}\sum_{m_1,\dots,m_j\lqs x^{1/4}}\mu(m_1)\cdots\mu(m_j)\sum_{m_1\cdots m_jn_1\cdots n_j=n}\log n_1 .
		\end{align} 
	\end{lemma}

	\begin{lemma}\label{lm:combine}
		Let $1/10<\sigma<1/2,$ and let $t_1,\dots,t_n\geqslant0$ such that $t_1+\cdots +t_n=1$. Then at least one of the following three conditions holds$:$\\
		$(1)$ There is a $t_i$ with $t_i\gqs 1/2+\sigma;$\\
		$(2)$ There is a partition $\{1,\dots,n\}=S\sqcup T$ such that
		\begin{align*}
			\frac{1}{2}-\sigma< \sum_{i\in S}t_i\lqs \sum_{i\in T}t_i<\frac{1}{2}+\sigma;
		\end{align*}
		$(3)$ There exist distinct $i,j,k$ with $2\sigma\lqs t_i\lqs t_j\lqs t_k\lqs 1/2-\sigma$ and
		\begin{align*}
			t_i+t_j,\,t_i+t_k,\,t_j+t_k\gqs \frac{1}{2}+\sigma.	
		\end{align*}
		Furthermore, if $\sigma > 1/6$, the third condition cannot occur.
	\end{lemma}
	\proof This is \cite[Lemma 3.1]{Po14}.\endproof
	\begin{lemma}\label{lm:FI}
		Let $7/24<\sigma<1/3,$ and let $t_1,\dots,t_n$ be nonnegative real numbers such that $t_1+\cdots +t_n=1$. Then at least one of the following two conditions holds$:$\\
		$(1)$ There is a $t_i$ with $t_i\gqs \sigma;$\\
		$(2)$ There is a subset $S\subset\{1,\dots,n\}$ such that
		\begin{align*}
			1-3\sigma< \sum_{i\in  S}t_i\lqs\frac{1}{4} .
		\end{align*}
	\end{lemma}
	\proof This follows from \cite[Lemma 4]{FI83}.\endproof

	\subsection{Congruences}
	We need the following reciprocity, which is an immediate consequence of the Chinese remainder theorem.
	
	\begin{lemma}[Bezout's identity]\label{lm:Bezout}
		Let $(r,s)=1$. For each integer $a,$ we have
		\begin{align*}
			\frac{a}{rs}\equiv\frac{a\overline{r}}{s}+\frac{a\overline{s}}{r}\bmod 1 .
		\end{align*}
	\end{lemma}

	\begin{lemma}\label{lm:change moduli}
		Let $q_0,q_1,q_2,r_0,r_1,r_2,t_1,t_2,s$ be pairwise coprime$,$ $(n_1,q_0q_1r_0r_1t_1s)=(n_2,q_0q_2r_0r_2t_2s)=(n_1,n_2)=1.$ Suppose $\varrho \bmod {q_0q_1q_2r_0r_1r_2t_1t_2s}$ is a common solution to the system  
		$$\begin{cases}
			\varrho  n_1\equiv a_1\bmod {r_0r_1t_1s},\\
			\varrho  n_2\equiv a_1\bmod {r_0r_2t_2s},\\
			\varrho  n_1\equiv a_2\bmod {q_0q_1t_2},\\
			\varrho  n_2\equiv a_2\bmod {q_0q_2t_1}.
		\end{cases}$$
		Assume that $n_1\equiv n_2\bmod{q_0r_0s},~a_2n_1\equiv a_1n_2\bmod{t_1}$ and $a_1n_1\equiv a_2n_2\bmod{t_2}.$ Then we have
		\begin{align*}
			\frac{\varrho}{q_0q_1q_2r_0r_1r_2t_1t_2s}
			&\equiv a_1\frac{n_1-n_2}{q_0r_0s}\frac{\overline{q_1r_1t_1n_2}}{q_2r_2t_2n_1}+(a_2-a_1)\frac{\overline{q_0q_2r_0r_1r_2t_1t_2sn_1}}{q_1}+(a_2-a_1)\frac{\overline{q_1r_0r_1r_2t_1t_2sn_2}}{q_0q_2}\notag\\
			&\ \ \ +\frac{a_1}{q_0q_1q_2r_0r_1r_2t_1t_2sn_1}\bmod 1.
		\end{align*}
	\end{lemma}
	\proof 
	The Chinese remainder theorem guarantees the solvability of the above system of congruence equations.
	From Lemma \ref{lm:Bezout} it follows that
	\begin{align*}
		\frac{\varrho}{q_0q_1q_2r_0r_1r_2t_1t_2s}
		&\equiv a_1\frac{\overline{q_0q_1q_2r_0r_2t_2sn_1}}{r_1t_1}+a_1\frac{\overline{q_0q_1q_2r_1t_1n_2}}{r_0r_2t_2s}+a_2\frac{\overline{q_0q_2r_0r_1r_2t_1t_2sn_1}}{q_1}+a_2\frac{\overline{q_1r_0r_1r_2t_1t_2sn_2}}{q_0q_2}\\
		&\equiv a_1\frac{\overline{q_0q_2r_0r_2t_2sn_1}}{q_1r_1t_1}+a_1\frac{\overline{q_1r_1t_1n_2}}{q_0q_2r_0r_2t_2s}
		+(a_2-a_1)\frac{\overline{q_0q_2r_0r_1r_2t_1t_2sn_1}}{q_1}\\
		&\ \ \ \ 
		+(a_2-a_1)\frac{\overline{q_1r_0r_1r_2t_1t_2sn_2}}{q_0q_2}\bmod1.
	\end{align*}
	Again, by Lemma \ref{lm:Bezout}, we transfer the first two fractions via
	\begin{align*}
		\frac{\overline{q_0q_2r_0r_2t_2sn_1}}{q_1r_1t_1}+\frac{\overline{q_1r_1t_1n_2}}{q_0q_2r_0r_2t_2s}
		&\equiv -\frac{\overline{q_1r_1t_1}}{q_0q_2r_0r_2t_2sn_1}+\frac{\overline{q_1r_1t_1n_2}}{q_0q_2r_0r_2t_2s}+\frac{1}{q_0q_1q_2r_0r_1r_2t_1t_2sn_1}\\
		&\equiv \frac{n_1-n_2}{q_0r_0s}\frac{\overline{q_1r_1t_1n_2}}{q_2r_2t_2n_1}+\frac{1}{q_0q_1q_2r_0r_1r_2t_1t_2sn_1}\bmod1,
	\end{align*}
	which finishes the proof.\endproof
	
	To see the transformation in Lemma \ref{lm:change moduli} more clearly, we assume $q_0=r_0=t_1=t_2=1,$ so that Lemma \ref{lm:change moduli} becomes
	\begin{align}\label{eq:easechange}
		\frac{\varrho}{q_1q_2r_1r_2s}
		&\equiv a_1\frac{n_1-n_2}{s}\frac{\overline{q_1r_1n_2}}{q_2r_2n_1}+(a_2-a_1)\frac{\overline{q_2r_1r_2sn_1}}{q_1}+(a_2-a_1)\frac{\overline{q_1r_1r_2sn_2}}{q_2}+\frac{a_1}{q_1q_2r_1r_2sn_1} \bmod1.
	\end{align}
	In practice, the above transformation in Lemma \ref{lm:change moduli} usually 
	appears in the treatment of exponential sums, when summing over $q_0,q_1,q_2,r_0,r_1,r_2,t_1,t_2,s\leqslant X,$ say. 
	We may impose the truncation
	\begin{align*}
		\max\{q_0,r_0,t_1,t_2\}\leqslant \log^CX
	\end{align*}
	with a suitably large $C>0,$ since otherwise one has very few terms in the summations, and trivial estimates for exponential sums would suffice.

	\subsection{Exponential sums}
	We start with the following Kloosterman sums defined over $\bZ/q\bZ:$
	\begin{align*}
		S(m,n;q)=\sideset{}{^*}\sum_{a\bmod q}\ue\Big(\frac{ma+n\overline{a}}{q}\Big).
	\end{align*}
	In particular, $S(m,0;q)$ is a Ramanujan sum, which satisfies the estimate
	\begin{align*}
		|S(m,0;q)|\lqs (m,q).
	\end{align*}
	We now recall Weil's bound for individual Kloosterman sums with general $m,n,q$.
	\begin{lemma}[Weil's bound]\label{lm:Weil} Let $m,n\in\bZ$ and $q\in\bZ^+.$ We have
		\begin{align*}
			|S(m,n;q)|\leqslant q^{1/2}(m,n,q)^{1/2}\tau(q).
		\end{align*}
	\end{lemma}
	\proof This result was first established by Weil \cite{We48} for prime $q$, and later generalized by Estermann \cite{Es61} to all $q\in\bZ^+.$
	\endproof

	Weil's bound is optimal in the sense that the exponent of $q$ in the upper bound cannot be replaced by any $\delta<1/2$. The following lemma presents non-correlations among different Kloosterman sums, which is stronger than the direct application of Weil's bound for individual sums.
	
	\begin{lemma}[Sum of products of Kloosterman sums]\label{lm:Kloostermanproduct}
		Let $p$ be a prime$,$ $(\xi,p)=1$ and let $s_1,\dots,s_k$ be integers$.$ Then we have 
		$$\sum_{t\bmod p}\ue\Big(\frac{bt}{p}\Big)\prod_{1\leqslant j\leqslant k}S(t+s_j,\xi;p)
		\ll \delta_p\cdot p^{k/2},$$	
		where $\delta_p=p$ if $b=0$ and all the $s_j$'s occur with even multiplicity modulo $p,$ and $\delta_p=p^{1/2}$ otherwise$.$
	\end{lemma}
	\proof This statement appeared on different occasions in various settings. It seems that \cite[Lemma 2.1]{FMRS04} is the earliest literature. See also \cite[Proposition 3.2]{FGKM14}, \cite[Lemma 4.4]{Ir15}, \cite[Corollary 1.6]{FKM15} or \cite[Lemma 4]{Xi17}.\endproof
	
	\begin{lemma}\label{lm:fourier}
		Let $q\in\bZ^+$ be square-free$,$ and $\xi,l_1,l_2\in\bZ.$ Denote by $Y$ the product of four Kloosterman sums$:$
		\begin{align*}
			Y(t)=S(t,\xi;q)S(t+l_1,\xi;q)S(t+l_2,\xi;q)S(t+l_1+l_2,\xi;q),
		\end{align*}
		and denote by $\widehat{Y}$ the Fourier transform of $Y$ modulo $q.$ Then we have
		\begin{align*}
			\widehat{Y}(v)\ll q^{2+\varepsilon}(l_1,v,q)^{1/2}(l_2,v,q)^{1/2}(\xi,l_1,l_2,q)^{3/2}.
		\end{align*}
	\end{lemma}
	\proof
	By the Chinese remainder theorem, we have
	\begin{align*}
		\widehat{Y}(v)&=\frac{1}{\sqrt{q}}\prod_{p\mid q}\sum_{t\bmod p}\eta(t,p),
	\end{align*}
	where
	\begin{align*}
		\eta(t,p)
		&=S(t,\xi_p;p)S(t+l_1,\xi_p;p)S(t+l_2,\xi_p;p)S(t+l_1+l_2,\xi_p;p)\ue\Big(\frac{-vt}{p}\Big)
	\end{align*}
	with $\xi_p\equiv \xi\overline{q/p}\bmod p.$
	
	For $p\nmid\xi$, we apply Lemma \ref{lm:Kloostermanproduct} to obtain
	\begin{align*}
		\sum_{t\bmod p}\eta(t,p)\ll p^{5/2}(l_1,v,p)^{1/2}(l_2,v,p)^{1/2}(l_1,l_2,v,p)^{-1/2}.
	\end{align*}
	
	For $p\mid\xi,$ we appeal to the estimates for Ramanujan sums to obtain 
	\begin{align*}
		|\eta(t,p)|\lqs (t,p)(t+l_1,p)(t+l_2,p)(t+l_1+l_2,p).
	\end{align*}
	Therefore, we observe the followings:
	\begin{itemize}
		\item $|\eta(t,p)|>1$ for at most four residue classes $t\Mod{p}.$
		\item If $p\nmid (l_1,l_2),$ then $|\eta(t,p)|\lqs p^2.$
		\item If $p\mid (l_1,l_2),$ then $|\eta(t,p)|\lqs p^4.$
	\end{itemize}
	Thus in the case $p\mid\xi,$ we have
	\begin{align*}
		\sum_{t\bmod p}\eta(t,p)\ll p^2+p^4\mathbf{1}_{p\mid(l_1,l_2)}\ll p^2(l_1,l_2,p)^2.
	\end{align*}
	
	We then finish the proof by combining both cases.
	\endproof

	Besides the classical Kloosterman sums, we also need the following hyper-Kloosterman sum:
	\begin{align*}
		\kl_3(a;q)= \frac{1}{q}\sum_{\substack{z_1,z_2,z_3\bmod q\\ z_1z_2z_3\equiv a\bmod q}}\ue\Big(\frac{z_1+z_2+z_3}{q}\Big). 
	\end{align*}
	From the Chinese remainder theorem and Deligne's bound for prime moduli (\cite{De77}), one has
	\begin{align}\label{eq:Deligne}
		|\kl_3(a;q)|\leqslant \tau_3(q).
	\end{align}
	for each square-free number $q$ with $(a,q)=1$.

	\begin{lemma}[Correlation of hyper-Kloosterman sums]\label{lm:correlation}
		Let $H\gqs 1,$ and $r,s_1,s_2$ be square-free numbers with $(r,s_1s_2)=1.$ For $(a_1,rs_1)=(a_2,rs_2)=1,$ we have
		\begin{align*}
			&\sum_{\substack{h\lqs H\\(h,rs_1s_2)=1}}\kl_3(a_1h;rs_1)\overline{\kl_3(a_2h;rs_2)}\\
			&\ll (Hrs_1s_2)^\varepsilon\Big(\frac{H}{r[s_1,s_2]}+1\Big)r^{1/2}[s_1,s_2]^{1/2}(a_2-a_1,s_1,s_2)^{1/2}(a_2s_1^3-a_1s_2^3,r)^{1/2} .
		\end{align*}
	\end{lemma}
	\proof This is \cite[Corollary 6.26]{Po14}. We have omitted the condition $H\ll (r[s_1,s_2])^{O(1)}$ here by incorporating a factor $H^\varepsilon$ in the bound.\endproof
	
	\begin{remark}
		In the typical case that $(s_1,s_2)=(a_2s_1^3-a_1s_2^3)=1,$ Lemma $\ref{lm:correlation}$ is non-trivial as long as $H>(rs_1s_2)^{1/2+\varepsilon},$ beyond the direct application of $\eqref{eq:Deligne}.$
	\end{remark}
	
	We also need the following estimate for complete algebraic exponential sums with general rational phases.

	\begin{lemma}\label{lm:algebraic type exponential sums}
		Let $\lambda=\lambda_1/\lambda_2$ with 
		\begin{align*}
			\lambda_1(x)=\sum_{0\lqs j\lqs d_1}r_jx^{j}\in\bZ[x],\ \lambda_2(x)=\sum_{0\lqs j\lqs d_2}t_jx^{j}\in\bZ[x].
		\end{align*}
		Define the degree of $\lambda$ by
		\begin{align*}
			d=d(\lambda)=\mathrm{deg}(\lambda_1)+\mathrm{deg}(\lambda_2)=d_1+d_2.
		\end{align*}
		Then for any $q\in\bZ^+$ satisfying $p^\alpha\parallel q\Rightarrow \alpha\lqs 2,$ we have
		\begin{align*}
			\sum_{\substack{a \bmod q\\(\lambda_2(a),q)= 1}}\ue\Big(\frac{\lambda(a)}{q}\Big)\ll q^{1/2}(\lambda,q_1)^{1/2}(\lambda,q_1)_*^{1/2}(\lambda',q_2)(2d)^{\omega(q)} ,
		\end{align*}
	\end{lemma}
	\noindent where
	\begin{align*}
		q_1=\prod_{p\parallel q}p,\ q_2=\prod_{p^2\parallel q}p,
	\end{align*}
	and
	\begin{align*}
		&(\lambda,q)_*= (r_0, r_1, r_2, \dots , r_{d_1}, q),\\
		&(\lambda,q)=(r_1, r_2, \dots , r_{d_1},t_1,t_2\dots t_{d_2}, q),\\
		&(\lambda',q)=(\lambda_1'\lambda_2-\lambda_1\lambda_2',q).
	\end{align*}
	\proof This is a special case of \cite[Theorem A.1]{WX21}.\endproof

	The bounds for algebraic exponential sums in Lemmas \ref{lm:Weil}, \ref{lm:Kloostermanproduct}, \ref{lm:fourier}, \ref{lm:correlation} and \ref{lm:algebraic type exponential sums} can be derived from $\ell$-adic cohomology. On the other hand, one can also appeal to the spectral theory of automorphic forms to study averages of complete/incomplete Kloosterman sums. This was initiated by Deshouillers and Iwaniec \cite{DI82} in the 1980s, and we recall one of their deep results as follows.

	\begin{lemma}[Kloostermania]\label{lm:kloostermania}
		Let $\Phi(\cdot,\cdot)$ be a smooth function with compact support on $\bR^+\times\bR^+$ and $U,V,W,C,D>0.$ For any complex number $\cB(u,v,w)$ we define
		\begin{align*}
			\sK(U,V,W,C,D)=\sum_{u\sim U}\sum_{v\sim V}\sum_{w\sim W}\cB(u,v,w)\sum_{c}\sum_{d}\Phi\Big(\frac{c}{C},\frac{d}{D}\Big)\ue\Big(w\frac{\overline{ud}}{vc}\Big).
		\end{align*}
		Then for any $\varepsilon>0$ we have
		\begin{align*}
			\sK(U,V,W,C,D)^2\ll(UVWCD)^\varepsilon\cK(U,V,W,C,D)\cdot \sum_{u\sim U}\sum_{v\sim V}\sum_{w\sim W}|\cB(u,v,w)|^2,
		\end{align*}
		where 
		\begin{align*}
			\cK(U,V,W,C,D)=CV(UV+W)(C+DU)+C^2DV\sqrt{(UV+W)U}+D^2WU, 
		\end{align*}
		and the constant implied in $\ll$ depends at most on $\varepsilon$ and $\Phi$.
	\end{lemma}
	
	\proof This is \cite[Theorem 12]{DI82}, corrected in \cite{BFI19}.
	\endproof

	\subsection{$q$-analogue of van der Corput method}
	
	Various incomplete algebraic exponential sums (like Lemma \ref{lm:correlation}) appear in analytic number theory, and one can apply the P\'olya--Vinogradov method to transform incomplete sums into complete ones, the latter of which might be controlled effectively with the aid of $\ell$-adic cohomology. This method works in many situations as long as the sums are not too short. On the other hand, if the modulus (of the algebraic exponential sum) has certain factorizations, one may beat the P\'olya--Vinogradov barrier using such factorizations and succeed in very short sums, as shown in the $q$-analogue of van der Corput method initiated by Heath-Brown \cite{HB78}.
	
	As in the classical van der Corput method, estimates for incomplete sums are derived from Weyl differencing and Poisson summation. These two techniques, known respectively as the $A$-process and the $B$-process, yield non-trivial bounds for short exponential sums (see \cite{WX21} for the general theory). We now introduce the relevant results concerning the incomplete sum
	\begin{align*}
		\sum_{n\in \cI}\Psi(n),
	\end{align*}
	where $\Psi:\bZ/q\bZ\rightarrow \bC$ is a function and $\cI$ is an interval.
	
	\begin{lemma}[$A$-process]\label{lm:A}
		Assume that $q=q_1q_2$ with $(q_1,q_2)=1$. For $i=1,2$ we define $\Psi_i:\bZ/q_i\bZ\rightarrow \bC$. Let $\Psi = \Psi_1\Psi_2.$ We then have
		\begin{align*}
			\Big|\sum_{n\in \cI}\Psi(n)\Big|^2\ll \frac{q_2}{|\cI|}\,\cU_2\Big(\cU_1+\sum_{0<l< |\cI|/q_2}\Big|\sum_{\substack{n,n+lq_2\in \cI}}\Psi_1(n)\overline{\Psi_1(n+lq_2)}\Big|\Big),
		\end{align*}
		where
		\begin{align*}
			\cU_1=\sum_{n\in\cI}\big|\Psi_1(n)\big|^2,\ \cU_2=\sum_{n\in\cI_1}\big|\Psi_2(n)\big|^2,
		\end{align*}
		with some interval $\cI_1\supseteq\cI$ such that $|\cI_1|\lqs 2|\cI|$.
	\end{lemma}
	\proof We follow the proof of \cite[Lemma 2.1]{Ir16}, which deals with the special case $\|\Psi_i\|_\infty\leqslant 1.$ If $|\cI|\lqs q_2$, the result follows from the trivial bound
	\begin{align*}
		\Big|\sum_{n\in \cI}\Psi(n)\Big|^2\lqs \cU_1\cU_2.
	\end{align*}
	Therefore, we assume that $L:=[|\cI|/q_2]\gqs 1$. By the decomposition $\Psi=\Psi_1\Psi_2$, and the fact that $\Psi_2$ is periodic of period $q_2$, we have
	\begin{align*}
		L\sum_{n\in \cI}\Psi(n)=\sum_{l\lqs L}\sum_{n+lq_2\in \cI}\Psi(n+lq_2)=\sum_{n\in \cI\cup(\cI-Lq_2)}\Psi_2(n)\sum_{\substack{l\lqs L\\n+lq_2\in \cI}}\Psi_1(n+lq_2),
	\end{align*}
	where $\cI-Lq_2=\{n-Lq_2:n\in\cI\}$. By Cauchy's inequality we have
	\begin{align*}
		L^2\Big|\sum_{n\in \cI}\Psi(n)\Big|^2&\lqs \sum_{n\in \cI\cup(\cI-Lq_2)}\big|\Psi_2(n)\big|^2\sum_{n\in \cI\cup(\cI-Lq_2)}\bigg|\sum_{\substack{l\lqs L\\n+lq_2\in \cI}}\Psi_1(n+lq_2)\bigg|^2\\
		&=\cU_2\sum_{l_1,l_2\lqs L}\sum_{\substack{n+l_1q_2\in\cI\\n+l_2q_2\in\cI}}\Psi_1(n+l_1q_2)\overline{\Psi_1(n+l_2q_2)}\\
		&=\cU_2\sum_{l_1,l_2\lqs L}\sum_{\substack{n\in\cI\\n+(l_2-l_1)q_2\in\cI}}\Psi_1(n)\overline{\Psi_1(n+(l_2-l_1)q_2)}.
	\end{align*}
	Thus, a change of variable $l=l_2-l_1$ gives
	\begin{align*}
		L^2\Big|\sum_{n\in \cI}\Psi(n)\Big|^2&\ll L\,\cU_2\sum_{0\lqs|l|\lqs L-1}\bigg|\sum_{n,n+lq_2\in \cI}\Psi_1(n)\overline{\Psi_1(n+lq_2)}\bigg|\\
		&= L\,\cU_2\bigg(\cU_1+\sum_{0<|l|< L}\bigg|\sum_{n,n+lq_2\in \cI}\Psi_1(n)\overline{\Psi_1(n+lq_2)}\bigg|\bigg),
	\end{align*}
	which deduces the final result.\endproof
	
	\begin{lemma}[$B$-process]\label{lm:B}
		There exist an $a\in \bZ$ and some interval $J$ not containing $0$ with $|J|\lqs q/|\cI|,$ such that 
		\begin{align*}
			\sum_{n\in \cI}\Psi(n)\ll \frac{|\cI||\widehat{\Psi}(0)|}{\sqrt{q}}+\frac{|\cI|\log q}{\sqrt{q}}\Big|\sum_{j\in J}\widehat{\Psi}(j)\ue\Big(\frac{aj}{q}\Big)\Big|,
		\end{align*}
		where $\widehat{\Psi}$ is the Fourier transform of $\Psi$ modulo $q$ as defined by $\eqref{eq:finiteFouriertransform}.$ 
	\end{lemma}
	\proof This is \cite[Lemma 3.1]{Ir16}.\endproof

	We also introduce another version of the $B$-process, which is less precise but is more efficient if the average is longer than the modulus $q.$ 
	\begin{lemma}\label{lm:Bcoro}
		We have
		\begin{align*}
			\sum_{n\in \cI}\Psi(n)\ll \frac{|\cI||\widehat{\Psi}(0)|}{\sqrt{q}}+\sqrt{q}\sum_{1\lqs |j|\lqs q/2}\frac{|\widehat{\Psi}(j)|}{|j|}.
		\end{align*}
	\end{lemma}
	\proof This is \cite[Lemma 3.2]{Ir16}.\endproof
	
	
	As one may see from Lemmas $\ref{lm:A}$ and $\ref{lm:B}$, the $A$- and $B$-processes transform the original sums to new averages, keeping in mind that the former reduces the modulus to a smaller one, and the latter does not alter the modulus but creates a dual sum (from a long sum to a short one, for instance). In some applications, we usually need suitable combinations of the two processes. For instance, it is possible to apply the $A$-process twice and then the $B$-process.
	
	\begin{lemma}[$A^2B$-process]\label{lm:Acoro}
		Suppose that $q=q_0q_1q_2.$ For $i=0,1,2,$ let $\Psi_i:\bZ/q_i\bZ\rightarrow \bC$. Let  $\Psi=\Psi_0\Psi_1\Psi_2$. Then we have
		\begin{align}\label{eq:Acoro1}
			\Big|\sum_{n\in \cI}\Psi(n)\Big|^4\ll &\frac{q^2_2}{|\cI|^2}\cU^2_{0,1}\,\cU^2_2+q_1\cU_0\cU_1\cU_2^2|\cI|^{\varepsilon-1}\notag\\
			&+\frac{q_0^{1/2}q_1q_2\cU_1\cU^2_2}{|\cI|}\sum_{\substack{0<|l_1|< |\cI|/q_1\\0<|l_2|< |\cI|/q_2}}\frac{1}{|\cI|-|l_2q_2|}\Big(\frac{|\cI||\widehat{Z}(0)|}{q_0}+\sum_{1\lqs |j|\lqs q_0/2}\frac{|\widehat{Z}(j)|}{|j|}\Big), 
		\end{align} 
		where $\cI_1,\cI_2$ are intervals satisfying $|\cI_1|,|\cI_2|\lqs 2|\cI|,$
		\begin{align*}
			&\cU_{0,1}=\sum_{j\in\cI}\big|\Psi_0(j)\Psi_1(j)\big|^2,\ \cU_0=\sum_{j\in\cI}\big|\Psi_0(j)\big|^2,\ 
			\cU_1=\sum_{j\in\cI_1}\big|\Psi_1(j)\big|^2,\ 
			\cU_2=\sum_{j\in\cI_2}\big|\Psi_2(j)\big|^2,\\
			&Z(j)=Z(j;l_1q_1,l_2q_2)=\Psi_0(j)\overline{\Psi_0(j+l_1q_1)}\overline{\Psi_0(j+l_2q_2)}\Psi_0(j+l_1q_1+l_2q_2),
		\end{align*}
		and $\widehat{Z}$ is the Fourier transform of $Z$ modulo $q_0.$ 
	\end{lemma}
	\proof This can be deduced by applying Lemma \ref{lm:A} twice and a single Lemma \ref{lm:B}. Related results can also be found in \cite[Lemma 2.2]{Ir16}, \cite[Theorem 8.1]{WX21}, and \cite[Theorem 8.3]{WX21}. \endproof
	
	\begin{remark}\label{briefly}
		If we assume $\Vert\Psi_i\Vert_\infty\ll 1$ for $i=0,1,2,$ and $\Vert \widehat{Z}\Vert_\infty\ll q_0^\varepsilon$, then Lemma $\ref{lm:Acoro}$ implies that
		\begin{align}\label{eq:Acoro2}
			\Big|\sum_{n\in \cI}\Psi(n)\Big|^4\ll q_0^\varepsilon(|\cI|^2q_2^2+|\cI|^3q_1+|\cI|^3q_0^{1/2}).
		\end{align}
		However, when applying the $A^2$-process in Section $\ref{thirdmethod}$, additional factors arise in the estimates for $\Psi_i.$ Although these factors are controlled by $\sL^A$ on average, they could be large pointwise. For this reason, we present the $A$-process in Lemma $\ref{lm:A}$ without assuming the condition $\Vert\Psi_i\Vert_\infty\ll 1$.
	\end{remark}

	\section{Dispersion method for the general bilinear form}\label{Dispersion}
	We first make trivial observation that to prove Theorem \ref{thm:meanvalue-primes}, it suffices to prove
	\begin{align*}
		\mathop{\sum_{q}\sum_{d}}_{(q,d)=1}\gamma_q\lambda_d\Big(\sum_{\substack{n\sim x\\P^-(n)>\sL^C\\n\equiv a_1\bmod d\\n\equiv a_2\bmod q}}\Lambda(n)-\frac{1}{\varphi(qd)}\sum_{\substack{n\sim x\\P^-(n)>\sL^C\\(n,qd)=1}}\Lambda(n)\Big)\ll x\sL^{-A}
	\end{align*} 
	for large constants $A,C>0.$ In what follows, one will see the advantages of the above new constraint $P^-(n)>\sL^C$. Roughly speaking, it would be convenient to study the least common multiple $[n_1,n_2]$ if $n_1n_2$ has no small prime factors. 
	
	By combining the two purely combinatorial results in Lemmas \ref{lm:combine} and \ref{lm:FI} with Heath-Brown identity (Lemma \ref{lm:HB}), and keeping in mind the applications of the well-factorable remainder terms of Iwaniec's linear sieve (Section \ref{bilinearsieve})$,$ we reduce the problem to establishing three types of mean value estimates$,$ namely Type I$,$ Type II and Type III estimates:
	\begin{align}
		\fS_{\textrm{I}}&:=\sum_{\substack{q,r,s\\ (q,rs)=1}}\gamma_q\delta_r\lambda_s\Big(\sum_{\substack{m,n\\P^-(mn)>\sL^C\\ mn\equiv a_1\bmod {rs}\\mn\equiv a_2\bmod q}}\beta_n-\frac{1}{\varphi(qrs)}\sum_{\substack{m,n\\P^-(mn)>\sL^C\\(mn,qrs)=1}}\beta_n\Big)\ll x\sL^{-A},\label{eq:TypeI}\\
		\fS_{\textrm{II}}&:=\sum_{\substack{q,r,s\\ (q,rs)=1}}\gamma_q\delta_r\lambda_s\bigg(\sum_{\substack{m,n\\P^-(mn)>\sL^C\\ mn\equiv a_1\bmod {rs}\\mn\equiv a_2\bmod q}}\alpha_m\beta_n-\frac{1}{\varphi(qrs)}\sum_{\substack{m,n\\P^-(mn)>\sL^C\\(mn,qrs)=1}}\alpha_m\beta_n\bigg)\ll x\sL^{-A},\label{eq:TypeII}\\
		\fS_{\textrm{III}}&:=\sum_{\substack{q,r,s\\ (q,rs)=1}}\gamma_q\delta_r\lambda_s\bigg(\sum_{\substack{m_1,m_2,m_3,n\\P^-(m_1m_2m_3n)>\sL^C\\ m_1m_2m_3n\equiv a_1\bmod {rs}\\m_1m_2m_3n\equiv a_2\bmod q}}\beta_n-\frac{1}{\varphi(qrs)}\sum_{\substack{m_1,m_2,m_3,n\\P^-(m_1m_2m_3n)>\sL^C\\(m_1m_2m_3n,qrs)=1}}\beta_n\bigg)\ll x\sL^{-A},\label{eq:TypeIII}
	\end{align}
	where $A,C>0$ are large constants. Although our arguments would work for complex coefficients, we assume that all the above coefficients are real-valued to simplify the presentation.
	
	Estimates for \eqref{eq:TypeI} and \eqref{eq:TypeIII} will be discussed in Sections \ref{sec:typeI} and \ref{sec:typeIII}, and in this section we focus on \eqref{eq:TypeII} under the following generic condition.
	
	\begin{condition}\label{cond:TypeII}
		Let $\eta\in~]0,10^{-2025}[.$ Suppose $MN=x,$ $x^\eta\lqs M\lqs x^{1-\eta},$ $|a_1|\lqs \sL^B$ for some $B>0$ and $|a_2|\lqs x$ with $a_1\neq a_2.$ Let $\bgamma=(\gamma_q),$ $\bdelta=(\delta_r),$ $\blambda=(\lambda_s),$ $\balpha=(\alpha_m)$ and $\bbeta=(\beta_n)$ be divisor-bounded sequences supported on positive integers with
		\begin{align*}
			q\sim Q, ~r\sim R, ~s\sim S, ~m\sim M, ~n\sim N,
		\end{align*}
		and
		\begin{align*}
			(rs,a_1)=(q,a_2)=(qr,s)=\mu^2(qr)=1.
		\end{align*}
		Assume that $\bbeta=(\beta_n)$ satisfies the Siegel--Walfisz condition $($Definition $\ref{Def:SW}).$
	\end{condition}

	Under Condition $\ref{cond:TypeII},$ we are able to estimate $\fS_{\textrm{II}}$ in three approaches.
	
	\begin{theorem}\label{thm:TypeII-1}
		Suppose $Q^2\lqs R.$ Under Condition $\ref{cond:TypeII},$ the desired bound $(\ref{eq:TypeII})$ holds provided 
		\begin{align*}
			x^\varepsilon S\lqs N\lqs x^{-\varepsilon}\min\{x^{2/3}Q^{-8/3}R^{-2/3}S^{-1/3},x^{3/2}Q^{-6}R^{-5}S^{-1},xQ^{-7/2}R^{-7/3}S^{-2/3},x^{1/2}\}.
		\end{align*}
	\end{theorem}
	\begin{theorem}\label{thm:TypeII-2}
		Suppose $Q^2\lqs S$. Under Condition $\ref{cond:TypeII},$ the desired bound $(\ref{eq:TypeII})$ holds provided 
		\begin{align*}
			x^\varepsilon S\lqs N \lqs x^{-\varepsilon}\min\{xQ^{-4}R^{-2},x^{1/3}Q^{-2}R^{-2/3}S^{2/3},x^{2/7}Q^{-10/7}R^{-2/7}S^{2/7},x^{1/3}Q^{-17/12}R^{-1/2}S^{1/3}\}.
		\end{align*}
	\end{theorem}
	\begin{theorem}\label{thm:TypeII-3}
		Under Condition $\ref{cond:TypeII},$ the desired bound $(\ref{eq:TypeII})$ holds provided $x^\varepsilon S\lqs N$ and one of the following two conditions$:$
		\begin{align*}
			&Q^4R\lqs S,~R^{2/3}S^{2/3} \lqs N\lqs x^{-\varepsilon}\min\{xR^{-1}S^{-1},x^{8/9}Q^{-22/9}R^{-17/9}S^{-5/9},x^{16/19}Q^{-40/19}R^{-34/19}S^{-10/19}\};\\
			&N\lqs x^{-\varepsilon}\min\{xR^{-1}S^{-1},xQ^{-5/2}R^{-9/4}S^{-3/4},xQ^{-9/4}R^{-2}S^{-7/8},x^{16/17}Q^{-32/17}R^{-32/17}S^{-14/17}\}.
		\end{align*}
	\end{theorem}
	
	We use Linnik's dispersion method \cite{Li63} to deal with $\fS_{\textrm{II}},$ and the arguments here are largely inspired by \cite{BFI86}. 
	
	We begin by imposing the following assumptions
	\begin{align}\label{eq:benzhi2}
		QRS\lqs x^{1-\varepsilon},\ S\lqs N^{1-\varepsilon}.
	\end{align}
	They ensure that there are enough terms in the underlying $(\ref{eq:dispersion})$ to produce considerable cancellation. To simplify the arguments, we write $\fS_\textrm{II}=\fS$ and further assume that $\beta_n=0$ whenever $n \mid a_2,$ since the contribution of such terms to \eqref{eq:TypeII} is $O(x^{1-\varepsilon})$ by \eqref{eq:benzhi2}.
	
	By Cauchy's inequality, we have (we will not display the constraint $P^-(n)>\sL^C$ for simplicity)
	\begin{align}\label{eq:dispersion}
		\fS^2&\lqs \sL^{D}MS\sum_{\substack{s,m\\(m,s)=1}}f(m)\bigg(\sum_{(qr,m)=1}\gamma_q\delta_r\bigg(\sum_{\substack{mn\equiv a_1\bmod {rs}\\ mn\equiv a_2\bmod q}}\beta_n-\frac{1}{\varphi(qrs)}\sum_{(n,qrs)=1}\beta_n\bigg)\bigg)^2
	\end{align}
	for some constant $D,$ where $f$ is smooth function supported on $[M/2,3M]$, satisfying
	$$\begin{cases}
		f(m)=1, & \text{for~}m\in [M,2M], \\
		f(m)\gqs 0, & \text{for~}m\in \bR,\\
		f^{(j)}(m)\ll_j M^{-j}, & \text{for~}\text{all}\ j\gqs 0.
	\end{cases}$$
	Squaring out we obtain
	\begin{align}\label{eq:lichayikaishi}
		\fS^2&\lqs \sL^{D}MS\cdot (\fS_1-2\fS_2+\fS_3),
	\end{align}
	where
	\begin{align*}
		&\fS_1=\sum_{\substack{s,m\\(m,s)=1}}f(m)\bigg(\sum_{(qr,m)=1}\gamma_q\delta_r\sum_{\substack{mn\equiv a_1\bmod {rs}\\ mn\equiv a_2\bmod q}}\beta_n\bigg)^2  ,\\ 
		&\fS_2=\sum_{\substack{s,m\\(m,s)=1}}f(m)\sum_{(q_1q_2r_1r_2,m)=1}\frac{\gamma_{q_1}\gamma_{q_2}\delta_{r_1}\delta_{r_2}}{\varphi(q_2r_2s)}\sum_{\substack{mn_1\equiv a_1\bmod {r_1s}\\mn_1\equiv a_2\bmod {q_1}\\(n_2,q_2r_2s)=1}}\beta_{n_1}\beta_{n_2} ,\\
		&\fS_3=\sum_{\substack{s,m\\(m,s)=1}}f(m)\bigg(\sum_{(qr,m)=1}\frac{\gamma_q\delta_r}{\varphi(qrs)}\sum_{(n,qrs)=1}\beta_n\bigg)^2 .
	\end{align*}
	
	Define the expected main term and the expected remainder term:
	\begin{align*}
		&\fS_{MT}=\widehat{f}(0)\sum_{q_i,r_i,s}\frac{\gamma_{q_1}\gamma_{q_2}\delta_{r_1}\delta_{r_2}\varphi(q_1q_2r_1r_2s)}{q_1q_2r_1r_2s\varphi(q_1r_1s)\varphi(q_2r_2s)}\sum_{(n_i,q_ir_is)=1}\beta_{n_1}\beta_{n_2} ,\\
		&\fR=xNS^{-1}\sL^{-A} .
	\end{align*}
	By $(\ref{eq:lichayikaishi})$, we expect that each $\fS_i$ $(i=1,2,3)$ can be expressed as
	\begin{align}\label{eq:expect}
		\fS_i=\fS_{MT}+O(\fR)
	\end{align}
	to obtain $(\ref{eq:TypeII})$.
	
	\subsection{Evaluation of $\fS_3$}
	\begin{lemma}[Poisson summation]
		Let $g$ be a continuous function with bounded variation satisfying that $g,\widehat{g}\in L^1(\bR)$. Then for fixed $q\in \bZ^+$ and $a\in \bZ/q\bZ$, we have
		\begin{align*}
			\sum_{m\equiv a\bmod q}g(m)=\frac{1}{q}\sum_{h\in\bZ}\widehat{g}\Big(\frac{h}{q}\Big)\ue\Big(\frac{ah}{q}\Big).
		\end{align*}
	\end{lemma}
	Integrating by parts, we find that the Fourier transform of $f$ satisfies 
	\begin{align*}
		\widehat{f}(v)\ll M(1+|v|M)^{-j}
	\end{align*}
	for any $j\gqs 0$.
	Thus, by Poisson summation, we have
	\begin{align*}
		\sum_{(m,q_1q_2r_1r_2s)=1}f(m)=\sum_{d\mid q_1q_2r_1r_2s}\frac{\mu(d)}{d}\sum_{h\in\bZ}\widehat{f}\Big(\frac{h}{d}\Big)=\frac{\varphi(q_1q_2r_1r_2s)}{q_1q_2r_1r_2s}\widehat{f}(0)+O(\tau(q_1q_2r_1r_2s)) .
	\end{align*}
	Since $N\lqs x^{1-\varepsilon}$, this yields
	\begin{align*}
		\fS_3&=\sum_{q_i,r_i,s}\frac{\gamma_{q_1}\gamma_{q_2}\delta_{r_1}\delta_{r_2}}{\varphi(q_1r_1s)\varphi(q_2r_2s)}\sum_{(n_i,q_ir_is)=1}\beta_{n_1}\beta_{n_2}\sum_{(m,q_1q_2r_1r_2s)=1}f(m)\notag\\
		&=\widehat{f}(0)\sum_{q_i,r_i,s}\frac{\gamma_{q_1}\gamma_{q_2}\delta_{r_1}\delta_{r_2}\varphi(q_1q_2r_1r_2s)}{q_1q_2r_1r_2s\varphi(q_1r_1s)\varphi(q_2r_2s)}\sum_{(n_i,q_ir_is)=1}\beta_{n_1}\beta_{n_2}+O(x^\varepsilon N^2S^{-1})\notag\\
		&=\fS_{MT}+O(\fR), 
	\end{align*}
	and so $(\ref{eq:expect})$ holds for $i=3$.
	~\\

	\subsection{Evaluation of $\fS_2$}
	
	Notice that $(q,rs)=1$ by Condition \ref{cond:TypeII}, we have
	\begin{align*}
		\fS_2=\sum_{\substack{q_i,r_i,s\\(q_i,r_is)=1}}\frac{\gamma_{q_1}\gamma_{q_2}\delta_{r_1}\delta_{r_2}}{\varphi(q_2r_2s)}\sum_{\substack{(n_i,q_ir_is)=1}}\beta_{n_1}\beta_{n_2}\sum_{\substack{mn_1\equiv a_1\bmod {r_1s}\\mn_1\equiv a_2\bmod {q_1}\\(m,q_2r_2s)=1}}f(m).
	\end{align*}
	The constraint $(m,q_2r_2s)=1$ can be relaxed by means of M\"obius inversion formula:
	\begin{align*}
		\sum_{\substack{mn_1\equiv a_1\bmod {r_1s}\\mn_1\equiv a_2\bmod {q_1}\\(m,q_2r_2s)=1}}f(m)=\sum_{\substack{v\mid q_2r_2s}}\mu(v)\sum_{\substack{vmn_1\equiv a_1\bmod {r_1s}\\vmn_1\equiv a_2\bmod {q_1}}}f(vm).
	\end{align*}
	Very similar to the second case of the coprime arguments in Section \ref{initial}, the contribution of $\fS_2$ with $v>x^\varepsilon$ can be bounded by $O(\fR)$. Since the inner sum vanishes when $(v,q_1r_1s)>1,$ Poisson summation yields
	\begin{align*}
		\sum_{\substack{mn_1\equiv a_1\bmod {r_1s}\\mn_1\equiv a_2\bmod {q_1}\\(m,q_2r_2s)=1}}f(m)&=\sum_{\substack{v\mid q_2r_2 ,v\lqs x^\varepsilon\\(v,q_1r_1s )=1}}\mu(v)\sum_{\substack{vmn_1\equiv a_1\bmod {r_1s}\\vmn_1\equiv a_2\bmod {q_1}}}f(vm)\\
		&=\sum_{\substack{v\mid q_2r_2,v\lqs x^\varepsilon\\(v,q_1r_1s)=1}}\frac{\mu(v)}{vq_1r_1s}\sum_{0\lqs |h|\lqs H_0}\widehat{f}\Big(\frac{h}{vq_1r_1s}\Big)\ue\Big(\frac{\varrho h}{q_1r_1s}\Big)+O(\fR) ,
	\end{align*}
	where $H_0=x^\varepsilon QRSM^{-1}$ and $\varrho \bmod {q_1r_1s}$ is a common solution of
	$$\begin{cases}
		v\varrho  n_1\equiv a_1\bmod {r_1s},\\  
		v\varrho  n_1\equiv a_2\bmod {q_1}.
	\end{cases}$$
	For $h=0$, we now sum up the terms and release the restriction $v\lqs x^\varepsilon$ at the cost of $O(\fR)$ to obtain the expected main term $\fS_{MT}$ , and to get
	\begin{align}\label{eq:S_2}
		\fS_2=\fS_{MT}+\sR_2+O(\fR) ,
	\end{align}
	where
	\begin{align*}
		\sR_2=\sum_{\substack{(n_i,q_ir_is)=1\\(q_1,r_1s)=1}}\frac{\gamma_{q_1}\gamma_{q_2}\delta_{r_1}\delta_{r_2}\beta_{n_1}\beta_{n_2}}{\varphi(q_2r_2s)q_1r_1s}\sum_{\substack{v\mid q_2r_2,v\lqs x^\varepsilon\\(v,q_1r_1s)=1}}\frac{\mu(v)}{v}\sum_{1\lqs|h|\lqs H_0}\widehat{f}\Big(\frac{h}{vq_1r_1s}\Big)\ue\Big(a_1h\frac{\overline{vq_1n_1}}{r_1s}+a_2h\frac{\overline{vr_1sn_1}}{q_1}\Big).
	\end{align*}
	By Bezout's identity $($Lemma \ref{lm:Bezout}$)$ we have
	\begin{align*}
		\frac{\overline{vq_1n_1}}{r_1s}\equiv -\frac{\overline{r_1s}}{vq_1n_1}+\frac{1}{vq_1r_1sn_1}\bmod 1,
	\end{align*}
	which implies that
	\begin{align*}
		\ue\Big(a_1h\frac{\overline{vq_1n_1}}{r_1s}+a_2h\frac{\overline{vr_1sn_1}}{q_1}\Big)=\ue\Big(\xi\frac{\overline{r_1s}}{vq_1n_1}\Big)+O\Big(\frac{|a_1h|}{vq_1r_1sn_1}\Big),
	\end{align*}
	where
	\begin{align*}
		\xi\equiv a_1h+a_2hvn_1\overline{vn_1}\bmod {vq_1n_1}, 
	\end{align*}
	and $vn_1\overline{vn_1}\equiv 1\bmod {q_1}$. 
	Now we use the formula
	\begin{align*}
		\frac{q_2r_2s}{\varphi(q_2r_2s)}=\frac{q_2r_2}{\varphi(q_2r_2)}\prod_{\substack{p\mid s\\ (p,q_2r_2)=1}}\Big(1-\frac{
			1}{p}\Big)^{-1}=\frac{q_2r_2}{\varphi(q_2r_2)}\sum_{\substack{p\mid \eta\Rightarrow p\mid s\\(\eta,q_2r_2)=1}}\frac{1}{\eta},
	\end{align*}
	together with partial summation and Weil's bound for Kloosterman sums $($Lemma \ref{lm:Weil}$)$, to obtain
	\begin{align*}
		\sR_2&\lqs x^{\varepsilon}\sum_{\substack{(n_i,q_ir_i)=1\\(q_i,r_i)=1}}\sum_{\substack{v\lqs x^\varepsilon \\ (v,q_1r_1)=1}}\frac{1}{vq_1q_2r_1r_2}\sum_{1\lqs |h|\lqs H_0}\bigg|\sum_{(s,vq_1q_2n_1n_2)=1}\frac{1}{s^2}\cdot\frac{q_2r_2s}{\varphi(q_2r_2s)}\widehat{f}\Big(\frac{h}{vq_1r_1s}\Big)\ue\Big(\xi\frac{\overline{r_1s}}{vq_1n_1}\Big)\bigg|\\
		&\ll x^{\varepsilon}\sum_{v,q_i,r_i,n_i}\frac{1}{vq_1q_2r_1r_2}\sum_{1\lqs h\lqs H_0}\frac{M}{S^2}\Big(\sqrt{q_1n_1}+\frac{(h,q_1n_1)S}{q_1n_1}\Big)\\
		&\ll x^{\varepsilon}(Q^{3/2}RS^{-1}N^{5/2}+RN) .
	\end{align*}
	Thus, by $(\ref{eq:S_2})$, we infer that $(\ref{eq:expect})$ holds for $i=2$ provided
	\begin{align}\label{eq:S2condition}
		R\lqs x^{-\varepsilon}\min\{xQ^{-3/2}N^{-3/2},xS^{-1}\}.
	\end{align}

	\subsection{Preliminary analysis of $\fS_1$}\label{initial}

	By expanding the square in $\fS_1$, we have
	\begin{align}\label{eq:S1}
		\fS_1=\sum_{\substack{s,m\\(m,s)=1}}f(m)\sum_{(q_1q_2r_1r_2,m)=1}\gamma_{q_1}\gamma_{q_2}\delta_{r_1}\delta_{r_2}\sum_{\substack{mn_1\equiv a_1\bmod {r_1s}\\ mn_1\equiv a_2\bmod {q_1}}}\beta_{n_1}\sum_{\substack{mn_2\equiv a_1\bmod {r_2s}\\ mn_2\equiv a_2\bmod {q_2}}}\beta_{n_2}.
	\end{align}
	
	By condition $\ref{cond:TypeII}$ and the congruences in $\fS_1$, we already have that $(q_1,r_1s)=(q_2,r_2s)=(r_1r_2,s)=(n_1,q_1r_1s)=(n_2,q_2r_2s)=1$ and $q_1,q_2,r_1,r_2$ are square-free. To make further transformations of $\fS_1$, we hope that there are additional restrictions that $(q_2r_2n_2,n_1)=1,$ $(q_1,q_2),(r_1,r_2),(q_1,r_2),$ $(q_2,r_1)\lqs\sL^C,$ $w(q_1),w(q_2),w(r_1),w(r_2)\lqs \sL^{1/3}$ and $n_1$ are square-free. To this end, we first introduce the following basic lemma.
	
	\begin{lemma}[Divisor functions on average]\label{lm:Shiu}
		Let $c>0,$ and $qx^\varepsilon<y\lqs x.$ Then for any $l\in\bZ$ we have
		\begin{align*}
			\sum_{\substack{x-y<n\lqs x\\n\equiv l\bmod q}}\tau^c(n)\ll \frac{y}{q}(\tau(q)\log x)^{2^c},
		\end{align*}
		where the constant implied in $\ll$ depends at most on $\varepsilon$ and $c.$
	\end{lemma}
	\proof Let $d=(l,q).$ Shiu's bound \cite[Theorem 2]{Sh80} deals with the case $d=1$. For $d>1$, we have
	\begin{align*}
		\sum_{\substack{x-y<n\lqs x\\n\equiv l\bmod q}}\tau^c(n)
		=\sum_{\substack{(x-y)/d<n\lqs x/d\\n\equiv l/d\bmod {q/d}}}\tau^c(dn)
		\leqslant \tau^c(d)\sum_{\substack{(x-y)/d<n\lqs x/d\\n\equiv l/d\bmod {q/d}}}\tau^c(n).
	\end{align*}
	Since $(l/d,q/d)=1,$ the desired bound now follows directly from \cite[Theorem 2]{Sh80}. 
	
	See also \cite[Lemma 12]{BFI86} and \cite[Lemma 8.7]{Ma25a} for alternative proofs.
	\endproof
	
	
	We are now ready to work out the above-mentioned conditions. In fact, the other terms of $\fS_1$ are negligible due to Condition $\ref{cond:TypeII}$ and Lemma \ref{lm:Shiu}. We will argue this point in four cases$:$
	
	\begin{enumerate}
		\item We may adopt the reasonable assumptions that $(r_1r_2,s)=\mu^2(r_1)=\mu^2(r_2)=1$ since $r,s$ come from the well-factorable weight $\blambda=(\lambda_d)$ supported on square-free numbers.
		
		\item Let $C$ be sufficiently large. Now we wish to prove, for example, that the terms with $n_0:=(n_1,n_2)>\sL^C$ contribute $O(\fR)$. By the triangle inequality, we can deduce that the contribution is bounded by
		\begin{align*}
			&\sum_{n_0>\sL^C}\sum_{s,m}\sum_{q_1,r_1}|\gamma_{q_1}||\delta_{r_1}|\sum_{\substack{n_0n_1\equiv a_1\overline{m}\bmod {r_1s}\\n_0n_1\equiv a_2\overline{m}\bmod {q_1}}}|\beta_{n_0n_1}|\sum_{q_2,r_2}|\gamma_{q_2}||\delta_{r_2}|\sum_{\substack{n_0n_2\equiv a_1\overline{m}\bmod{r_2s}\\n_0n_2\equiv a_2\overline{m}\bmod {q_2}}}|\beta_{n_0n_2}|.
		\end{align*}
		Now by Lemma $\ref{lm:Shiu}$, H\"older's inequality and assumption $(\ref{eq:benzhi2}),$ we have (recall that $mn\neq a_1,a_2$ since $|a_1|\lqs \log^Bx$ and $n\nmid a_2$)
		\begin{align*}
			\sum_{q_2,r_2}|\gamma_{q_2}||\delta_{r_2}|\sum_{\substack{n_0n_2\equiv a_1\overline{m}\bmod {r_2s}\\n_0n_2\equiv a_2\overline{m}\bmod {q_2}}}|\beta_{n_0n_2}|&\ll \sum_{\substack{x<n\lqs 4x,~n\neq a_1,a_2\\n\equiv a_1\bmod {s}\\n\equiv 0\bmod {mn_0}}}\tau^A(|n-a_1|)\tau^A(|n-a_2|)\tau^A(n)\\
			&\ll N(n_0s)^{-1}(\tau(mn_0s)\sL)^{O_A(1)}+x^{\varepsilon/2}\\
			&\ll NS^{-1}(\tau(mn_0s))^{O_A(1)}\sL^{-C/2}
		\end{align*}
		for some constant $A.$ Again, by H\"older's inequality and assumption $(\ref{eq:benzhi2}),$ we see that the contribution of the terms with $n_0:=(n_1,n_2)>\sL^C$ can be bounded by
		\begin{align*}
			&\ \ \ \  NS^{-1}\sL^{-C/2}\sum_{n_0,q_1,r_1,s}|\gamma_{q_1}||\delta_{r_1}|\sum_{\substack{mn_0n_1\equiv a_1\bmod {r_1s}\\mn_0n_1\equiv a_2\bmod {q_1}}}|\beta_{n_0n_1}|\tau^{O_{A}(1)}(mn_0s)\\
			&\ll NS^{-1}\sL^{-C/2}\sum_{\substack{x<n\lqs 4x\\n\neq a_1,a_2}}\tau(|n-a_1|)^{O_A(1)}\tau(|n-a_2|)^{O_A(1)}\tau(n)^{O_A(1)}\\
			&\ll xNS^{-1}\sL^{-C/4}\ll \fR.
		\end{align*}
		Using the same method, we see that the contributions of the terms with $(q_1,q_2)>\sL^C,$ $(r_1,r_2)>\sL^C,$ $(q_1,r_2)>\sL^C,$ $(q_2,r_1)>\sL^C,$ $(q_2,n_1)>\sL^C,$ $(r_2,n_1)>\sL^C$ are bounded by $O(\fR).$ Recalling that $P^-(n)>\sL^C$ by Condition \ref{cond:TypeII}, we may restrict our attention to the case $(q_2r_2n_2,n_1)=1,$ incurring only a negligible error term.
		
		\item Similarly, by Lemma $\ref{lm:Shiu},$ H\"older's inequality and assumption $(\ref{eq:benzhi2}),$ the sums over $q_1$, $r_1$ and $n_1$ with $\mu^2(n_1)=0$ can be bounded by $($recall that $P^-(n)>\sL^C)$
		\begin{align*}
			\sum_{p_0>\sL^C}\sum_{\substack{x<n\lqs 4x,~n\neq a_1,a_2\\n\equiv a_1\bmod{s}\\n\equiv 0\bmod{mp_0^2}}}\tau^A(|n-a_1|)\tau^A(|n-a_2|)\tau^A(np^2_0)\ll NS^{-1}\tau^{O_A(1)}(ms)\sL^{-C/2}
		\end{align*}
		for some constants $A.$ Again, by H\"older's inequality and assumption $(\ref{eq:benzhi2})$, the above expression contributes in total
		\begin{align*}
			NS^{-1}\sL^{-C/2}\sum_{q_2,r_2,s}\tau^{O_A(1)}(q_2r_2s)\sum_{\substack{mn_2\neq a_1,a_2\\mn_2\equiv a_1\bmod {r_2s}\\mn_2\equiv a_2\bmod {q_2}}}\tau^{O_{A}(1)}(mn_2)=O(\fR).
		\end{align*}
		Therefore, we may restrict our attention to the case $\mu^2(n_1)=1,$ incurring only a negligible error term.
		
		\item The contribution of terms with $\omega(r_2)> \sL^{1/3}$ can be bounded by
		\begin{align*}
			\qquad\qquad&\exp(\sL^{-1/3})\sum_{q_1,r_1,s}|\gamma_{q_1}\delta_{r_1}|\sum_{\substack{mn_1\equiv a_1\bmod {r_1s}\\ mn_1\equiv a_2\bmod {q_1}}}|\beta_{n_1}|\sum_{\substack{mn_2\equiv a_1\bmod {s}}}|\beta_{n_2}|\sum_{\substack{q_2\mid mn_2-a_2\\r_2\mid mn_2-a_1}}|\gamma_{q_2}\delta_{r_2}|\exp(\omega(r_2))\\
			\ll &\exp(\sL^{-1/3})\sum_{q_1,r_1,s}\tau^A(q_1r_1)\sum_{\substack{mn_1\neq a_1,a_2\\mn_1\equiv a_1\bmod {r_1s}\\ mn_1\equiv a_2\bmod {q_1}}}\tau^A(n_1)\sum_{\substack{x<n\lqs 4x,~n\neq a_1,a_2\\n\equiv a_1\bmod {s}\\n\equiv 0\bmod m}}\tau^A(n)\tau^A(|n-a_1|)\tau^A(|n-a_2|).
		\end{align*}
		Similarly, we obtain that the above expression is $O(\fR).$
		By the same method, we may impose the constraints $\omega(q_1),\omega(q_2),\omega(r_1),\omega(r_2)\lqs\sL^{1/3}$
		at the cost of a negligible error term.
	\end{enumerate}
	
	Therefore, based on the preceding arguments, we can incorporate the coprime and square-free conditions with an admissible error term, namely,
	\begin{align}\label{eq:S1*}
		\fS_1=\fS_1^*+O(\fR),
	\end{align}
	where $\fS_1^*$ denotes the summations in $\fS_1$ with the restrictions that 
	\begin{align*}
		(q_1,r_1s)=(q_2,r_2s)=(r_1r_2,s)=(q_2r_2n_2,n_1)=1,  
	\end{align*}
	\begin{align*}
		\max\{(q_1,q_2),(r_1,r_2),(q_1,r_2),(q_2,r_1)\}\lqs\sL^C,\ \ \max\{\omega(q_1),\omega(q_2),\omega(r_1),\omega(r_2)\}\lqs\sL^{1/3},  
	\end{align*}
	and that $q_1,q_2,r_1,r_2,n_1$ are square-free. Thanks to these conditions, we can later perform convenient transformations, such as those relying on the Chinese remainder theorem. In the following arguments, we do not display these restrictions explicitly to simplify the presentation. 
	
	Note that
	\begin{align*}
		\fS_1^*&=\sum_{q_0,r_0,t_1,t_2\lqs\sL^C}\sum_{\substack{q_1,q_2,r_1,r_2,s\\(q_1,q_2)=q_0,(r_1,r_2)=r_0\\(r_1,q_2)=t_1,(r_2,q_1)=t_2}}\gamma_{q_1}\gamma_{q_2}\delta_{r_1}\delta_{r_2}\sum_{\substack{n_1\equiv n_2\bmod {q_0r_0s}\\a_2n_1\equiv a_1n_2\bmod{t_1}\\a_1n_1\equiv a_2n_2\bmod{t_2}}}\beta_{n_1}\beta_{n_2}\\
		&\ \ \ \ \times\sum_{m\equiv \varrho\bmod {q_0q_1q_2r_0r_1r_2t_1t_2s}}f(m)+O(\fR),
	\end{align*}
	where $\varrho\bmod {q_0q_1q_2r_0r_1r_2t_1t_2s}$ is a common solution of
	$$\begin{cases}
		\varrho n_1\equiv a_1\bmod {r_0r_1t_1s},\\
		\varrho n_2\equiv a_1\bmod {r_0r_2t_2s},\\
		\varrho n_1\equiv a_2\bmod {q_0q_1t_2},\\
		\varrho n_2\equiv a_2\bmod {q_0q_2t_1}.
	\end{cases}$$
	From Poisson summation it follows that
	\begin{equation*}
		\begin{split}\fS_1^*
			&=\sum_{q_0,r_0,t_1,t_2\lqs\sL^C}\sum_{\substack{q_1,q_2,r_1,r_2,s}}\frac{\gamma_{q_0q_1t_2}\gamma_{q_0q_2t_1}\delta_{r_0r_1t_1}\delta_{r_0r_2t_2}}{q_0q_1q_2r_0r_1r_2t_1t_2s}\sum_{\substack{n_1\equiv n_2\bmod {q_0r_0s}\\a_2n_1\equiv a_1n_2\bmod{t_1}\\a_1n_1\equiv a_2n_2\bmod{t_2}}}\beta_{n_1}\beta_{n_2}\\
			&\quad\ \times\sum_{0\lqs |h|\lqs H}\widehat{f}\Big(\frac{h}{q_0q_1q_2r_0r_1r_2t_1t_2s}\Big)\ue\Big(\frac{\varrho h}{q_0q_1q_2r_0r_1r_2t_1t_2s}\Big)+O(\fR),
		\end{split}
	\end{equation*}
	where $H=x^\varepsilon Q^2R^2SM^{-1}.$ 
	
	For the terms with $h=0$, we handle the congruence of $n_1,n_2$ by applying the following result, which is \cite[Theorem 9.16]{FI10}.
	\begin{lemma}[Generalized Bombieri--Vinogradov theorem]\label{lm:lemma3}
		Let $MN=x,$ $x^\varepsilon <M< x^{1-\varepsilon}$ and fix $A>0$. Let $(\alpha_m),(\beta_n)$ be divisor-bounded sequences supported on $[M,2M]$ and $[N,2N],$ respectively$,$ with $(\beta_n)$ satisfying the Siegel--Walfisz condition $($Definition $\ref{Def:SW})$. Then for some sufficiently large $B,$ any $b,c\in\bZ^+$ and $Q\lqs x^{1/2}\sL^{-B}$ we have 
		\begin{align*}
			\sum_{q\lqs Q}\max_{(a,q)=1}\Bigg|\sum_{\substack{m\overline{n}\equiv a\bmod q\\(m,bq)=(n,cq)=1}}\alpha_m\beta_n-\frac{1}{\varphi(q)}\sum_{\substack{(m,bq)=(n,cq)=1}}\alpha_m\beta_n\Bigg|\ll x\log^{-A}x .
		\end{align*} 
	\end{lemma}
	
	Let
	\begin{align*}
		\varrho_1\equiv t_1t_2\overline{t_1t_2}+a_1\overline{a_2}q_0r_0t_2s\overline{q_0r_0t_2s}+\overline{a_1}a_2q_0r_0t_1s\overline{q_0r_0t_1s}\bmod{q_0r_0t_1t_2s}.
	\end{align*}
	By applying Lemma $\ref{lm:lemma3}$ and extending the summation to all $n_0,$ the contribution of the terms with $h=0$ is given by (there are implicit restrictions $(n_1,q_0q_1r_0r_1t_1)=(n_2,q_0q_2r_0r_2t_2)=1$ in the summation over $n_1,n_2$)
	\begin{align*}
		&\quad\  \widehat{f}(0)\sum_{\substack{q_0,r_0,t_1,t_2\lqs\sL^C\\q_1,q_2,r_1,r_2,s}}\frac{\gamma_{q_0q_1t_2}\gamma_{q_0q_2t_1}\delta_{r_0r_1t_1}\delta_{r_0r_2t_2}}{q_0q_1q_2r_0r_1r_2t_1t_2s}\sum_{n_1\overline{n_2}\equiv \varrho_1\bmod {q_0r_0t_1t_2s}}\beta_{n_1}\beta_{n_2}\\
		&=\widehat{f}(0)\sum_{\substack{q_0,r_0,t_1,t_2\lqs\sL^C\\q_1,q_2,r_1,r_2,s}}\frac{\gamma_{q_0q_1t_2}\gamma_{q_0q_2t_1}\delta_{r_0r_1t_1}\delta_{r_0r_2t_2}}{q_0q_1q_2r_0r_1r_2t_1t_2s\varphi(q_0r_0t_1t_2s)}\sum_{(n_1n_2,q_0r_0t_1t_2s)=1}\beta_{n_1}\beta_{n_2}+O(\fR) 
	\end{align*}
	provided that $(\ref{eq:benzhi2})$ holds. Extending the summation to all $q_0,r_0,t_1,t_2$ shows that the above expression equals $\fS_{MT}+O(\fR).$
	Therefore, under the assumption $(\ref{eq:benzhi2})$, we deduce that
	\begin{align}\label{eq:S_1}
		\fS_1^*=\fS_{MT}+\sR^*_1+O(\fR) ,
	\end{align}
	where
	\begin{align*}
		\sR^*_1=&\sum_{\substack{q_0,r_0,t_1,t_2\lqs\sL^{C}\\q_1,q_2,r_1,r_2,s}}\frac{\gamma_{q_0q_1t_2}\gamma_{q_0q_2t_1}\delta_{r_0r_1t_1}\delta_{r_0r_2t_2}}{q_0q_1q_2r_0r_1r_2t_1t_2s}\sum_{\substack{n_1\equiv n_2\bmod {q_0r_0s}\\a_2n_1\equiv a_1n_2\bmod{t_1}\\a_1n_1\equiv a_2n_2\bmod{t_2}}}\beta_{n_1}\beta_{n_2}\\
		&\times\sum_{1\lqs |h|\lqs H}\widehat{f}\Big(\frac{h}{q_0q_1q_2r_0r_1r_2t_1t_2s}\Big)\Big(\frac{\varrho h}{q_0q_1q_2r_0r_1r_2t_1t_2s}\Big).
	\end{align*}
	
	To estimate $\sR^*_1$, we focus on the variables $q_1,q_2,r_1,r_2,s,n_1,n_2,h$. Since we aim to bound $\sR^*_1$ with a saving of arbitrarily large power of $\sL$, it should be harmless if summing over $q_0,r_0,t_1,t_2\lqs\sL^{C}$ trivially.
	Hence we may assume $q_0=r_0=t_1=t_2=1$ to simplify the expressions. By Lemma $\ref{lm:change moduli}$ and \eqref{eq:easechange}, the problem is therefore reduced to estimating
	\begin{align*}
		\sR_1=&\sum_{q_1,q_2,r_1,r_2,s}\frac{\gamma_{q_1}\gamma_{q_2}\delta_{r_1}\delta_{r_2}}{q_1q_2r_1r_2s}\sum_{n_1\equiv n_2\bmod {s}}\beta_{n_1}\beta_{n_2}\sum_{1\lqs |h|\lqs H}\widehat{f}\Big(\frac{h}{q_1q_2r_1r_2s}\Big)\\
		&\times \ue\Big (a_1h\frac{n_1-n_2}{s}\frac{\overline{q_1r_1n_2}}{q_2r_2n_1}+(a_2-a_1)h\Big(\frac{\overline{q_2r_1r_2sn_1}}{q_1}+\frac{\overline{q_1r_1r_2sn_2}}{q_2}\Big)\Big).
	\end{align*}
	In the next three sections, we use three methods to prove
	\begin{align}\label{eq:sR1}
		\sR_1=O(\fR)
	\end{align}
	when the lengths of summations are suitably balanced.
	This result, combined with $(\ref{eq:S1*})$ and $(\ref{eq:S_1})$, implies that $(\ref{eq:expect})$ holds for $i=1$. Consequently, the specific constraints required by each method, together with assumptions $(\ref{eq:S2condition})$, $(\ref{eq:benzhi2})$, imply the corresponding conditions in Theorems $\ref{thm:TypeII-1},$ $\ref{thm:TypeII-2},$ and $\ref{thm:TypeII-3},$ respectively.
	
	\begin{remark}
		Unlike the situation in \cite{FI83} and \cite{BFI86} $(q_1=q_2=1),$ there is one more congruence condition $mn\equiv a_2\bmod q$ in $(\ref{eq:TypeII})$, which makes the variables entangled in a very complicated way $($as one may see from Lemma $\ref{lm:change moduli})$. As a result, we need to take more care in separating the variables in order to apply Kloostermania, which produces an additional loss from the residue classes modulo $q_1q_2$.
	\end{remark}

	\section{Estimation of $\sR_1$: the first method}\label{firstmethod}
	\subsection{Preparation for Kloostermania}\label{Fou1} In this method we plan to use Cauchy's inequality in the form
	\begin{align*}
		\sum\ldots\sum\Big|\sum_{r_1}\sum_{r_2}\Big|,
	\end{align*}
	and then apply Kloostermania.
	
	We expand $\widehat{f}$ by its definition, and divide the sums over $r_1,r_2$ according to residue classes modulo $q_1q_2$, to get
	\begin{align*}
		\sR_1&\ll  x^\varepsilon\sum_{\substack{q_1,q_2\\1\lqs |h|\lqs H}}\sum_{\substack{a_3\bmod {q_1q_2}\\a_4\bmod {q_1q_2}}}\sum_{\substack{k,s,n_1,n_2\\n_1-n_2=ks}}\int \frac{f(wq_1q_2s/h)}{|h|}\Big|\sum_{\substack{r_1\equiv a_3\bmod {q_1q_2}\\r_2\equiv a_4\bmod {q_1q_2}}}\frac{\delta_{r_1}\delta_{r_2}}{r_1r_2}\ue\Big(\frac{-w}{r_1r_2}\Big)\ue\Big(a_1hk\frac{\overline{q_1r_1n_2}}{q_2r_2n_1}\Big)\Big|\ud w,
	\end{align*}
	where $k$ is restricted to $[-K,K]\setminus\{0\}$ $($notice that $k\neq 0$ since $(n_1,n_2)=1)$ with $K=N/S$. Originally, the restriction is $n_1\equiv n_2\bmod {s}$. Now we can switch the roles of $s$ and $k$ as $n_1\equiv n_2\bmod {k}$. In order to apply Kloostermania $($Lemma  $\ref{lm:kloostermania})$, we attach a smooth weight function $\Phi(\cdot,\cdot)$ which majorizes the indicator function of the box $[1,2]^2$ to the sums over $n_1,n_2$. Therefore, we deduce that 
	\begin{align}\label{eq:trick1}
		\sR_1&\ll  \frac{x^\varepsilon M}{Q^2R^2S}\sum_{q_1,q_2}\sum_{\substack{1\lqs |h|\lqs H\\1\lqs|k|\lqs K}}\sum_{\substack{a_3\bmod {q_1q_2}\\a_4\bmod {q_1q_2}}}\sum_{n_1\equiv n_2\bmod {k}}\Phi\Big(\frac{n_1}{N},\frac{n_2}{N}\Big)\Big|\sum_{\substack{r_1\equiv a_3\bmod {q_1q_2}\\r_2\equiv a_4\bmod {q_1q_2}}}\eta_1(r_1,r_2)\ue\Big(a_1hk\frac{\overline{q_1r_1n_2}}{q_2r_2n_1}\Big)\Big|
	\end{align}
	with some $\Vert\eta_1\Vert_\infty\lqs 1$.
	Now we apply Cauchy's inequality and then remove the congruence condition $n_1\equiv n_2\bmod {k}$, to obtain 
	\begin{align}\label{eq:Klresult1}
		\sR^2_1	\ll \frac{Q^2Hx^{2+\varepsilon}}{R^{4}S^2}|\sN_1|,
	\end{align}
	where
	\begin{align*}
		\sN_1=\sum_{q_1,q_2}\sum_{\substack{1\lqs |h|\lqs H\\1\lqs|k|\lqs K}}\sum_{\substack{r_1\equiv \widetilde{r_1}\bmod {q_1q_2}\\r_2\equiv \widetilde{r_2}\bmod {q_1q_2}}}\eta_1(r_1,r_2)\overline{\eta_1}(\widetilde{r_1},\widetilde{r_2})\sum_{n_1}\sum_{n_2}\Phi\Big(\frac{n_1}{N},\frac{n_2}{N}\Big)\ue\Big( a_1hk\Big(\frac{\overline{q_1r_1n_2}}{q_2r_2n_1}-\frac{\overline{q_1\widetilde{r_1}n_2}}{q_2\widetilde{r_2}n_1}\Big)\Big).
	\end{align*}

	In addition to the conditions in Section \ref{initial}, the subsequent transformations would require further coprime constraints $(r_1,\widetilde{r_2})=(r_2,\widetilde{r_1})=1.$ To this end, we appeal to an elegant trick due to Fouvry \cite[Lemma 6]{Fo84a}.

	\begin{lemma}[Splitting into coprime sets]\label{lm:splitting into coprime sets}
		For $\nu\in\bZ^+$ and $Q,R\gqs 2,$ we define
		\begin{align*}
			\sG_\nu(Q,R):=\{(q,r)\in[1,Q]\times[1,R]:(q,r)=1\ and\ \omega(q),\omega(r)\lqs \nu\}.
		\end{align*} 	
		Then there exists a partition of $\sG_\nu(Q,R)$ into at most $(2\log(3QR))^{\nu^2}$ subsets $\sG^*_\nu(Q,R)$ having the following property$:$
		\begin{align*}
			(q_1,r_1),(q_2,r_2)\in \sG^*_\nu(Q,R)\Rightarrow (q_1,r_2)=(q_2,r_1)=1.
		\end{align*}
	\end{lemma}
	
	Recall that, according to the preceding arguments, the number of prime factors of $r_1,r_2$ is at most $\sL^{1/3}$. Consequently in \eqref{eq:trick1}, we apply Fouvry's trick (Lemma \ref{lm:splitting into coprime sets}) to partition the summations over $r_1,r_2$ into at most $\exp(\sL ^{2/3+\varepsilon})$ sub-sums. Then we apply Cauchy's inequality to each sub-sum to obtain \eqref{eq:Klresult1}, with the constraints $(r_1,\widetilde{r_2})=(r_2,\widetilde{r_1})=1.$ This partitioning is efficient since the subsequent saving in the exponential sum estimates is at least $x^\varepsilon$, which dominates the cost $\exp(\sL^{2/3+\varepsilon})$. It therefore suffices to estimate a single sub-sum.
	To simplify the presentation, we still use the notation $\sN_{1}$, with the restrictions $(r_1,\widetilde{r_2})=(r_2,\widetilde{r_1})=1$ now in force. Combining these with the previously established coprime conditions from Section \ref{initial}, we proceed to simplify the fraction in $\sN_{1}$:
	
	\begin{align*}
		\sN_1=\sum_{q_1,q_2}\sum_{\substack{1\lqs |h|\lqs H\\1\lqs|k|\lqs K}}\sum_{\substack{r_1\equiv \widetilde{r_1}\bmod {q_1q_2}\\r_2\equiv \widetilde{r_2}\bmod {q_1q_2}}}\eta_1(r_1,r_2)\overline{\eta_1}(\widetilde{r_1},\widetilde{r_2})\sum_{n_1}\sum_{n_2}\Phi\Big(\frac{n_1}{N},\frac{n_2}{N}\Big)\ue\Big(a_1hk(\widetilde{r_1}\widetilde{r_2}-r_1r_2)\frac{\overline{q_1r_1\widetilde{r_1}n_2}}{q_2r_2\widetilde{r_2}n_1}\Big).
	\end{align*}

	\subsection{Application of Kloostermania}
	For $QR^2\lqs u,v\lqs 8QR^2$ and $0\lqs w\lqs 3R^2HK\sL^B$, we put
	\begin{align*}
		\cB_1(u,v,w)=\sum_{\substack{q_1r_1\widetilde{r_1}=u,q_2r_2\widetilde{r_2}=v\\r_1\equiv \widetilde{r_1}\bmod {q_1q_2},r_2\equiv \widetilde{r_2}\bmod {q_1q_2}}}\sum_{|a_1hk(r_1r_2-\widetilde{r_1}\widetilde{r_2})|=w}\eta_1(r_1,r_2)\overline{\eta_1}(\widetilde{r_1},\widetilde{r_2}).
	\end{align*}
	Then  we can rewrite $\sN_1$ as
	\begin{align*}
		\sN_1=\sum_{\substack{QR^2\lqs u\lqs 8QR^2\\QR^2\lqs v\lqs 8QR^2}}\sum_{0\lqs w\lqs 3R^2HK\sL^B}\cB_1(u,v,w)\sum_{n_1}\sum_{n_2}\Phi\Big(\frac{n_1}{N},\frac{n_2}{N}\Big)\ue\Big(w\frac{\overline{un_2}}{vn_1}\Big),
	\end{align*}
	which is of type $\sK(4QR^2,4QR^2,3R^2HK\sL^B,N,N)$ in Lemma $\ref{lm:kloostermania}$.
	
	The diagonal terms with $w=0$ contribute
	\begin{align}\label{eq:diagonal1}
		\sN_1{(w=0)}\ll x^\varepsilon Q^2R^2N^2HK.
	\end{align}
	
	For the off-diagonal terms with $w\neq 0$, we apply Cauchy's inequality to deduce that
	\begin{align*}
		\sum_{u,v,w}|\cB_1(u,v,w)|^2&\ll (QRK)^\varepsilon \sum_{k,q_1,q_2}\sum_{\substack{r_1\equiv \widetilde{r_1}\bmod {q_1q_2}\\r_2\equiv \widetilde{r_2}\bmod {q_1q_2}}}\sum_{w> 0}\Big(\sum_{\substack{1\lqs |h|\lqs H\\|a_1hk(r_1r_2-\widetilde{r_1}\widetilde{r_2})|=w\\}}1\Big)^2\\
		&\ll x^\varepsilon Q^{-2}R^4HK 
	\end{align*}
	assuming
	\begin{align}\label{eq:RQ2}
		R\gqs Q^2.
	\end{align}
	Thus, by Lemma $\ref{lm:kloostermania}$ we have 
	\begin{align}\label{eq:kloostermania1}
		\sN_1^2(w\neq 0)&\ll x^\varepsilon \big(Q^4R^8N^2+Q^{5/2}R^5N^3\big)Q^{-2}R^4HK\notag\\
		&\ll x^\varepsilon(Q^4R^{14}M^{-1}N^3+Q^{5/2}R^{11}M^{-1}N^4)
	\end{align}
	assuming that
	\begin{align}\label{eq:assume1}
		N^2\lqs x^{1-\varepsilon}.
	\end{align}
	
	Combining $(\ref{eq:Klresult1})$, $(\ref{eq:diagonal1})$, $(\ref{eq:kloostermania1})$ we have (recall that $MN=x$)
	\begin{align*}
		\sR_1\ll &x^{\varepsilon}(Q^{4}RN^{5/2}S^{-1/2}+Q^{3}R^{5/2}S^{-1/2}N^{3/2}x^{1/4}+Q^{21/8}R^{7/4}S^{-1/2}N^{7/4}x^{1/4}).
	\end{align*}
	To guarantee $(\ref{eq:sR1})$, we require
	\begin{align*}
		N\lqs x^{-\varepsilon}\min\{x^{2/3}Q^{-8/3}R^{-2/3}S^{-1/3},x^{3/2}Q^{-6}R^{-5}S^{-1},xQ^{-7/2}R^{-7/3}S^{-2/3}\}.
	\end{align*}
	This requirement, together with assumptions $(\ref{eq:S2condition})$, $(\ref{eq:benzhi2})$, $(\ref{eq:RQ2})$, $(\ref{eq:assume1})$, establishes Theorem \ref{thm:TypeII-1}.

	\section{Estimation of $\sR_1$: the second method}\label{secondmethod}
	\subsection{Preparation for Kloostermania}
	In this method we plan to use Cauchy's inequality in the form
	\begin{align*}
		\sum\ldots\sum\Big|\sum_{n_2}\sum_{h}\Big|,
	\end{align*}
	and then apply Kloostermania.
	
	To separate the variables, we require more effort than in Section \ref{firstmethod} (the first method). We first appeal to the M\"obius formula to get rid of the implicit restriction $(s,a_1r_1r_2)=1,$ and to get
	\begin{align*}
		\sR_1=&\sum_{q_1,q_2,r_1,r_2}\frac{\gamma_{q_1}\gamma_{q_2}\delta_{r_1}\delta_{r_2}}{q_1q_2r_1r_2}\sum_{\delta\mid a_1r_1r_2}\frac{\mu(\delta)}{\delta}\sum_{s\sim S/\delta}\frac{1}{s}\sum_{n_1\equiv n_2\bmod {\delta s}}\beta_{n_1}\beta_{n_2}\sum_{1\lqs |h|\lqs H}\widehat{f}\Big(\frac{h}{q_1q_2r_1r_2\delta s}\Big)\\
		&\times \ue\Big (a_1h\frac{n_1-n_2}{\delta s}\frac{\overline{q_1r_1n_2}}{q_2r_2n_1}+(a_2-a_1)h\Big(\frac{\overline{q_2r_1r_2\delta sn_1}}{q_1}+\frac{\overline{q_1r_1r_2\delta sn_2}}{q_2}\Big)\Big).
	\end{align*}
	We now switch the roles of $s$ and $k$ in the relation $n_1-n_2=k\delta s$, and obtain the constraints (recall that $(n_1,n_2)=1$)
	\begin{align}\label{eq:condition1}
		(n_1n_2,k\delta )=1,\ n_1\equiv n_2\bmod {k\delta }
	\end{align}
	and
	\begin{align}\label{eq:condition2}
		|n_1-n_2|\sim |k|S.
	\end{align}
	In particular, we have 
	\begin{align*}
		1\lqs |k| \lqs K=N/S.
	\end{align*}
	First, we detect the condition $(\ref{eq:condition1})$ by means of multiplicative characters $\chi\bmod k$ and $\psi\bmod {\delta}$:
	
	$$\frac{1}{\varphi(k)}\sum_{\chi\bmod k}\overline{\chi}(n_1)\chi(n_2)=\begin{cases}
		1,& \text{if~}(n_1n_2,k)=1\text{~and~}n_1\equiv n_2\bmod k,\\
		0,& \text{otherwise},
	\end{cases}$$
	$$\frac{1}{\varphi(\delta)}\sum_{\psi\bmod {\delta}}\overline{\psi}(n_1)\psi(n_2)=\begin{cases}
		1,& \text{if~}(n_1n_2,\delta)=1\text{~and~}n_1\equiv n_2\bmod {\delta},\\
		0,& \text{otherwise}.
	\end{cases}$$
	Next, we detect the condition $(\ref{eq:condition2})$ by means of additive characters:
	$$\int_{0}^{1}\ue((n_1-n_2)\alpha)F(\alpha)\ud\alpha=\begin{cases}
		1,& \text{if~}|n_1-n_2|\sim |k|S,\\
		0,& \text{otherwise},
	\end{cases}$$
	where
	\begin{align}
		F(\alpha):=\sum_{|l|\sim |k|S}\ue(l\alpha).
	\end{align}
	Notice that
	\begin{align*}
		\int_{0}^{1}|F(\alpha)|\ud\alpha\ll\log 2N.
	\end{align*}
	Finally, we separate the variables in $\widehat{f}$ by 
	\begin{align*}
		\frac{1}{q_1q_2r_1r_2\delta s}\widehat{f}\Big(\frac{h}{q_1q_2r_1r_2\delta s}\Big)
		&=\int_{0}^{\frac{8M}{Q^2R^2S}}f(q_1q_2r_1r_2\delta su)\ue(-hu)\ud u\\
		&=\int_{0}^{\frac{8M}{Q^2R^2S}}\ue(-hu)\ud u\int_{\bR}\widehat{f}(v)\ue(q_1q_2r_1r_2\delta suv)\ud v\\
		&=\frac{k}{q_1q_2r_1r_2}\int_{0}^{\frac{8M}{Q^2R^2S}}\ue(-hu)\ud u\int_{\bR}\widehat{f}\Big(\frac{kv}{q_1q_2r_1r_2}\Big)\ue((n_1-n_2)uv)\ud v,
	\end{align*}
	and notice that
	\begin{align*}
		\int_{\bR}\Big|\widehat{f}\Big(\frac{kv}{q_1q_2r_1r_2}\Big)\Big|\ud v\ll \frac{q_1q_2r_1r_2}{|k|}.
	\end{align*}
	
	The argument above separates the variables $n_1$ and $n_2$. Based on the residue classes $n_2\Mod{kq_1q_2}$ and $h\Mod{q_1q_2}$, we have
	\begin{align*}\label{eq:Klresult12}
		\sR_1&=\sum_{\substack{q_1,q_2,r_1,r_2,n_1\\1\lqs |k|\lqs K}}\frac{\gamma_{q_1}\gamma_{q_2}\delta_{r_1}\delta_{r_2}\beta_{n_1}}{q_1q_2r_1r_2}\sum_{\substack{a_5\bmod {kq_1q_2}\\a_6\bmod {q_1q_2}}}\ue\Big((a_2-a_1)a_6\Big(\frac{\overline{q_2r_1r_2n_1(n_1-a_5)/k}}{q_1}+\frac{\overline{q_1r_1r_2a_5(n_1-a_5)/k}}{q_2}\Big)\Big) \\
		&\ \ \ \qquad\times\sum_{\delta\mid a_1r_1r_2}\frac{k}{\varphi(k)\varphi(\delta)}  \sum_{\substack{\chi\bmod k\\ \psi\bmod {\delta}}}\overline{\chi\psi}(n_1)\int_{0}^{\frac{M}{Q^2R^2S}}\ud u\int_{\bR}\widehat{f}\Big(\frac{kv}{q_1q_2r_1r_2}\Big)\ud v\int_{0}^{1}F(\alpha
		)\ue(n_1(uv+\alpha))\\
		&\ \ \ \qquad\times\sum_{\substack{n_2\equiv a_5\bmod {kq_1q_2}}}\sum_{\substack{1\lqs |h|\lqs H\\h\equiv a_6(q_1q_2)}}\beta_{n_2}\ue(-n_2(uv+\alpha))\ue(-hu)\chi\psi(n_2)\ue\Big(a_1hk\frac{\overline{q_1r_1n_2}}{q_2r_2n_1}\Big) \ud\alpha\\
		&\ll  \frac{Mx^\varepsilon}{Q^2R^2S}\sum_{\substack{q_1,q_2,r_1,r_2\\}}\sum_{k,n_1}\frac{1}{\varphi(k)}\sum_{\chi\bmod k}\sum_{\substack{a_5\bmod {kq_1q_2}\\a_6\bmod {q_1q_2}}}\Big|\sum_{\substack{1\lqs |h|\lqs H\\h\equiv a_6\bmod {q_1q_2}\\n_2\equiv a_5\bmod {kq_1q_2}}}\eta_2(h,n_2)\chi(n_2)\ue\Big(a_1hk\frac{\overline{q_1r_1n_2}}{q_2r_2n_1}\Big)\Big|
	\end{align*}
	with some $|\eta_2(h,n_2)|\lqs 1$. Now by Cauchy's inequality and combining the variables as $r_2n_1=c$, we obtain
	\begin{align}
		\sR_1^2\ll \frac{Q^2K^2Mx^{1+\varepsilon}}{R^2S^2}|\sN_2|,
	\end{align}
	where
	\begin{align*}
		\sN_2=\sum_{k,q_1,q_2}\sum_{\substack{n_2\equiv n_2'\bmod {kq_1q_2}\\h\equiv h'\bmod {q_1q_2}}}\eta_2(h,n_2)\overline{\eta_2}(h',n_2')\sum_{c}\sum_{r_1}\Phi\Big(\frac{c}{RN},\frac{r_1}{R}\Big)\ue\Big(a_1k(hn_2'-h'n_2)\frac{\overline{q_1n_2n_2'r_1}}{q_2c}\Big),
	\end{align*}
	and $\Phi(\cdot,\cdot)$ is a smooth weight function that majorizes the indicator function of the box $[1,4]\times[1,2]$.
	
	\subsection{Application of Kloostermania}\label{appli2}
	For $QN^2\lqs u\lqs 8QN^2,v\sim Q$ and $0\lqs w\lqs 4NHK\sL^B$, we put 
	\begin{align*}
		\cB_2(u,v,w)=\sum_{\substack{k\mid w,q_1n_2n_2'=u,q_2=v\\n_2\equiv n_2'\bmod {kq_1q_2}}}\sum_{\substack{|a_1k(hn_2'-h'n_2)|=w\\h\equiv h'\bmod {q_1q_2}}}\eta_2(h,n_2)\overline{\eta_2}(h',n_2').
	\end{align*}
	Then we can rewrite $\sN_2$ as
	\begin{align*}
		\sN_2=\sum_{\substack{QN^2\lqs u\lqs 8QN^2\\v\sim Q}}\sum_{0\lqs w\lqs 4NHK\sL^B}\cB_2(u,v,w)\sum_{c}\sum_{r_1}\Phi\Big(\frac{c}{RN},\frac{r_1}{R}\Big)\ue\Big(w\frac{\overline{ur_1}}{vc}\Big),
	\end{align*}
	which is of type $\sK(4QN^2,Q,4NHK\sL^B,RN,R)$ in Lemma $\ref{lm:kloostermania}$.
	
	The diagonal terms with $w=0$ contribute
	\begin{align}\label{eq:diagonal2}
		\sN_2(w=0)\ll x^\varepsilon KQ^2R^2N^2H.
	\end{align}
	
	For the off-diagonal terms with $w\neq 0$, we apply Cauchy's inequality to deduce
	\begin{align*}
		\sum_{u,v,w}|\cB_2(u,v,w)|^2&\ll (QRHKN)^\varepsilon \sum_{k,q_1,q_2}\sum_{n_2\equiv n_2'\bmod {kq_1q_2}}\sum_{w'> 0}\Big(\sum_{\substack{h\equiv h'\bmod {q_1q_2}\\|(hn_2'-h'n_2)|=w'}}1\Big)^2\notag\\
		&\ll x^\varepsilon\sum_{k,q_1,q_2}\sum_{\substack{h_1\equiv h_2\bmod {q_1q_2}\\h_3\equiv h_4\bmod {q_1q_2}}}\sum_{\substack{n_2\equiv n_2'\bmod {kq_1q_2}\\(h_1-h_3)n_2'=(h_2-h_4)n_2}}1 \notag\\
		&\ll x^\varepsilon(Q^{-2}N^2H^2+KNH^3+N^2H+KNH^2Q^2)
	\end{align*}
	assuming that
	\begin{align}\label{eq:HQ2}
		Q^2\lqs S.
	\end{align}
	
	If we further assume that $Q^2\lqs H$, we have
	\begin{align*}
		\sum_{u,v,w}|\cB_2(u,v,w)|^2\ll x^\varepsilon(Q^{-2}N^2H^2+KNH^3).
	\end{align*}
	Now by Lemma $\ref{lm:kloostermania}$ we have 
	\begin{align}\label{eq:kloostermania12}
		\sN_2^2(w\neq 0)&\ll x^\varepsilon(Q^4R^2N^5+Q^{5/2}R^3N^4)(Q^{-2}N^2H^2+KNH^3)\notag\\
		&\ll x^\varepsilon(Q^6R^{6}S^2M^{-2}N^{7}+Q^{9/2}R^{7}S^2M^{-2}N^{6})
	\end{align}
	assuming further that
	\begin{align}\label{eq:assume2}
		Q^4R^2N\lqs x^{1-\varepsilon}.
	\end{align}
	
	Combining $(\ref{eq:Klresult12})$, $(\ref{eq:diagonal2})$, $(\ref{eq:kloostermania12})$ and recalling the assumption $Q^2\lqs H$ we have (recall that $MN=x$)
	\begin{align}\label{eq:R_12}
		\sR_1\ll x^{\varepsilon}(Q^{3}RS^{-2}N^{5/2}x^{1/2}+Q^{5/2}R^{1/2}S^{-3/2}N^{11/4}x^{1/2}+Q^{17/8}R^{3/4}S^{-3/2}N^{5/2}x^{1/2}).
	\end{align}
	
	It remains to consider the case $Q^2>H$. Note that in this case, placing the sum over $h$ inside the absolute value would yield no saving for the diagonal terms, due to the extra summation over $a_6\Mod{q_1q_2}$ introduced when separating the variables. Instead, we may apply Cauchy's inequality starting from the form
	\begin{align*}
		\sum\ldots\sum\Big|\sum_{n_2}\Big|,
	\end{align*}
	which leads to a bound superior to $(\ref{eq:R_12})$ when $Q^2>H$:
	\begin{align*}
		\sR_1\ll x^{\varepsilon}(Q^3R^2N^3S^{-3/2}+Q^{5/2}RS^{-5/4}N^3x^{1/4}+Q^{17/8}R^{5/4}S^{-5/4}N^{11/4}x^{1/4}).
	\end{align*} 
	
	For simplicity, we employ the unified bound $(\ref{eq:R_12})$ in both cases.
	Hence, to guarantee the desired result $(\ref{eq:sR1})$, we require 
	\begin{align*}
		N \lqs x^{-\varepsilon}\min\{x^{1/3}Q^{-2}R^{-2/3}S^{2/3},x^{2/7}Q^{-10/7}R^{-2/7}S^{2/7},x^{1/3}Q^{-17/12}R^{-1/2}S^{1/3}\}.
	\end{align*}
	This requirement, together with assumptions $(\ref{eq:S2condition})$, $(\ref{eq:benzhi2})$, $(\ref{eq:HQ2})$, $(\ref{eq:assume2})$, establishes Theorem \ref{thm:TypeII-2}.

	\section{Estimation of $\sR_1$: the third method}\label{thirdmethod}
	
	In Sections \ref{firstmethod} and \ref{secondmethod}, we apply $\it{Kloostermania}$ to exploit the average over the moduli, at the cost of some loss introduced by the separation of variables (such as additional summations over residue classes $\Mod{q_1q_2}$). This loss becomes significant when $q$ is large, and we have not found an effective alternative to separate these variables. Therefore, for large $q$, we instead leverage the favorable factorization properties of the moduli to estimate the resultant exponential sums.
	
	\subsection{Preparation for the $q$-analogue of van der Corput method}\label{Fou2}
	Similar to previous arguments, we put $n_2=n_1-ks$. After a change of variable, we find $\sR_1$ is bounded by
	\begin{align}
		\frac{x^\varepsilon M}{RS}\sum_{\substack{1\lqs|k|\lqs K\\r_2,s,n_1}}\bigg|\sum_{\substack{q_1,q_2\\r_1,h}}&\frac{\eta_3(h,q_1,q_2,r_1)}{q_1q_2r_1}\ue\Big (a_1hk\frac{\overline{q_1r_1(n_1-ks)}}{q_2r_2n_1}\notag\\
		&+(a_2-a_1)h\frac{\overline{q_2r_1r_2sn_1}}{q_1}+(a_2-a_1)h\frac{\overline{q_1r_1r_2s(n_1-ks)}}{q_2}\Big)\bigg|\label{eq:trick2}
	\end{align}
	for some $\Vert\eta_3\Vert_\infty\lqs 1$. Then by Cauchy's inequality we have
	\begin{align}\label{eq:M}
		\sR_1^2\ll \frac{ MKx^{1+\varepsilon}}{Q^4R^3S}\sum_{\substack{k,q_1,\widetilde{q_1},q_2,\widetilde{q_2}\\r_1,\widetilde{r_1},r_2,h,\tilde{h},n_1}}\big|\sM\big|,
	\end{align}
	where \begin{align*}
		\sM=\sum_{s\sim S}\ue\Big(&a_1hk\frac{\overline{q_1r_1(n_1-ks)}}{q_2r_2n_1}+(a_2-a_1)h\frac{\overline{q_2r_1r_2n_1s}}{q_1}+(a_2-a_1)h\frac{\overline{q_1r_1r_2s(n_1-ks)}}{q_2}\\
		&-a_1\widetilde{h}k\frac{\overline{\widetilde{q_1}\widetilde{r_1}(n_1-ks)}}{\widetilde{q_2}r_2n_1}-(a_2-a_1)\widetilde{h}\frac{\overline{\widetilde{q_2}\widetilde{r_1}r_2n_1s}}{\widetilde{q_1}}-(a_2-a_1)\widetilde{h}\frac{\overline{\widetilde{q_1}\widetilde{r_1}r_2s(n_1-ks)}}{\widetilde{q_2}}\Big).
	\end{align*}
	
	To apply the $q$-analogue of van der Corput method to $\sM$, we hope that $q_1\widetilde{q_1},q_2\widetilde{q_2},r_2,n_1,s$ are pairwise coprime, that $(r_1\widetilde{r_1},q_1\widetilde{q_1}q_2\widetilde{q_2}r_2n_1)=1$ and that $q_1,\widetilde{q_1},q_2,\widetilde{q_2},r_2,n_1$ are square-free on the right hand side of $(\ref{eq:M}).$ As argued in Section $\ref{initial}$, it suffices to ensure $(q_1,\widetilde{q_2})=(q_2,\widetilde{q_1})=1$ (as well as $(q_1,\widetilde{r_1})=(r_1,\widetilde{q_1})=(q_2,\widetilde{r_1})=(r_1,\widetilde{q_2})=1$, which can be deduced by the same method). Similar to the application of Fouvry's trick (Lemma \ref{lm:splitting into coprime sets}) in Section $\ref{Fou1}$, we partition the summations over $q_1,q_2$ in \eqref{eq:trick2} into at most $\exp(\sL^{2/3+\varepsilon})$ $($since $w(q)\lqs \sL^{1/3}$ by Condition \ref{cond:TypeII}$)$ sub-sums. Then we apply Cauchy's inequality to each sub-sum to obtain \eqref{eq:M}, with the constraints that $(q_1,\widetilde{q_2})=(q_2,\widetilde{q_1})=1$ now in force. For simplification, we proceed to estimate $\sM$ with these constraints. 
	
	By applying the Chinese remainder theorem, it suffices to estimate an exponential sum in the single variable $s$ with modulus $q_1\widetilde{q_1}q_2\widetilde{q_2}r_2n_1:$
	\begin{align*}
		\sM=\sum_{s\sim S}\Psi(s),
	\end{align*}
	where
	\begin{align*}
		\Psi(s)=\ue\Big(&-a_1\xi\frac{\overline{q_1\widetilde{q_1}q_2\widetilde{q_2}r_1\widetilde{r_1}r_2s}}{n_1}+(a_2-a_1)\xi\frac{\overline{q_2\widetilde{q_2}r_1\widetilde{r_1}r_2n_1s}}{q_1\widetilde{q_1}}\\
		&+a_1k\xi\frac{\overline{q_1\widetilde{q_1}r_1\widetilde{r_1}n_1(n_1-ks)}}{q_2\widetilde{q_2}r_2}+(a_2-a_1)\xi \frac{\overline{q_1\widetilde{q_1}r_1\widetilde{r_1}r_2s(n_1-ks)}}{q_2\widetilde{q_2}}\Big),
	\end{align*}
	and
	\begin{align*}
		\xi&= h\widetilde{q_1}\widetilde{q_2}\widetilde{r_1}-\widetilde{h}q_1q_2r_1,~r_2\overline{r_2}\equiv 1\bmod {q_2\widetilde{q_2}}.
	\end{align*}
	
	For the diagonal terms with $\xi=0$ $(h\widetilde{q_1}\widetilde{q_2}\widetilde{r_1}=\widetilde{h}q_1q_2r_1)$, we estimate $\sM$ trivially, yielding a total contribution of
	\begin{align}\label{eq:ABdiagonal}
		\sR_{1}({\xi=0})\ll x^{1/2+\varepsilon}R^{1/2}S^{1/2}N^{1/2}K.
	\end{align}

	\subsection{Application of the $q$-analogue of van der Corput method ($BA^2B$-process)}
	In this subsection, we treat the off-diagonal terms via the $q$-analogue of van der Corput method for algebraic exponential sums. Especially, we apply the $BA^2B$-process to the case $Q^4RN\lqs S^2$. 
	
	For the first step $($the $B$-process$)$, we denote by $\widehat{\Psi}$ the Fourier transform of $\Psi$ $\it{w.r.t.}$ modulus $q_1\widetilde{q_1}q_2\widetilde{q_2}r_2n_1.$ Applying Lemma \ref{lm:B} yields
	\begin{align*}
		\sM\ll \frac{S|\widehat{\Psi}(0)|}{\sqrt{q_1\widetilde{q_1}q_2\widetilde{q_2}r_2n_1}}+\frac{S\sL }{\sqrt{q_1\widetilde{q_1}q_2\widetilde{q_2}r_2n_1}}\Big|\sum_{j\in \cJ}\frac{1}{\sqrt{q_1\widetilde{q_1}q_2\widetilde{q_2}r_2n_1}}\sum_{s\bmod {q_1\widetilde{q_1}q_2\widetilde{q_2}r_2n_1}}\Psi(s)\ue\Big(\frac{(a_7-s)j}{q_1\widetilde{q_1}q_2\widetilde{q_2}r_2n_1}\Big)\Big|
	\end{align*}
	for some integer $a_7$ and an interval $\cJ$ with $0\notin\cJ$ and $|\cJ|\ll Q^4RN/S$. An application of the Chinese remainder theorem then yields
	\begin{align}\label{eq:B1}
		\sM\ll \frac{S|\widehat{\Psi}(0)|}{\sqrt{q_1\widetilde{q_1}q_2\widetilde{q_2}r_2n_1}}+\frac{S\sL}{\sqrt{q_1\widetilde{q_1}q_2\widetilde{q_2}r_2n_1}}\Big|\sum_{j\in \cJ}\Psi_0(j)\Psi_1(j)\Psi_2(j)\Big|,
	\end{align} 
	where
	\begin{align*}
		&\Psi_0(j)=\frac{1}{\sqrt{n_1}}\sum_{s_0\bmod {n_1}}\ue\Big(-a_1\xi\frac{\overline{q_1^2{\widetilde{q_1}}^2q_2^2\widetilde{q_2}^2r_1\widetilde{r_1}r_2^2s_0}}{n_1}+\frac{(\overline{q_1\widetilde{q_1}q_2\widetilde{q_2}r_2}a_7-s_0)j}{n_1}\Big),\\
		&\Psi_1(j)=\frac{1}{\sqrt{q_1\widetilde{q_1}}}\sum_{s_1\bmod {q_1\widetilde{q_1}}}\ue\Big((a_2-a_1)\xi\frac{\overline{q_2^2\widetilde{q_2}^2r_1\widetilde{r_1}r_2^2n_1^2s_1}}{q_1\widetilde{q_1}}+\frac{(\overline{q_2\widetilde{q_2}r_2n_1}a_7-s_1)j}{q_1\widetilde{q_1}}\Big),\\
		&\Psi_2(j)=\frac{1}{\sqrt{q_2\widetilde{q_2}r_2}}\sum_{s_2\bmod {q_2\widetilde{q_2}r_2}}\ue\Big(a_1k\xi\frac{\overline{q_1\widetilde{q_1}r_1\widetilde{r_1}n_1(n_1-kq_1\widetilde{q_1}n_1s_2)}}{q_2\widetilde{q_2}r_2}\\ &\quad\quad\quad\ \  +(a_2-a_1)\xi r_2\overline{r_2}\frac{\overline{q_1^2\widetilde{q_1}^2r_1\widetilde{r_1}n_1s_2(n_1-kq_1\widetilde{q_1}n_1s_2)}}{q_2\widetilde{q_2}r_2}+\frac{(\overline{q_1\widetilde{q_1}n_1}a_7-s_2)j}{q_2\widetilde{q_2}r_2}\Big),
	\end{align*}
	and $r_2\overline{r_2}\equiv 1\bmod{q_2\widetilde{q_2}}$ in $\Psi_2(j).$
	Observe that both $\Psi_0(j)$ and $\Psi_1(j)$ are Kloosterman sums attached with an additive character:
	\begin{align*}
		\Psi_0(j)=S(j,\xi_0;n_1)\ue\Big(\frac{\overline{q_1\widetilde{q_1}q_2\widetilde{q_2}r_2}a_7j}{n_1}\Big),~\Psi_1(j)=S(j,\xi_1;q_1\widetilde{q_1})\ue\Big(\frac{\overline{q_2\widetilde{q_2}r_2n_1}a_7j}{q_1\widetilde{q_1}}\Big),
	\end{align*}
	where $\xi_0\equiv a_1\xi\overline{q_1^2{\widetilde{q_1}}^2q_2^2\widetilde{q_2}^2r_1\widetilde{r_1}r_2^2}\bmod{n_1}$ and $\xi_1\equiv (a_1-a_2)\xi\overline{q_2^2\widetilde{q_2}^2r_1\widetilde{r_1}r_2^2n_1^2}\bmod{q_1\widetilde{q_1}}.$
	
	For the second step (the $A^2B$-process), by Lemma \ref{lm:Acoro} we have 
	\begin{align}\label{eq:A2}
		\Big|\sum_{j\in \cJ}&\Psi_0(j)\Psi_1(j)\Psi_2(j)\Big|^4\ll \frac{Q^4R^2}{|\cJ|^2}\cU^2_{0,1}\,\cU^2_2+Q^2\cU_0\cU_1\cU_2^2|\cJ|^{\varepsilon-1}\notag\\
		&+\frac{Q^4RN^{1/2}\cU_1\cU^2_2}{|\cJ|}\sum_{\substack{0<|l_1|< |\cJ|/(q_1\widetilde{q_1})\\0<|l_2|< |\cJ|/(q_2\widetilde{q_2}r_2)}}\frac{1}{|\cJ|-|l_2q_2\widetilde{q_2}r_2|}\Big(\frac{|\cJ||\widehat{Z}(0)|}{N}+\sum_{1\lqs |v|\lqs n_1/2}\frac{|\widehat{Z}(v)|}{|v|}\Big) ,
	\end{align} 
	where
	\begin{align*}
		&\cU_{0,1}=\sum_{j\in \cJ}\big|\Psi_0(j)\Psi_1(j)\big|^2,\ \cU_0=\sum_{j\in \cJ}\big|\Psi_0(j)\big|^2,\ 
		\cU_1=\sum_{j\in\cJ_1}\big|\Psi_1(j)\big|^2,\ 
		\cU_2=\sum_{j\in\cJ_2}\big|\Psi_2(j)\big|^2,\\
		&Z(v):=\Psi_0(v)\overline{\Psi_0(v+l_1q_1\widetilde{q_1})}\overline{\Psi_0(v+l_2q_2\widetilde{q_2}r_2)}\Psi_0(v+l_1q_1\widetilde{q_1}+l_2q_2\widetilde{q_2}r_2)\\
		&\ \qquad =n_1^{-2}S(v,\xi_0;n_1)S(v+l_1q_1\widetilde{q_1},\xi_0;n_1)S(v+l_2q_2\widetilde{q_2}r_2,\xi_0;n_1)S(v+l_1q_1\widetilde{q_1}+l_2q_2\widetilde{q_2}r_2,\xi_0;n_1),
	\end{align*}
	$\cJ_1,\cJ_2$ are intervals satisfying $|\cJ_1|,|\cJ_2|\lqs 2|\cJ|,$ and $\widehat{Z}$ is the Fourier transform of $Z$ modulo $n_1.$ Now it remains to estimate $\cU_{0,1},\cU_{0},\cU_{1},\cU_{2},$ and 
	\begin{align*}
		\widehat{Z}(v)&=n_1^{-5/2}\sum_{t\bmod{n_1}}Z(t)\ue\Big(\frac{-vt}{n_1}\Big)\\
		&=n_1^{-5/2}\sum_{t\bmod{n_1}}S(t,\xi_0;n_1)S(t+l_1q_1\widetilde{q_1},\xi_0;n_1)S(t+l_2q_2\widetilde{q_2}r_2,\xi_0;n_1)\\
		&\ \ \qquad\qquad\qquad\cdot S(t+l_1q_1\widetilde{q_1}+l_2q_2\widetilde{q_2}r_2,\xi_0;n_1)\ue\Big(\frac{-vt}{n_1}\Big)
	\end{align*}
	for $0\lqs|v|\lqs n_1/2$. To this end, we apply Lemma $\ref{lm:fourier}$ to obtain 
	\begin{align}\label{eq:Z}
		|\widehat{Z}(v)|\ll x^\varepsilon(l_1,v,n_1)^{1/2}(l_2,v,n_1)^{1/2}(\xi,l_1,l_2,n_1)^{3/2},
	\end{align}
	and the following result to bound $\cU_{0,1},~\cU_{0},~\cU_{1},~\cU_{2}.$
	\begin{lemma}\label{lm:l2norm} Let the notation and conditions be as above. We have
		\begin{align*}
			\cU_{0,1},~\cU_{0},~\cU_{1}\ll x^\varepsilon|\cJ|,
		\end{align*}
		and
		\begin{align*}
			\cU_{2}\ll x^\varepsilon|\cJ|(q_2,\widetilde{q_2}).
		\end{align*}
	\end{lemma}
	\proof Applying Lemma \ref{lm:Weil} we have
	\begin{align*}
		\Psi_0(j)\ll x^{\varepsilon/10} (j,\xi,n_1)^{1/2},~\Psi_1(j)\ll x^{\varepsilon/10} (j,(a_1-a_2)\xi,q_1\widetilde{q_1})^{1/2}.
	\end{align*}
	Applying Lemma \ref{lm:algebraic type exponential sums} by choosing 
	\begin{align*}
		&\lambda_1(y)=(a_2-a_1)\xi r_2\overline{r_2}+a_1k\xi q_1\widetilde{q_1}y-jq_1^2\widetilde{q_1}^2r_1\widetilde{r_1}n_1^2y^2+jkq_1^3\widetilde{q_1}^3r_1\widetilde{r_1}n_1^2y^3,\\
		&\lambda_2(y)=q_1^2\widetilde{q_1}^2r_1\widetilde{r_1}n_1^2y-kq_1^3\widetilde{q_1}^3r_1\widetilde{r_1}n_1^2y^2,
	\end{align*}
	where $r_2\overline{r_2}\equiv 1\bmod{q_2\widetilde{q_2}},$ we have
	\begin{align*}
		\Psi_2(j)\ll x^{\varepsilon/10}(j,q_2\widetilde{q_2}r_2)^{1/2}(q_2,\widetilde{q_2})^{1/2}.
	\end{align*}
	Thus, we finish the proof by the following estimates:
	\begin{align*}
		&\cU_{0,1}\ll x^{\varepsilon/2}\sum_{j\in\cJ}(j,n_1)(j,q_1\widetilde{q_1})=x^{\varepsilon/2}\sum_{j\in\cJ}(j,n_1q_1\widetilde{q_1})\ll x^\varepsilon|\cJ|,\\
		&\cU_{0}\ll x^{\varepsilon/2}\sum_{j\in\cJ}(j,n_1)\ll x^\varepsilon|\cJ|,\\
		&\cU_{1}\ll x^{\varepsilon/2}\sum_{j\in\cJ_1}(j,q_1\widetilde{q_1})\ll x^\varepsilon|\cJ|,\\
		&\cU_{2}\ll x^{\varepsilon/2}(q_2,\widetilde{q_2})\sum_{j\in\cJ_2}(j,q_2\widetilde{q_2}r_2)\ll x^\varepsilon|\cJ|(q_2,\widetilde{q_2}).
	\end{align*}
	\endproof
	
	Applying \eqref{eq:Z} and partial summation, we have
	\begin{align}\label{eq:fourier2}
		\sum_{\substack{0<|l_1|< |\cJ|/(q_1\widetilde{q_1})\\0<|l_2|< |\cJ|/(q_2\widetilde{q_2}r_2)}}\frac{1}{|\cJ|-|l_2q_2\widetilde{q_2}r_2|}\Big(\frac{|\cJ||\widehat{Z}(0)|}{N}+\sum_{1\lqs |v|\lqs n_1/2}\frac{|\widehat{Z}(v)|}{|v|}\Big)\ll  \frac{x^\varepsilon|\cJ|(n_1,\xi)^{1/2}}{Q^4RN}(|\cJ|+N).
	\end{align}

	In conclusion, combining $\eqref{eq:M}$-$\eqref{eq:A2}$, \eqref{eq:fourier2} and Lemma \ref{lm:l2norm}, we have
	\begin{align}\label{eq:R11}
		\sR_1&\ll x^{1/2+\varepsilon}R^{1/2}S^{1/2}N^{1/2}K+\frac{x^\varepsilon \sqrt{MKx}}{\sqrt{Q^4R^3S}}\Big(KQ^4R^3H^2N\cdot\frac{S}{\sqrt{Q^4RN}}\Big(1+\Big(\frac{Q^4RN}{S}\Big)^{1/2}(Q^2R)^{1/2}\notag\\
		&\ \ \ \ +\Big(\frac{Q^4RN}{S}\Big)^{3/4}(Q^2)^{1/4}+\Big(\frac{Q^4RN}{S}\Big)^{3/4}(N)^{1/8}\Big) \Big)^{1/2}\notag\\
		&\ll x^\varepsilon\big(R^{1/2}S^{-1/2}N^{3/2}x^{1/2}+Q^{11/4}R^{17/8}S^{-3/8}N^{17/8}+Q^{5/2}R^{17/8}S^{-3/8}N^{35/16}\big)
	\end{align}
	assuming that
	\begin{align}\label{eq:assume8}
		RS\lqs N^{3/2},~Q^4R\lqs S.
	\end{align}
	\begin{remark}\label{rm:simple}
		To see the $A^2B$-process above more briefly, one may assume that $\Vert\Psi_i\Vert_\infty\ll 1$ and $\Vert \widehat{Z}\Vert_\infty\ll x^\varepsilon,$ to apply $(\ref{eq:Acoro2})$ in Remark $\ref{briefly}$ instead of complicated $(\ref{eq:Acoro1}),$ which leads to
		\begin{align*}
			\Big|\sum_{j\in \cJ}\Psi_0(j)\Psi_1(j)\Psi_2(j)\Big|^4&\ll x^\varepsilon(|\cJ|^2(Q^2R)^2+|\cJ|^3Q^2+|\cJ|^3N^{1/2}).
		\end{align*}
		Thus we can see that $BA^2B$-process gives the estimate
		\begin{align*}
			\sM\ll\frac{x^\varepsilon S}{\sqrt{Q^4RN}}\Big(1+\Big(\frac{Q^4RN}{S}\Big)^{1/2}(Q^2R)^{1/2}+\Big(\frac{Q^4RN}{S}\Big)^{3/4}(Q^2)^{1/4}+\Big(\frac{Q^4RN}{S}\Big)^{3/4}(N)^{1/8}\Big). 
		\end{align*}
	\end{remark}
	
	\subsection{Application of the $q$-analogue of van der Corput method ($A^2B$-process)}
	In this subsection, we treat the off-diagonal terms with $Q^4RN>S^2,$ still via the $q$-analogue of van der Corput method. Observe that in this case, the dual sum is longer than the original one $(Q^4RN/S>S).$ We therefore apply the $A^2B$-process instead. 
	
	We rewrite
	\begin{align*}
		\sM=\sum_{s\sim S}\Upsilon_1(s)\Upsilon_2(s)\Upsilon_3(s),
	\end{align*}
	where
	\begin{align*}
		&\Upsilon_0(s):=\ue\Big(-a_1\xi\frac{\overline{q_1\widetilde{q_1}q_2\widetilde{q_2}r_1\widetilde{r_1}r_2s}}{n_1}\Big),~\Upsilon_1(s):=\ue\Big((a_2-a_1)\xi\frac{\overline{q_2\widetilde{q_2}r_1\widetilde{r_1}r_2n_1s}}{q_1\widetilde{q_1}}\Big),\\
		&\Upsilon_2(s):=\ue\Big(a_1k\xi\frac{\overline{q_1\widetilde{q_1}r_1\widetilde{r_1}n_1(n_1-ks)}}{q_2\widetilde{q_2}r_2}+(a_2-a_1)\xi\frac{\overline{q_1\widetilde{q_1}r_1\widetilde{r_1}r_2s(n_1-ks)}}{q_2\widetilde{q_2}}\Big)
	\end{align*}
	By Lemma \ref{lm:Acoro} we have 
	\begin{align*}
		|\sM|^4\ll Q^4R^2S^2+Q^2S^{3+\varepsilon}+Q^4RS^2N^{1/2}\sum_{\substack{0<|l_1|< S/(q_1\widetilde{q_1})\\0<|l_2|< S/(q_2\widetilde{q_2}r_2)}}\frac{1}{S-|l_2q_2\widetilde{q_2}r_2|}\Big(\frac{S|\widehat{Y}(0)|}{N}+\sum_{1\lqs |v|\lqs n_1/2}\frac{|\widehat{Y}(v)|}{|v|}\Big) ,
	\end{align*}
	where
	\begin{align*}
		&Y(v):=\Upsilon_0(v)\overline{\Upsilon_0(v+l_1q_1\widetilde{q_1})}\overline{\Upsilon_0(v+l_2q_2\widetilde{q_2}r_2)}\Upsilon_0(v+l_1q_1\widetilde{q_1}+l_2q_2\widetilde{q_2}r_2),
	\end{align*}
	and $\widehat{Y}$ is the Fourier transform of $Y$ modulo $n_1.$ 
	Let $\xi_2\equiv -a_1\xi\overline{q_1\widetilde{q_1}q_2\widetilde{q_2}r_1\widetilde{r_1}r_2}\bmod{n_1}.$
	Applying Lemma \ref{lm:algebraic type exponential sums}, we have
	\begin{align*}
		\widehat{Y}(v)&=\frac{1}{\sqrt{n_1}}\sum_{t\bmod {n_1}}\ue\Big(\xi_2\frac{\overline{t}-\overline{t+l_1q_1\widetilde{q_1}}-\overline{t+l_2q_2\widetilde{q_2}r_2}+\overline{t+l_1q_1\widetilde{q_1}+l_2q_2\widetilde{q_2}r_2}}{n_1}-\frac{vt}{n_1}\Big)\notag\\
		&\ll x^\varepsilon(v,\xi l_1l_2,n_1)^{1/2}.
	\end{align*}
	Hence
	\begin{align*}
		\sM\ll x^\varepsilon(Q(SR)^{1/2}+Q^{1/2}S^{3/4}+N^{1/8}S^{3/4}).
	\end{align*}
	(This can be compared with Remark \ref{rm:simple}.)

	Recalling \eqref{eq:M}, \eqref{eq:ABdiagonal}, we now obtain that
	\begin{align}\label{eq:R12}
		\sR_1&\ll x^{\varepsilon}(xRSN)^{1/2}K+\frac{x^\varepsilon \sqrt{MKx}}{\sqrt{Q^4R^3S}}\Big(KQ^4R^3H^2N(Q(SR)^{1/2}+Q^{1/2}S^{3/4}+N^{1/8}S^{3/4})\Big)^{1/2}\notag\\
		&\ll x^\varepsilon\big(R^{1/2}S^{-1/2}N^{3/2}x^{1/2}+Q^{5/2}R^{9/4}S^{-1/4}N^{2}+Q^{9/4}R^{2}S^{-1/8}N^{2}+Q^{2}R^{2}S^{-1/8}N^{33/16}\big).
	\end{align}
	
	In conclusion, by \eqref{eq:R11}, \eqref{eq:R12} and assumption $(\ref{eq:assume8})$, we require
	\begin{align*}
		RS\lqs N^{3/2},~Q^4R\lqs S,~N\lqs x^{-\varepsilon}\min\{xR^{-1}S^{-1},x^{8/9}Q^{-22/9}R^{-17/9}S^{-5/9},x^{16/19}Q^{-40/19}R^{-34/19}S^{-10/19}\};
	\end{align*}
	or
	\begin{align*}
		N\lqs x^{-\varepsilon}\min\{xR^{-1}S^{-1},xQ^{-5/2}R^{-9/4}S^{-3/4},xQ^{-9/4}R^{-2}S^{-7/8},x^{16/17}Q^{-32/17}R^{-32/17}S^{-14/17}\}.
	\end{align*}
	to guarantee $(\ref{eq:sR1})$. These requirements, together with assumptions $(\ref{eq:S2condition})$ and $(\ref{eq:benzhi2})$,  establish Theorem \ref{thm:TypeII-3}.

	\section{The Type \textrm{I} estimate}\label{sec:typeI}
	In this section, we work with the sum 
	\begin{align*}
		\fS_{\textrm{I}}&=\sum_{(q,rs)=1}\gamma_q\delta_r\lambda_s\bigg(\sum_{\substack{P^-(mn)>\sL^C\\mn\equiv a_1\bmod {rs}\\mn\equiv a_2\bmod q}}\beta_n-\frac{1}{\varphi(qrs)}\sum_{\substack{P^-(mn)>\sL^C\\(mn,qrs)=1}}\beta_n\bigg)
	\end{align*}
	under the following condition.
	\begin{condition}\label{cond:TypeI}
		Suppose $MN=x,$ $|a_1|\lqs \sL^B$ for some $B>0$ and $|a_2|\lqs x$ with $a_1\neq a_2.$ Let $\bgamma=(\gamma_q),$ $\bdelta=(\delta_r),$ $\blambda=(\lambda_s)$ and $\bbeta=(\beta_n)$ be divisor-bounded sequences supported on positive integers with
		\begin{align*}
			q\sim Q, ~r\sim R, ~s\sim S, ~n\sim N\lqs x^{1-\varepsilon},~(rs,a_1)=(q,a_2rs)=1,
		\end{align*}
		and the variable $m\sim M$ is assumed.
	\end{condition}
	
	\begin{theorem}\label{thm:TypeI}
		Under Condition $\ref{cond:TypeI},$ the desired bound $(\ref{eq:TypeI})$ holds provided that $M\gqs QRSx^{\varepsilon}$ or 
		\begin{align*}
			M\gqs x^{\varepsilon}\max\{Q^{1/2}R^{-1/2}S^{-1/2}x^{1/2},Q^3S,Q^5R^4Sx^{-1},Q^{4}R^{3}S^2x^{-1}\},\ Q\lqs R,\ Q^{2}R^2S\lqs x^{1-\varepsilon}.
		\end{align*}
	\end{theorem}
	
	In order to apply Poisson summation directly, we need to remove the restriction $P^-(m)>\sL^C.$ To this end, we introduce the ``fundamental lemma" of sieve theory, which was coined by Halberstam and Richert \cite[Theorem 2.5]{HR74}.
	\begin{lemma}[Fundamental lemma]\label{lm:fundamentallemma}
		Let $z\gqs 2$ and $D=z^u.$ Let $\cP$ be a set of primes and
		\begin{align*}
			\cP(z)=\prod_{\substack{p\in \cP,p<z}}p.
		\end{align*}
		There exist two real sequences $(\varpi^{\pm}_d),$ supported on $d\mid \cP(z)$ with $d\lqs D,$ satisfying the following properties$:$
		\begin{itemize}
			\item $|\varpi^\pm_1|=1$ and $|\varpi^{\pm}_d|\lqs 1.$
			\item For every $l\mid \cP(z)$ with $l>1,$
			\begin{align*}
				\sum_{d\mid l}\varpi^-_d\lqs 0\lqs \sum_{d\mid l}\varpi^+_d.
			\end{align*}
			\item 
			\begin{align*}
				\sum_{d\mid \cP(z)}\frac{\varpi^\pm_d}{d}=\prod_{p\mid \cP(z)}\Big(1-\frac{1}{p}\Big)\{1+O(e^{-u})\}.
			\end{align*}
		\end{itemize}
	\end{lemma}
	\proof This follows from \cite[Lemma 5]{FI78}. \endproof
	Without loss of generality, we may assume that $\beta_n\gqs 0$ (otherwise, we treat $\beta^+_n:=\max\{0,\beta_n\}$ and $\beta^-_n:=\min\{0,\beta_n\}$ separately) and $\beta_n=0$ for $P^-(n)\lqs \sL^C.$
	Now we take $\cP$ to be the set of primes not dividing $qrs$ and $z=\sL^C.$
	Following the notation in Lemma \ref{lm:fundamentallemma}, we obtain
	\begin{align*}
		\cE(q,r,s)&:=\sum_{\substack{P^-(m)>\sL^C\\mn\equiv a_1\bmod {rs}\\mn\equiv a_2\bmod q}}\beta_n-\frac{1}{\varphi(qrs)}\sum_{\substack{P^-(m)>\sL^C\\(mn,qrs)=1}}\beta_n\\
		&\lqs\sum_{\substack{d|\cP(z)}}\Big(\sum_{\substack{dm\sim M,n\sim N\\dmn\equiv a_1\bmod {rs}\\dmn\equiv a_2\bmod q}}\varpi^+_d\beta_n-\frac{1}{\varphi(qrs)}\sum_{\substack{dm\sim M,n\sim N\\(mn,qrs)=1}}\varpi^-_d\beta_n\Big)\\
		&=\cE^+(q,r,s)+\Delta(q,r,s),
	\end{align*}
	where we write
	\begin{align*}
		&\cE^+(q,r,s)=\sum_{\substack{d|\cP(z)}}\varpi^+_d\Big(\sum_{\substack{dm\sim M,n\sim N\\dmn\equiv a_1\bmod {rs}\\dmn\equiv a_2\bmod q}}\beta_n-\frac{1}{\varphi(qrs)}\sum_{\substack{dm\sim M,n\sim N\\(mn,qrs)=1}}\beta_n\Big),\\
		&\Delta(q,r,s)=\sum_{\substack{d|\cP(z)}}\frac{\varpi^+_d-\varpi^-_d}{\varphi(qrs)}\sum_{\substack{dm\sim M,n\sim N\\(mn,qrs)=1}}\beta_n.
	\end{align*}
	Analogously, we have the lower bound
	\begin{align*}
		\cE(q,r,s)\gqs\cE^-(q,r,s)-\Delta(q,r,s).
	\end{align*}
	Now we choose $D=x^{\varepsilon/2}$ in Lemma \ref{lm:fundamentallemma}. Together with the elementary estimate
	\begin{align*}
		\frac{1}{\varphi(qrs)}\sum_{\substack{m\sim M/d\\(m,qrs)=1}}=\frac{M}{dqrs}+O\Big(\frac{\tau(qrs)}{\varphi(qrs)}\Big),
	\end{align*}
	Lemma \ref{lm:fundamentallemma} yields
	\begin{align*}
		\Delta(q,r,s)&=\bigg\{\frac{M}{qrs}\sum_{\substack{d|\cP(z)}}\frac{\varpi^+_d-\varpi^-_d}{d}+O\Big(\frac{\tau(qrs)}{\varphi(qrs)}\Big)\bigg\}\sum_{(n,qrs)=1}\beta_n\\
		&\ll \Big(\frac{M}{qrs}\exp(\sL^{\varepsilon-1})+\frac{\tau(qrs)}{\varphi(qrs)}\Big)\sum_{(n,qrs)=1}|\beta_n|.
	\end{align*}
	Summing over $q,r,s,$ we obtain (without loss of generality, we assume that $(\gamma_q),$ $(\delta_r)$ and $(\lambda_s)$ are real sequences)
	\begin{align*}
		\fS_{\textrm{I}}\lqs& \sum_{q,r,s}(\gamma^+_q\delta^+_r\lambda^+_s+\gamma^+_q\delta^-_r\lambda^-_s+\gamma^-_q\delta^+_r\lambda^-_s+\gamma^-_q\delta^-_r\lambda^+_s)\cE^+\\
		&+\sum_{q,r,s}(\gamma^-_q\delta^-_r\lambda^-_s+\gamma^-_q\delta^+_r\lambda^+_s+\gamma^+_q\delta^-_r\lambda^+_s+\gamma^+_q\delta^+_r\lambda^-_s)\cE^-+O(x\sL^{-A}),
	\end{align*}
	where the notation $\gamma^\pm_q,\delta^\pm_r,\lambda^\pm_s$ is analogous to that for $\beta^\pm_n.$
	It now remains to prove, for example, that
	\begin{align*}
		\sum_{d\lqs x^{\varepsilon/2}}\Big|\sum_{(qrs,d)=1}\gamma^+_q\delta^+_r\lambda^+_s\Big(\sum_{\substack{m\sim M/d,n\sim dN\\mn\equiv a_1\bmod {rs}\\mn\equiv a_2\bmod q}}\beta_{n/d}-\frac{1}{\varphi(qrs)}\sum_{\substack{m\sim M/d,n\sim dN\\(mn,qrs)=1}}\beta_{n/d}\Big)\Big|=O(x^{1-\varepsilon}).
	\end{align*}
	To simplify the expression, we continue to write $M,N,\gamma_q,\delta_r,\lambda_s,\beta_n,\fS_{\textrm{I}}$ (the corresponding rescaling of coefficients are understood). Thus, it suffices to establish the estimate:
	\begin{align*}
		\fS_{\textrm{I}}:=\sum_{q,r,s}\gamma_q\delta_r\lambda_s\Big(\sum_{\substack{m\sim M,n\sim N\\mn\equiv a_1\bmod {rs}\\mn\equiv a_2\bmod q}}\beta_n-\frac{1}{\varphi(qrs)}\sum_{\substack{m\sim M,n\sim N\\(mn,qrs)=1}}\beta_n\Big)=O(x^{1-\varepsilon}).
	\end{align*}

	Let $f$ be a smooth function supported on $[M-M^{1-\varepsilon},2M+M^{1-\varepsilon}]$ satisfying
	$$\begin{cases}
		f(m)=1, & \text{for~}m\in [M,2M], \\
		f(m)\gqs 0, & \text{for~}m\in \bR,\\
		f^{(j)}(m)\ll_j M^{(\varepsilon-1)j}, & \text{for~all~}j\gqs 0.
	\end{cases}$$
	By Poisson summation we have 
	\begin{align*}
		\fS_{\textrm{I}}=\sum_{(q,rs)=1}\frac{\gamma_q\delta_r\lambda_s}{qrs}\sum_{(n,qrs)=1}\beta_n\sum_{1\lqs |h|\lqs H_0}\widehat{f}\Big(\frac{h}{qrs}\Big)\ue\Big(a_1h\frac{\overline{qn}}{rs}+a_2h\frac{\overline{rsn}}{q}\Big)+O(x^{1-\varepsilon}),
	\end{align*}
	where $H_0=x^\varepsilon QRSM^{-1}$. 
	Note that the term $O(x^{1-\varepsilon})$ comes from the rapid decay of $\widehat{f},$ and we are done if $M>QRS x^\varepsilon.$ We now henceforth assume that $M\lqs QRS x^\varepsilon.$
	From Cauchy's inequality it follows that
	\begin{align}\label{eq:guanghua0}
		\fS^2_{\textrm{I}} &\ll \frac{x^{1+\varepsilon}QM}{R^2S}\sum_{q,s,n}\sum_{\substack{a_8\bmod q\\a_9\bmod q}}\bigg|\sum_{r\equiv a_8\bmod q}\sum_{\substack{1\lqs |h|\lqs H_0\\h\equiv a_9\bmod q}}\eta_5(r,h)\ue\Big(a_1h\frac{\overline{qn}}{rs}\Big)\bigg|^2+x^{1-\varepsilon} \notag\\
		&:= \frac{x^{1+\varepsilon}QM}{R^2S}\sN_0 +x^{1-\varepsilon},
	\end{align}
	where 
	\begin{align*}
		\sN_0=\sum_{u\sim Q}\sum_{v\sim R^2}\sum_{\substack{0\lqs w\lqs 4a_1RH_0\sL^B}}	\cB_0(u,v,w)\sum_{s}\sum_{n}\Phi\Big(\frac{s}{S},\frac{n}{N}\Big)\ue\Big(w\frac{\overline{un}}{vs}\Big)
	\end{align*}
	with
	\begin{align*}
		\cB_0(u,v,w)=\sum_{q=u}\sum_{\substack{rr'=v\\r\equiv r'\bmod q}}\sum_{\substack{|a_1(h'r-hr')|=w\\h\equiv h'\bmod q}}\eta_5(r,h)\overline{\eta_5}(r',h'),
	\end{align*}
	a smooth weight function $\Phi(\cdot,\cdot)$ that majorizes the indicator function of the square $[1,2]^2,$ and some $\Vert\eta_5\Vert_\infty\lqs 1$. Still $\sN_0$ is of type $\sK(Q,R^2,4RH_0\sL^B,S,N)$ in Lemma \ref{lm:kloostermania}. 
	
	The diagonal terms with $w=0$ contribute
	\begin{align}\label{eq:diagonal3}
		\sN_0(w=0)\ll x^\varepsilon QRSNH_0.
	\end{align}
	
	For off-diagonal terms we calculate
	\begin{align*}
		\sum_{u,v,w}|\cB_0(u,v,w)|^2&\ll (QRH_0)^\varepsilon\sum_{q}\sum_{\substack{r\equiv r'\bmod q}}\sum_{w\gqs 0}\Big(\sum_{\substack{|a_1(h'r-hr')|=w\\h\equiv h'\bmod q}}1\Big)^2\\
		&\ll x^\varepsilon(Q^{-1}R^2H_0^2+RH_0^3)
	\end{align*}
	assuming that
	\begin{align}\label{eq:assume4}
		Q\lqs R
	\end{align}
	and $Q\lqs H_0$. Similar to the arguments in Section \ref{appli2}, if $Q> H_0$, we may apply Cauchy's inequality starting from the form
	\begin{align*}
		\sum\ldots\sum\Big|\sum_{r}\Big|,
	\end{align*}
	rather than using $(\ref{eq:guanghua0})$. This approach yields a better bound compared to $(\ref{eq:3})$, allowing us to omit the assumption $Q\lqs H_0$.
	
	Now, by Lemma $\ref{lm:kloostermania}$ we have
	\begin{align}\label{eq:guanghuakl2}
		\sN_0(w\neq 0)^2&\ll x^\varepsilon \big(QR^4S(S+QN)+QR^3S^2N\big)Q^{-1}R^2H_0^2\notag\\
		&= x^\varepsilon\big(Q^2R^8S^4M^{-2}+Q^3R^8S^3M^{-2}N+Q^2R^7S^4M^{-2}N\big)
	\end{align}	
	assuming that
	\begin{align}\label{eq:assume3}
		M>x^{\varepsilon}\{Q^2S,QR^{-1}S^{-1}N\}.
	\end{align}
	Finally, by $(\ref{eq:guanghua0})$, $(\ref{eq:diagonal3})$, $(\ref{eq:guanghuakl2})$ and assumptions $(\ref{eq:assume4})$, $(\ref{eq:assume3})$ we have $($recall that $MN=x$ and $H_0=x^\varepsilon QRSM^{-1})$
	\begin{align}\label{eq:3}
		\fS_{\textrm{I}}\ll x^\varepsilon\big(Q^{3/2}S^{1/2}M^{1/2}N+QRS^{1/2}M^{1/2}N^{1/2}+Q^{5/4}RS^{1/4}M^{1/2}N^{3/4}+QR^{3/4}S^{1/2}M^{1/2}N^{3/4}\big),
	\end{align}
	which, yields Theorem \ref{thm:TypeI} (recall that the case $M> QRSx^{\varepsilon}$ is trivial).

	\section{The Type \textrm{III} estimate}\label{sec:typeIII}
	Throughout this section, we mainly follow the idea of \cite{Po14} to work with the sum (similar to Section \ref{sec:typeI}, we remove the constraint $P^-(m_1m_2m_3n)>\sL^C$)
	\begin{align*}
		\fS_{\textrm{III}}=\sum_{(q,rs)=1}\gamma_q\delta_r\lambda_s\bigg(\sum_{\substack{m_1m_2m_3n\equiv a_1\bmod {rs}\\m_1m_2m_3n\equiv a_2\bmod q}}\beta_n-\frac{1}{\varphi(qrs)}\sum_{(m_1m_2m_3n,qrs)=1}\beta_n\bigg)
	\end{align*}
	under the following condition.
	
	\begin{condition}\label{cond:TypeIII}
		Suppose $M_1M_2M_3N=x,$ $|a_1|\lqs \sL^B$ for some $B>0$ and $|a_2|\lqs x$ with $a_1\neq a_2.$ Let $\bgamma=(\gamma_q),$ $\bdelta=(\delta_r),$ $\blambda=(\lambda_s)$ and $\bbeta=(\beta_n)$ be divisor-bounded sequences supported on positive integers with
		\begin{align*}
			q\sim Q, ~r\sim R, ~s\sim S, ~n\sim N,~(rs,a_1)=(q,a_2)=\mu^2(qrs)=1,
		\end{align*}
		and the variables $m_1\sim M_1,m_2\sim M_2,m_3\sim M_3$ are assumed.
		Assume further that $\bbeta$ satisfies the Siegel--Walfisz condition $($Definition $\ref{Def:SW})$.
	\end{condition}

	\begin{theorem}\label{thm:TypeIII}
		Under Condition $\ref{cond:TypeIII},$ the desired bound $(\ref{eq:TypeIII})$ holds provided that
		\begin{align*}
			\min\{M_1,M_2,M_3\}\gqs x^\varepsilon,\ \ M_1M_2M_3\gqs x^{\varepsilon}\max\{Q^4R^4S^{3}x^{-1},Q^{7/4}R^{7/4}S^{3/2},Q^{3/2}R^{3/2}S^2\}.
		\end{align*}
	\end{theorem}
	
	To prove Theorem \ref{thm:TypeIII}, we now let
	$\min\{M_1,M_2,M_3\}\gqs x^\varepsilon.$ For $i=1,2,3$, denote by $f_i$ a smooth function supported on $[M_i-M_i^{1-\varepsilon},2M_i+M_i^{1-\varepsilon}]$ such that
	$$\begin{cases}
		f_i(m)=1, & \text{for~}m\in [M_i,2M_i], \\
		f_i(m)\gqs 0, & \text{for~}m\in \bR,\\
		f^{(j)}_i(m)\ll_j M_i^{(\varepsilon-1)j}, & \text{for~all~} j\gqs 0.
	\end{cases}$$
	Put $f(m_1,m_2,m_3)=f_1(m_1)f_2(m_2)f_3(m_3).$
	
	Therefore, we may write
	\begin{align}\label{eq:ST}
		\fS_{\textrm{III}}&=\cT+O(x^{1-\varepsilon}),
	\end{align}
	where
	\begin{align*}
		\cT&=\sum_{q,r,s}\gamma_q\delta_r\lambda_s\bigg(\sum_{\substack{m_1m_2m_3n\equiv a_1\bmod {rs}\\m_1m_2m_3n\equiv a_2\bmod q}}f(m_1,m_2,m_3)\beta_n-\frac{1}{\varphi(qrs)}\sum_{(m_1m_2m_3n,qrs)=1}f(m_1,m_2,m_3)\beta_n\bigg).
	\end{align*}
	By Poisson summation it then follows that
	\begin{align*}
		\cT
		=\sum_{q,r,s,n}\frac{\gamma_q\delta_r\lambda_s\beta_n}{(qrs)^3}\sum_{1\lqs |h_i|\lqs H_i}\widehat{f_1}\Big(\frac{h_1}{qrs}\Big)\widehat{f_2}\Big(\frac{h_2}{qrs}\Big)\widehat{f_3}\Big(\frac{h_3}{qrs}\Big)\sum_{\substack{z_1,z_2,z_3\bmod {qrs}\\ z_1z_2z_3\equiv \varrho\bmod {qrs}}}\ue\Big(\frac{h_1z_1+h_2z_2+h_3z_3}{qrs}\Big)+O(1),
	\end{align*}
	where $\varrho\equiv a_1q\overline{qn}+a_2rs\overline{rsn}\bmod {qrs}$ and $H_i=QRSM_i^{\varepsilon-1}$.
	
	\subsection{A simple model problem}\label{sketch}
	The $``$simple model" in the title refers to the case $(h_1h_2h_3,qrs)=1$. Under this condition we have the identity
	\begin{align*}
		\sum_{\substack{z_1,z_2,z_3\bmod {qrs}\\ z_1z_2z_3\equiv \varrho\bmod {qrs}}}\ue\Big(\frac{h_1z_1+h_2z_2+h_3z_3}{qrs}\Big)=qrs\cdot\kl_3(\varrho h_1h_2h_3;qrs).
	\end{align*}
	This allows us to combine the three variables $h_1,h_2,h_3$ very conveniently. Put $H=H_1H_2H_3$. Then by Cauchy's inequality, there exists a weight $\eta_6$ with $\Vert\eta_6\Vert_\infty\lqs 1,$ such that
	\begin{align*}
		\cT^2&\ll \frac{x^{\varepsilon}Q^3R^3S^2}{H}\sum_{(q,r)=1}\sum_{\substack{1\lqs |h|\lqs H_1H_2H_3\\(h,qr)=1}}\Big|\sum_{\substack{(s,qrhn)=1\\(n,qr)=1}}\eta_6(s,n)\kl_3(\varrho h;qrs)\Big|^2\\
		&\ll\frac{x^\varepsilon Q^3R^3S^2}{H}\sum_{\substack{(qr,s_1s_2)=1\\(q,r)=1}}\sum_{\substack{(n_1,qrs_1)=1\\(n_2,qrs_2)=1}}\Big|\sum_{\substack{1\lqs |h| \lqs H\\(h,qrs_1s_2)=1}}\kl_3(\varrho_1 h;qrs_1)\overline{\kl_3(\varrho_2 h;qrs_2)} \Big| ,
	\end{align*}
	where $\varrho_1\equiv a_1q\overline{qn_1}+a_2rs_1\overline{rs_1n_1}\bmod {qrs_1}$ and $\varrho_2\equiv a_1q\overline{qn_2}+a_2rs_2\overline{rs_2n_2}\bmod {qrs_2}$. Applying Lemma \ref{lm:correlation}, we have that
	\begin{align}\label{eq:handle}
		\cT^2&\ll \frac{x^\varepsilon Q^3R^3S^2}{H}\Big(\sum_{\substack{q,r,s_1,s_2,n_1,n_2}}\Big(\frac{H(s_1,s_2)}{qrs_1s_2}+1\Big)(qr[s_1,s_2])^{1/2}(\varrho_1-\varrho_2,s_1,s_2)^{1/2}(\varrho_1s_2^3-\varrho_2s_1^3,qr)^{1/2}\Big)\notag\\
		&\ll x^\varepsilon\big(Q^4R^4S^{3}N+Q^{7/2}R^{7/2}S^{3}N^2+Q^{3/2}R^{3/2}S^{2}M_1M_2M_3N^2\big).
	\end{align}
	Comparing $(\ref{eq:ST})$ $(\ref{eq:handle})$ with the expected remainder $O(x^{1-\varepsilon})$, the expected bound $(\ref{eq:TypeIII})$ holds provided that
	\begin{align*}
		M_1M_2M_3\gqs x^\varepsilon\max\{Q^4R^4S^3x^{-1},Q^{7/4}R^{7/4}S^{3/2},Q^{3/2}R^{3/2}S^2\}.
	\end{align*}
	Therefore, we finish the proof of Theorem \ref{thm:TypeIII} with the case that $h$ is always coprime to $qrs$. In the next subsection we will handle possible common factors of $h$ and $qrs$, and obtain the same result as $(\ref{eq:handle})$. The reader may skip the next subsection with technical details and just use $(\ref{eq:handle})$ for simplicity.

	\subsection{Handling common factors}
	
	Let $\bh=(h_1,h_2,h_3)\in(\bZ/q\bZ)^3$ and write $c\bh:=(ch_1,ch_2,ch_3)$ for any $c\in\bR$. For $a\in(\bZ/q\bZ)^\times$ we define 
	\begin{align*}
		F(\bh,a;q):=\frac{1}{q}\sum_{\substack{z_1,z_2,z_3\bmod {q}\\ z_1z_2z_3\equiv a\bmod {q}}}\ue\Big(\frac{h_1z_1+h_2z_2+h_3z_3}{q}\Big).
	\end{align*}
	Similarly, we now turn to consider 
	\begin{align}\label{eq:T}
		\cT:=\frac{QRS}{H}\sum_{q,r,s,n}\gamma_q\delta_r\lambda_s\beta_n\sum_{\substack{1\lqs |h_i|\lqs H_i}}F(\bh,\varrho;qrs).
	\end{align}
	For $i=1,2,3$, we write
	\begin{align*}
		h_i=u_iv_iw_i,
	\end{align*}
	where $(w_i,qrs)=1$, $u_i\mid (qr)^\infty$ and $v_i\mid s^\infty$ $($i.e., $u_i$ and $v_i$ is the product of all the primes in $h_i$, with multiplicity, that also divide $qr$ and $s$, respectively$)$. We also write (recall that $qrs$ are square-free numbers)
	\begin{align*}
		q=q_1q_2,\ r=r_1r_2,\ s=s_1s_2,
	\end{align*}
	where
	\begin{align*}
		q_1r_1=\prod_{p\mid u_1u_2u_3}p=(h_1h_2h_3,qr),\ s_1=\prod_{p\mid v_1v_2v_3}p=(h_1h_2h_3,s).
	\end{align*}
	By Lemma \ref{lm:Bezout} we have that
	\begin{align}\label{eq:F}
		F(\bh,\varrho;qrs)=F(\overline{q_2r_2s}\bh,\varrho;q_1r_1)F(\overline{qrs_2}\bh,\varrho;s_1)F(\overline{q_1r_1s_1}\bh,\varrho;q_2r_2s_2).
	\end{align}
	Let $n^\flat$ be the largest square-free divisor of an integer $n\gqs 1$ $($the square-free radical of $n)$ and let $u=u_1u_2u_3,v=v_1v_2v_3,w=w_1w_2w_3.$ Notice that $(h_1h_2h_3,q_2r_2s_2)=1$, which implies
	\begin{align}\label{eq:F1}
		F(\overline{q_1r_1s_1}\bh,\varrho;q_2r_2s_2)=\kl_3\big(\overline{q_1^3r_1^3s_1^3}\varrho uvw;q_2r_2s_2\big).
	\end{align}
	
	To deal with the remaining terms on the right hand side of $(\ref{eq:F})$, we introduce two results as following (note that $q_1r_1=(u_1u_2u_3)^\flat$ and $s_1=(v_1v_2v_3)^\flat$). 
	\begin{lemma}\label{lm:F2}
		Let the notation and conditions be as above. We have that both $F(q_2r_2s\bh,\varrho;q_1r_1)$ and $F(qrs_2\bh,\varrho;s_1)$ are independent of $s_2,n,$ and 
		\begin{align*}
			|F(\overline{q_2r_2s}\bh,\varrho;q_1r_1)|\lqs \frac{u^\flat_1u^\flat_2u^\flat_3}{(q_1r_1)^2},\ |F(\overline{qrs_2}\bh,\varrho;s_1)|\lqs \frac{v^\flat_1v^\flat_2v^\flat_3}{s_1^2},
		\end{align*}
	\end{lemma}
	\proof  This is \cite[Lemma 7.3]{Po14}.\endproof
	\begin{lemma}\label{lm:F3}
		Let the notation and hypotheses be as above. we have that each of
		\begin{align*}
			\sum_{u_1,u_2,u_3\gqs 1}\frac{u^\flat_1u^\flat_2u^\flat_3}{(u_1u_2u_3)^{1/2}(q_1r_1)^3}\,,\ \sum_{v_1,v_2,v_3\gqs 1}\frac{v^\flat_1v^\flat_2v^\flat_3}{(v_1v_2v_3)^{1/2}s_1^3}
		\end{align*}
		converges to a finite value.
	\end{lemma}
	\proof  This is \cite[Lemma 7.4]{Po14}.\endproof

	Now by $(\ref{eq:T})$, $(\ref{eq:F})$, $(\ref{eq:F1})$ and Lemma \ref{lm:F2}, we have
	\begin{align}\label{eq:F4}
		\cT&\ll \frac{x^\varepsilon QRS}{H}\sum_{\substack{u_i,v_i\\i=1,2,3}}\frac{u^\flat_1u^\flat_2u^\flat_3}{(q_1r_1)^2}\,\frac{v^\flat_1v^\flat_2v^\flat_3}{s_1^2}\sum_{q_2,r_2}\sum_{\substack{0<|w|\lqs H/(uv)\\(w,q_2r_2)=1}}\bigg|\sum_{\substack{s_2,n\\(s_2,w)=1}}\eta_7(s_2,n)\kl_3\big(\overline{q_1^3r_1^3s_1^3}\varrho uvw;q_2r_2s_2\big)\bigg|\notag\\
		&:=\frac{x^\varepsilon QRS}{H}\sum_{\substack{u_i,v_i\\i=1,2,3}}\frac{u^\flat_1u^\flat_2u^\flat_3}{(q_1r_1)^2}\,\frac{v^\flat_1v^\flat_2v^\flat_3}{s_1^2}\,\widetilde{\cT}
	\end{align}
	for some $\Vert\eta_7\Vert_\infty\lqs 1.$
	Then by Cauchy's inequality we have
	\begin{align}\label{eq:F5}
		\widetilde{\cT}^2\ll \frac{x^\varepsilon H}{uv}\sum_{q_2,r_2}q_2r_2\sum_{\substack{s_2,s_3\\n_1,n_2}}\bigg|\sum_{\substack{0<|w|\lqs H/(uv)\\(w,q_2r_2s_2s_3)=1}}\kl_3\big(\overline{q_1^3r_1^3s_1^3}\varrho_1 uvw;q_2r_2s_2\big)\overline{\kl_3\big(\overline{q_1^3r_1^3s_1^3}\varrho_2 uvw;q_2r_2s_3\big)}\bigg|,
	\end{align}
	where $\varrho_1\equiv a_1q\overline{qn_1}+a_2rs_1s_2\overline{rs_1s_2n_1}\bmod {qrs_2}$ and $\varrho_2\equiv a_1q\overline{qn_2}+a_2rs_1s_3\overline{rs_1s_3n_2}\bmod {qrs_3}$.
	Applying Lemma \ref{lm:correlation}, we obtain that
	\begin{align}\label{eq:F7}
		\widetilde{\cT}^2&\ll \frac{x^\varepsilon H}{uv}\sum_{q_2,r_2}q_2r_2\sum_{\substack{s_2,s_3\\n_1,n_2}}\bigg(\frac{H/(uv)}{q_2r_2[s_2,s_3]}+1\bigg)\big(q_2r_2[s_2,s_3]\big)^{1/2}\big(
		\varrho_1-\varrho_2,s_2,s_3\big)\big(
		\varrho_1s_3^3-\varrho_2s_2^3,q_2r_2\big)\notag\\
		&\ll \frac{x^\varepsilon H}{uv}\frac{Q^2R^2}{(q_1r_1)^2}\frac{NS}{s_1}\frac{H}{uv}+\frac{x^\varepsilon Q^{3/2}R^{3/2}SH^2N^2}{(uv)^2(q_1r_1)^{3/2}s_1}+\frac{x^\varepsilon Q^{5/2}R^{5/2}S^{3}HN^2}{uv(q_1r_1)^{5/2}s_1^3}.
	\end{align}
	Finally by $(\ref{eq:F4})$-$(\ref{eq:F7})$ and Lemma \ref{lm:F3} we have
	\begin{align}\label{eq:handle1}
		\cT&\ll \frac{x^\varepsilon QRS}{H}\sum_{\substack{u_i,v_i\\i=1,2,3}}\frac{u^\flat_1u^\flat_2u^\flat_3v^\flat_1v^\flat_2v^\flat_3}{(uv)^{1/2}(q_1r_1s_1)^3}\bigg(\frac{QRS^{1/2}HN^{1/2}}{(uv/s_1)^{1/2}}+\frac{Q^{3/4}R^{3/4}S^{1/2}HN}{(uv)^{1/2}(q_1r_1s_1^2)^{-1/4}}+\frac{Q^{5/4}R^{5/4}S^{3/2}H^{1/2}N}{(q_1r_1s^2_1)^{1/4}}\bigg)\\
		&\ll x^\varepsilon\big(Q^{2}R^{2}S^{3/2}N^{1/2}+Q^{7/4}R^{7/4}S^{3/2}N+Q^{3/4}R^{3/4}S(M_1M_2M_3)^{1/2}N\big),\notag
	\end{align}
	which is the same as $(\ref{eq:handle})$.
	\begin{remark}
		The first line of $(\ref{eq:handle1})$ shows that the most difficult part is when $u,v$ are closed to $1$, which is nearly the case in Section $\ref{sketch}$ that $(h,qrs)=1$.
	\end{remark}
	
	\section{Applying the Heath-Brown identity}\label{useHB}
	
	In this section, we apply the Heath-Brown identity (Lemma \ref{lm:HB}) to connect distributions of primes in arithmetic progressions with multilinear form estimates from earlier. In order to apply the Heath-Brown identity more efficiently, we consider the following three cases separately: 
	\begin{itemize}
		\item narrow range: $Q\lqs x^{1/52};$
		\item intermediate range: $Q\lqs x^{1/35};$ 
		\item wide range: $Q\lqs x^{7/36}$.
	\end{itemize}
	The division is strategic since the narrow range yields a refined bound, while the wide range accommodates large $q\sim Q$ and allows the expected modulus $QRS$ to exceed $x^{1/2}$.
	
	To begin with, we make two new conventions in this section. We say $n$ a {\it smooth variable} if the coefficient attached to $n$ is a smooth function; say the variable $m$ is {\it S-W} if the coefficient attached to $m$ satisfies the Siegel--Walfisz condition (Definition \ref{Def:SW}).
	\subsection{The narrow range}\label{tinyQ}
	In this subsection, we follow the idea in \cite{FI83}.  Let
	\begin{align*}
		Q\lqs x^{1/52}.
	\end{align*}
	\begin{corollary}\label{coro:decomposition3}
		Fix $j$ on the right hand side of Heath-Brown identity $(\ref{eq:HB})$, then at least one of the following two conditions holds$:$
		\\[3pt]
		$(1)$ There exists a smooth variable $n_i\gqs x^{4/13};$
		\\[3pt]
		$(2)$ There is a partial product $P_1$ of $\{m_1,\dots,m_j,n_1,\dots,n_j\}$ satisfying $x^{1/13}\lqs P_1\lqs x^{1/4}.$
		\\[3pt]
		\indent Furthermore, all the variables mentioned above are S-W.
	\end{corollary}
	\proof The result can be obtained directly by taking $\theta=4/13$ in Lemma $\ref{lm:FI}$.\endproof
	
	\noindent\textbf{For condition (1)} We use the Type \uppercase\expandafter{\romannumeral1} estimate. By Theorem \ref{thm:TypeI}, we have $(\ref{eq:TypeI})$ provided that
	\begin{equation*}
		\begin{cases}
			Q\lqs R,\\
			Q^{2}R^{2}S\lqs x^{1-\varepsilon},\\
			QR^{-1}S^{-1}\lqs x^{-1-\varepsilon}M^2,\\
			Q^{3}S\lqs x^{-\varepsilon}M,\\
			Q^5R^4S\lqs x^{1-\varepsilon}M,\\
			Q^4R^3S^2\lqs x^{1-\varepsilon}M,
		\end{cases}
	\end{equation*}
	which leads to
	\begin{align*}
		x^{1+\varepsilon}M^{-2}Q^2\lqs QRS\lqs x^{-\varepsilon}\min\{x^{1/2}M^{1/2}Q^{-3/2},x^{1/4}MQ^{-5/2},x^{1/3}M^{2/3}Q^{-4/3}\}
	\end{align*}
	by choosing $S=x^{-\varepsilon}MQ^{-3}$. Now, taking the variable $m\sim M$ as the smooth variable $n_i$ in Corollary \ref{coro:decomposition3} completes the Type I estimate \eqref{eq:TypeI}, provided that
	\begin{align*}
		S=x^{-\varepsilon}MQ^{-3},\ x^{5/13+\varepsilon}Q^2\lqs QRS\lqs x^{-\varepsilon}\min\{x^{29/52}Q^{-5/2},x^{7/13}Q^{-4/3}\}.
	\end{align*}
	
	\noindent\textbf{For condition (2)} We apply the Type \uppercase\expandafter{\romannumeral2} estimate (the second method, see Section \ref{secondmethod}). Taking $S=N^{1-\varepsilon}.$ By Theorem \ref{thm:TypeII-2}, we have $(\ref{eq:TypeII})$ provided that
	$$\begin{cases}
		Q^2\lqs S,\\
		Q^{2}R^{2/3}S^{2/3}\lqs x^{1/3-\varepsilon}N^{1/3},\\
		Q^{10/7}R^{2/7}S^{2/7}\lqs x^{2/7-\varepsilon}N^{-3/7},\\
		Q^{17/12}R^{1/2}S^{1/2}\lqs x^{1/3-\varepsilon}N^{-1/6}.
	\end{cases}$$
	Taking the variable $n\sim N$ as the partial product $P_1$ mentioned in Corollary \ref{coro:decomposition3}, we obtain $(\ref{eq:TypeII})$ provided that
	\begin{align*}
		S=P_1^{1-\varepsilon},\ QRS\lqs x^{7/13-\varepsilon}Q^{-2}.
	\end{align*}
	
	Since the Bombieri--Vinogradov theorem can treat the case $QRS\lqs x^{1/2-\varepsilon}$, we conclude the following result.
	
	\begin{proposition}\label{prop:smallq}
		Let $\varepsilon>0$ and $0\lqs\theta\leqslant 1/52.$ Suppose $|a_1|\lqs \sL^B$ for some $B>0,$ $|a_2|\leqslant x$ and $a_1\neq a_2$. 
		Let $\bgamma=(\gamma_q)$ and $\blambda=(\lambda_d)$ be divisor-bounded sequences supported on positive integers with
		\begin{align*}
			q\sim x^\theta,~d\lqs x^{\cL(\theta)-\varepsilon},~(d,a_1)=(q,a_2)=\mu^2(d)=\mu^2(q)=1.
		\end{align*}
		Assume further that $\blambda$ is well-factorable of level $x^{\cL(\theta)-\varepsilon}.$ Then for  $\cL(\theta)=7/13-3\theta,$
		we have
		\begin{align*}
			\sum_{\substack{q,d\\(q,d)=1}}\gamma_q\lambda_d\bigg(\sum_{\substack{n\lqs x\\n\equiv a_1\bmod d\\n\equiv a_2\bmod q}}\Lambda(n)-\frac{1}{\varphi(qd)}\sum_{\substack{n\lqs x\\(n,qd)=1}}\Lambda(n)\bigg)\ll x\sL^{-A} .
		\end{align*}      
		for any $A>0,$ where the implied constant depends at most on $(\varepsilon,A,B).$
	\end{proposition}
	
	\begin{remark}
		When $Q<x^{\varepsilon}$, the expected level $QD$ reaches $x^{1/2+1/26-\varepsilon}$. Note that the problem reduces to that in \cite{FI83} when $Q=1$, in which case one may expect that the result could be as good as that in \cite{FI83}. In fact, a more efficient application of Cauchy's inequality yields a stronger result. Specifically, applying it in the form
		$$\sum\ldots\sum\Big|\sum_{n_2}\sum_{h}\Big|$$ gives the level $1/2+1/26,$ as opposed to the form $$\sum\ldots\sum\Big|\sum_{n_2}\Big|$$ used in  \cite{FI83}, which leads to the level $1/2+1/34.$
		
	\end{remark}

	\subsection{The intermediate range}\label{prop:smallQ}
	In this subsection, we assume
	\begin{align*}
		Q\lqs x^{1/35}.
	\end{align*}
	We begin with the following statement.
	
	\begin{corollary}\label{coro:decomposition1}
		Fix $j$ on the right hand side of Heath-Brown identity $(\ref{eq:HB})$, then at least one of the following four conditions holds$:$
		\\[3pt]
		$(1)$ There exists a smooth variable $n_i\gqs x^{8/13}Q^{-7/13};$
		\\[3pt]
		$(2)$ There is a partial product $P_2$ of $\{m_1,\dots,m_j,n_1,\dots,n_j\}$ satisfying $x^{5/13}Q^{7/13}\lqs P_2\lqs x^{6/13}Q^{-32/13};$
		\\[3pt]
		$(3)$ There is a partial product $P_3$ of $\{m_1,\dots,m_j,n_1,\dots,n_j\}$ satisfying $x^{6/13}Q^{-32/13}\lqs P_3\lqs x^{1/2};$
		\\[3pt]
		$(4)$ There exist distinct smooth variables $n_{i_1},n_{i_2},n_{i_3}\gqs 3/13-14\log Q/(13\sL),$ such that
		\begin{align*}
			n_{i_1}n_{i_2}n_{i_3}\gqs x^{12/13}Q^{-21/26}.
		\end{align*} 
		Furthermore, all the variables mentioned above are S-W.
	\end{corollary}
	
	\proof This lemma can be easily proved by choosing $\sigma= 3/26-7\log Q/(13\sL )$ in Lemma \ref{lm:combine}. For  $Q\in ~]x^{1/39},x^{1/35}]$ condition $(2)$ might not exist, in which case our proof still works in this subsection.\endproof

	\noindent\textbf{For condition (1)} We use the Type \uppercase\expandafter{\romannumeral1} estimate. By taking $m\sim M$ in Theorem \ref{thm:TypeI} as the smooth variable $n_i$ in Corollary \ref{coro:decomposition1}, we obtain \eqref{eq:TypeI} provided that
	\begin{align*}\label{eq:teshu1}
		QRS\lqs x^{8/13-\varepsilon}Q^{-7/13}.
	\end{align*}
	
	\noindent\textbf{For condition (2)} We use the Type \uppercase\expandafter{\romannumeral2} estimate (the first method, see Section \ref{firstmethod}). Taking $S=N^{1-\varepsilon}$ (recall that $MN=x$) and applying Theorem \ref{thm:TypeII-1}, we have $(\ref{eq:TypeII})$
	as soon as
	\begin{equation*}
		\begin{cases}
			Q^2\lqs R,\\
			N^2\lqs x^{1-\varepsilon},\\
			Q^{4}RS\lqs x^{1-\varepsilon}N^{-1},\\
			Q^{6}R^{5}S^{5}\lqs x^{3/2-\varepsilon}N^{3},\\
			Q^{7/2}R^{7/3}S^{7/3}\lqs x^{1-\varepsilon}N^{2/3}.
		\end{cases}
	\end{equation*}
	Taking the variable $n\sim N$ as the partial product $P_2$ (notice that $Q^3P_2\lqs x^{1/2-\varepsilon}$) mentioned above, we obtain \eqref{eq:TypeII} provided that
	\begin{align*}
		S=P_2^{1-\varepsilon}   ,\ x^{1/2-\varepsilon}\lqs QRS\lqs x^{-\varepsilon}\min\{x^{7/13}Q^{-7/13},x^{69/130}Q^{8/65},x^{7/13}Q^{-9/26}\}.
	\end{align*}
	Note here that the range $QRS\lqs x^{1/2-\varepsilon}$ is covered directly by the Bombieri--Vinogradov theorem.
	
	\noindent\textbf{For condition (3)} We use the Type \uppercase\expandafter{\romannumeral2} estimate (the third method, see Section \ref{thirdmethod}). Taking  $S=N^{1-\varepsilon}$ and applying Theorem \ref{thm:TypeII-3}, we have $(\ref{eq:TypeII})$
	as soon as
	\begin{equation*}
		\begin{cases}
			Q^4RS\lqs x^{-\varepsilon}N^2,\\
			R^{2/3}S^{2/3}\lqs x^{-\varepsilon}N,\\
			RS\lqs x^{1-\varepsilon}N^{-1},\\
			Q^{22/9}R^{17/9}S^{17/9}\lqs x^{8/9-\varepsilon}N^{1/3},\\
			Q^{40/19}R^{34/19}S^{34/19}\lqs x^{16/19-\varepsilon}N^{5/19}.
		\end{cases}
	\end{equation*}
	Taking the variable $n\sim N$ as the partial product $P_3$ (notice that $Q^4\lqs P_3$ and $(\min P_3)^{3/2}\gqs x(\max P_3)^{-1}$) mentioned above, we obtain \eqref{eq:TypeII} provided that
	\begin{align}
		S=P_3^{1-\varepsilon}   ,\ QRS\lqs x^{-\varepsilon}\min\{x^{1/2}Q,x^{7/13}Q^{-7/13}\} .
	\end{align}
	\noindent\textbf{For condition (4)} We use the Type \uppercase\expandafter{\romannumeral3} estimate. By Theorem \ref{thm:TypeIII}, we have $(\ref{eq:TypeIII})$ assuming that
	\begin{equation*}
		\begin{cases}
			\min\{M_1,M_2,M_3\}\gqs x^\varepsilon,\\
			Q^4R^4S^{3}x^{-1}\lqs (M_1M_2M_3)^{1-\varepsilon},\\
			Q^{7/4}R^{7/4}S^{3/2}\lqs (M_1M_2M_3)^{1-\varepsilon},\\
			Q^{3/2}R^{3/2}S^{2}\lqs (M_1M_2M_3)^{1-\varepsilon}.
		\end{cases}
	\end{equation*}
	Taking the variable $m_1\sim M_1,m_2\sim M_2,m_3\sim M_3$ as $n_{i_1},n_{i_2},n_{i_3}$ mentioned in Corollary \ref{coro:decomposition1}, respectively, we complete the Type III estimate \eqref{eq:TypeIII} provided that
	\begin{align*}
		S=x(QR)^{-5/2}   ,\ QRS\lqs x^{-\varepsilon}\min\{x^{1/7}(M_1M_2M_3)^{3/7},x^{-1/8}(M_1M_2M_3)^{3/4}\} \lqs x^{7/13-\varepsilon}Q^{-9/26}.
	\end{align*}
	
	In summary, we obtain the following result.
	\begin{proposition}\label{prop:largeq}
		Let $0\lqs\theta\leqslant 1/35.$ Follow the notation and other assumptions in Proposition $\ref{prop:smallq}$. Then for all $\cL(\theta)$ satisfying
		$$\cL(\theta)=\begin{cases}
			1/2, & \theta\in [0,1/40], \\
			7/13-20\theta/13, & \theta\in [1/40,1/35],
		\end{cases}$$
		we have
		\begin{align*}
			\sum_{\substack{q,d\\(q,d)=1}}\gamma_q\lambda_d\bigg(\sum_{\substack{n\lqs x\\n\equiv a_1\bmod d\\n\equiv a_2\bmod q}}\Lambda(n)-\frac{1}{\varphi(qd)}\sum_{\substack{n\lqs x\\(n,qd)=1}}\Lambda(n)\bigg)\ll x\sL^{-A} .
		\end{align*}
	\end{proposition}
	\begin{remark}
		When $Q=x^{1/40}$, the expected level $QD$ reaches $x^{1/2+1/40-\varepsilon}$.
	\end{remark}

	\subsection{The wide range}\label{prop:largeQ}
	In this subsection, we assume
	\begin{align*}
		Q\lqs x^{7/36}.
	\end{align*}
	\begin{corollary}\label{coro:decomposition2}
		Fix $j$ on the right hand side of Heath-Brown identity $(\ref{eq:HB})$, then at least one of the following three conditions holds$:$
		\\[3pt]
		$(1)$ There exists a smooth variable $n_i\gqs x^{1/2+1/9};$
		\\[3pt]
		$(2)$ There is a partial product $P_4$ of $\{m_1,\dots,m_j,n_1,\dots,n_j\}$ satisfying $ x^{7/18}\lqs P_4\lqs x^{1/2};$
		\\[3pt]
		$(3)$ There exist distinct smooth variables $n_{i_1},n_{i_2},n_{i_3}\gqs 2/9,$ such that
		\begin{align}
			n_{i_1}n_{i_2}n_{i_3}\gqs x^{3/4+1/6}.
		\end{align} 
		Furthermore, all the variables mentioned above are S-W.
	\end{corollary}
	\proof This lemma can be easily proved by choosing $\sigma= 1/9$ in Lemma \ref{lm:combine}.\endproof
	
	\noindent\textbf{For condition (1)} We use the Type \uppercase\expandafter{\romannumeral1} estimate. By taking $m\sim M$ in Theorem \ref{thm:TypeI} as the smooth variable $n_i$ in Corollary \ref{coro:decomposition1}, we have \eqref{eq:TypeI} provided that
	\begin{align}\label{eq:teshu2}
		QRS\lqs x^{11/18-\varepsilon}
	\end{align}
	to obtain $(\ref{eq:TypeI})$.
	
	\noindent\textbf{For condition (2)} The same as the previous section, we apply Theorem $\ref{thm:TypeII-3}$ to give the Type \uppercase\expandafter{\romannumeral2} estimate (the third method, see Section \ref{thirdmethod}). Taking $S=N^{1-\varepsilon},$ we have \eqref{eq:TypeII} provided that
	\begin{equation*}
		\begin{cases}
			Q^4R\lqs x^{-\varepsilon}N,\\
			R^{2/3}S^{2/3}\lqs x^{-\varepsilon}N,\\
			RS\lqs x^{1-\varepsilon}N^{-1},\\
			Q^{22/9}R^{17/9}S^{17/9}\lqs x^{8/9-\varepsilon}N^{1/3},\\
			Q^{40/19}R^{34/19}S^{34/19}\lqs x^{16/19-\varepsilon}N^{5/19},
		\end{cases}
	\end{equation*}
	or
	\begin{equation*}
		\begin{cases}
			RS\lqs x^{1-\varepsilon}N^{-1},\\
			Q^{5/2}R^{9/4}S^{9/4}\lqs x^{1-\varepsilon}N^{1/2},\\
			Q^{9/4}R^{2}S^{2}\lqs x^{1-\varepsilon}N^{1/8},\\
			Q^{32/17}R^{32/17}S^{32/17}\lqs x^{16/17-\varepsilon}N^{1/17}.
		\end{cases}
	\end{equation*}
	Taking the variable $n\sim N$ as the partial product $P_4$ mentioned above and $S=P_4^{1-\varepsilon},$ the Type II estimate \eqref{eq:TypeII} follows under the assumptions that
	\begin{align*}
		QRS\lqs x^{-\varepsilon}\min\{x^{1/2}Q,x^{19/36}Q^{-3/17}\}
	\end{align*}
	for $Q\lqs x^{17/192},$ and
	\begin{align*}
		QRS\lqs x^{-\varepsilon}\min\{x^{151/288}Q^{-1/8},x^{295/576}\}
	\end{align*}
	for $x^{17/192}<Q\lqs x^{7/36}.$
	
	\noindent\textbf{For condition (3)} Similarly we use the Type \uppercase\expandafter{\romannumeral3} estimate.  
	Taking the variables $m_1\sim M_1,m_2\sim M_2,m_3\sim M_3$ in Theorem \ref{thm:TypeIII} to be $n_{i_1},n_{i_2},n_{i_3}$ mentioned in Corollary \ref{coro:decomposition2}, respectively, we have $(\ref{eq:TypeIII})$ assuming that
	\begin{align}
		S=x(QR)^{-5/2}   ,\ QRS\lqs x^{-\varepsilon}\min\{x^{1/7}(M_1M_2M_3)^{3/7},x^{-1/8}(M_1M_2M_3)^{3/4}\} \lqs x^{15/28-\varepsilon} .
	\end{align}
	
	In summary, we obtain the following result.
	\begin{proposition}\label{prop:largerq}
		Let $0\lqs\theta\leqslant 7/36.$ Follow the notation and other assumptions in Proposition $\ref{prop:smallq}$. Then for all $\cL(\theta)$ satisfying
		$$\cL(\theta)=\begin{cases}
			1/2, & \theta\in [0,17/720], \\
			19/36-20\theta/17, & \theta\in [17/720,17/192],\\
			295/576-\theta, & \theta\in [17/192,7/72],\\
			151/288-9\theta/8, & \theta\in [7/72,7/36],
		\end{cases}$$
		we have
		\begin{align*}
			\sum_{\substack{q,d\\(q,d)=1}}\gamma_q\lambda_d\bigg(\sum_{\substack{n\lqs x\\n\equiv a_1\bmod d\\n\equiv a_2\bmod q}}\Lambda(n)-\frac{1}{\varphi(qd)}\sum_{\substack{n\lqs x\\(n,qd)=1}}\Lambda(n)\bigg)\ll x\sL^{-A} .
		\end{align*}
	\end{proposition}
	\begin{remark}
		When $Q=x^{17/720}$, the expected level $QD$ reaches $x^{1/2+17/720-\varepsilon}$.
	\end{remark}

	\section{Proof of the mean value theorems}
	\subsection{Proof of Theorem \ref{thm:meanvalue-primes}}
	By Proposition \ref{prop:smallq} we need 
	\begin{align*}
		\cL(\theta)= 7/13-3\theta,~ \text{for~}\theta\in [0,1/52].
	\end{align*}
	By Proposition \ref{prop:largeq} we need 
	$$\begin{cases}
		\cL(\theta)= 1/2, & \text{for~}\theta\in [0,1/40],\\
		\cL(\theta)\lqs 7/13-20\theta/13, & \text{for~}\theta\in [1/40,1/35].
	\end{cases}$$
	By Proposition \ref{prop:largerq} we need 
	$$\begin{cases}
		\cL(\theta)\lqs 1/2, & \text{for~}\theta\in [0,17/720],\\
		\cL(\theta)\lqs 19/36-20\theta/17, & \text{for~}\theta\in [17/720,17/192],\\
		\cL(\theta)\lqs 295/576-\theta, & \text{for~}\theta\in [17/192,7/72],\\
		\cL(\theta)\lqs 151/288-9\theta/8, & \text{for~}\theta\in [7/72,7/36].
	\end{cases}$$
	We then finish the proof by comparing the above terms.
	\begin{remark}
		One can improve this result by a more refined choice of $\sigma$ $($depending on $\log Q/\sL)$ in Section $\ref{useHB}$, at the cost of greater complexity.
	\end{remark}

	\subsection{Proof of Theorem \ref{thm:meanvalue-multilinearform}}
	For the first condition, recall that we use the Type \uppercase\expandafter{\romannumeral2} estimate (Theorem \ref{thm:TypeII-2}) to deduce $(\ref{eq:TypeII})$ by choosing $S=N^{1-\varepsilon}$, provided that $Q^2\lqs S$ and
	\begin{align*}
		RS \lqs x^{-\varepsilon}\min\{x^{1/2}N^{1/2}Q^{-3},xN^{-3/2}Q^{-5},x^{2/3}N^{-1/3}Q^{-17/6}\}.
	\end{align*}
	For the second condition, recall that we use the Type \uppercase\expandafter{\romannumeral2} estimate (Theorem \ref{thm:TypeII-1}) to deduce $(\ref{eq:TypeII})$ provided that $N\lqs x^{1/2}$ and
	\begin{align*}
		Q^2N\lqs RS \lqs x^{-\varepsilon}\min\{xN^{-1}Q^{-4},x^{3/10}N^{3/5}Q^{-6/5},x^{3/7}N^{2/7}Q^{-3/2}\}.
	\end{align*}
	For the last condition, recall that we use the Type \uppercase\expandafter{\romannumeral2} estimate (use the first condition in Theorem \ref{thm:TypeII-3}, and remove the condition $Q^4R\lqs S$ by the second condition$)$ to deduce $(\ref{eq:TypeII})$ provided that
	\begin{align*} 
		RS \lqs x^{-\varepsilon}\min\{N^{3/2},xN^{-1},x^{8/17}N^{3/17}Q^{-22/17},x^{8/17}N^{5/34}Q^{-20/17}\}.
	\end{align*}
	Finally by the Bombieri--Vinogradov theorem, which can deal with the part $\theta+\cL(\theta,\nu)<1/2$, we finish the proof.

	\section{The greatest prime factor of $p+6$ for Chen primes $p$}\label{greatest}
	
	\subsection{Linear sieve with well-factorable remainder}\label{bilinearsieve}
	We now present some preliminaries from sieve theory. Let $\sA$ be a finite set of integers. Recall the definitions that $\cP$ is a set of primes and
	\begin{align*}
		\cP(z):=\prod_{\substack{p\in\cP,p<z}}p.
	\end{align*}
	We then can define the sifting function
	\begin{align*}
		S(\sA,\cP;z)=\sum_{\substack{n\in\sA\\(n,\cP(z))=1}}1. 
	\end{align*}
	The basic input for the sieve method comes via the subsequences $\sA(d):=\{n\in\sA:n\equiv 0\bmod d\}$ and in particular, from estimates for the congruence sums
	\begin{align*}
		|\sA(d)|=\sum_{\substack{n\in\sA\\n\equiv 0\bmod d}}1.
	\end{align*}
	
	Now we define a multiplicative function $g$ (density function) supported on square-free numbers with all its prime factors belonging to $\cP$ verifying $0<g(p)<1$, and satisfying that there exists a constant $K>1$ such that
	\begin{align*}
		\frac{V(g,z_1)}{V(g,z_2)}\lqs K\frac{\log z_2}{\log z_1}\ \text{with}\ V(g,z):=\prod_{p<z,p\in \cP}(1-g(p))
	\end{align*}
	for $2\lqs z_1<z_2.$
	Let $X$ be an approximation to $|\sA|$.
	It is expected that $\sA$ is equidistributed in the arithmetic progression $n\equiv 0\bmod d$, which means that the remainder
	\begin{align*}
		r(\sA,d)=|\sA(d)|-g(d)X
	\end{align*}
	is small.
	
	Let $F$ and $f$ be the continuous solutions to the system
	$$\begin{cases}
		sF(s)=2e^{C_0}, & \text{for~}1\lqs s\lqs 3,\\
		f(s)=0, & \text{for~}0<s\lqs 2,\\
		(sF(s))'=f(s-1), & \text{for~}s>3,\\
		(sf(s))'=F(s-1), & \text{for~}s>2,
	\end{cases}$$
	where $C_0$ is the Euler constant.
	
	With the above notation, we introduce the following result of Iwaniec \cite{Iw80}.
	\begin{proposition}[Linear sieve with well-factorable remainder]\label{linearsieve}
		Let $\varepsilon>0$ be sufficiently small$,$ 
		\begin{align*}
			S(\sA,\cP;z)\lqs XV(g,z)\bigg(F\Big(\frac{\log D}{\log z}\Big)+O(\varepsilon)\bigg)+\sum_{j\lqs J(\varepsilon)}\sum_{d\mid \cP(z)}\lambda^+_{j}(d)r(\sA,d),\\
			S(\sA,\cP;z)\gqs XV(g,z)\bigg(f\Big(\frac{\log D}{\log z}\Big)+O(\varepsilon)\bigg)+\sum_{j\lqs J(\varepsilon)}\sum_{d\mid \cP(z)}\lambda^-_{j}(d)r(\sA,d),
		\end{align*}
		where the implied constant depends only on $K$, $J(\varepsilon)$ depends only on $\varepsilon$, and all $\lambda^\pm_{j}$ are well-factorable of level $D$.
	\end{proposition}
	\begin{remark}\label{linearsieveremark}
		Let $\varepsilon>0$, $D=RS$ with $R,S\gqs 1$, $u=D^{\varepsilon^2}$. For $1\lqs j\lqs J(\varepsilon):=e^{-1/\varepsilon^3}$ there are $1$-bounded sequences $(\delta^\pm_j)$,  $(\lambda^\pm_s)$, which depend on $\varepsilon,R,S$, such that the remainder terms of the linear sieve of level $D$ can be rewritten as
		\begin{align*}
			\sum_{1\lqs j\lqs J(\varepsilon)}\sum_{\substack{r\lqs R\\r\mid \cP(z)}}\sum_{\substack{s\lqs S\\(s,P(u))=1}}\delta^\pm_{j}(r)\lambda^\pm_{j}(s)r(\sA,rs).
		\end{align*}
	\end{remark}
	Note that both Proposition \ref{linearsieve} and Remark \ref{linearsieveremark} are given in \cite[Chapter 12]{FI10}. 
	
	One may observe that $F(s)$ decays to $1$ and $f(s)$ increases to $1$ .Thus in applications, one should find the admissible level $D$ as large as possible while the remainder can be controlled.

	\subsection{Weighted sieve and Chebyshev--Hooley method}\label{CH}
	Now we consider the greatest prime factor of $p+6$ for Chen primes $p$. The first step is to reformulate the information about Chen primes. Unlike the problem of $p+2\in\sP_2$, where a positive lower bound of the counting function suffices, we require a weighted sieve to ensure a large prime factor in $p+6$. 
	
	As mentioned in Section \ref{sec:shiftedChenprimes}, there are many works devoted to enhance the admissible constant in Chen's theorem \eqref{eq:Chen}. Besides the clear applications of various mean value theorems on primes (or structured sequences) in arithmetic progressions, all these works rely on careful constructions of weighted sieves, as well as Buchstab iterations. To prove Theorem \ref{thm:p+6},  we adopt the weighted sieve of Liu \cite{Li89} for simplicity and effectiveness.
	
	\begin{proposition}[Liu \cite{Li89}]\label{prop:Liu1}
		Let $z>0.$ Let $\cP$ be the set of all odd primes, and correspondingly
		\begin{align*}
			\cP(z)=\prod_{2<p<z}p.
		\end{align*}
		Let $\vartheta_1=0.1058,\vartheta_2=0.2981$. Then we have
		\begin{align*}
			\sum_{\substack{p\lqs x\\p+2\in\sP_2}}1&\gqs \sum_{\substack{p\lqs x\\(p+2,\cP(x^{\vartheta_1}))=1}}1-\frac{1}{2}\sum_{x^{\vartheta_1}\lqs p_1<x^{\vartheta_2}}\sum_{\substack{p\lqs x,\ p_1\mid p+2\\(p+2,\cP(x^{\vartheta_1}))=1}}1-\frac{1}{2}\sum_{\substack{p+2=p_1p_2p_3\lqs x\\x^{\vartheta_1}\lqs p_1<x^{\vartheta_2}\lqs p_2<p_3}}1-\sum_{\substack{p+2=p_1p_2p_3\lqs x\\x^{\vartheta_2}\lqs p_1<p_2<p_3}}1\\
			&\gqs 1.015\sC_0\frac{x}{\sL^2},
		\end{align*}
		where
		\begin{align*}
			\sC_0=2\prod_{p>2}\Big(1-\frac{1}{(p-1)^2}\Big).
		\end{align*}
	\end{proposition}

	The second step is to find the greatest prime factor. Our goal is to prove a lower bound of the shape
	\begin{align*}
		\sum_{\substack{p\lqs x\\p+2\in\sP_2}}\sum_{\substack{p_0\mid p+6\\p_0>x^{\vartheta}}}\log p_0\gqs  \varepsilon \sC_0\pi(x),
	\end{align*}
	for some $\vartheta$ as large as possible, which implies that for many Chen primes $p\lqs x$ such that $p+6$ has a prime factor $>x^{\vartheta}$. 
	Motivated by the Chebyshev--Hooley method, we need the lower bound of the full sum
	\begin{align*}
		\sum_{\substack{p\lqs x\\p+2\in\sP_2}}\sum_{\substack{p_0\mid p+6}}\log p_0\sim \sum_{\substack{p\lqs x\\p+2\in\sP_2}}\sum_{\substack{q\mid p+6}}\Lambda(q)=\sum_{\substack{p\lqs x\\p+2\in\sP_2}}\log(p+6)\gqs 1.015\sC_0\pi(x),
	\end{align*}
	which can be deduced from Proposition \ref{prop:Liu1} immediately. It then remains to give the upper bound
	\begin{align}\label{eq:application}
		\sum_{\substack{p\lqs x\\p+2\in\sP_2}}\sum_{\substack{p_0\mid p+6\\p_0\lqs x^{\vartheta}}}\log p_0=\sum_{\substack{p_0\lqs x^{\vartheta}}}\log p_0\sum_{\substack{p\lqs x\\p+2\in\sP_2\\p+6\equiv 0\bmod{p_0}}}1\lqs (1.015-\varepsilon)\sC_0\pi(x).
	\end{align}
	
	In fact, Liu's work essentially proves the following result with the weight $\log(p+6)$.
	
	\begin{corollary}\label{coro:1}
		Let $\vartheta_1,\vartheta_2,\sC_0$ follow the notation in Proposition $\ref{prop:Liu1}$. For sufficiently large $x$ we define
		\begin{align*}
			&\cC^*_1=\sum_{\substack{p\lqs x\\(p+2,\cP(x^{\vartheta_1}))=1}}\log (p+6),\\
			&\cC^*_2=\sum_{x^{\vartheta_1}\lqs p_1<x^{\vartheta_2}}\sum_{\substack{p\lqs x,\ p_1\mid p+2\\(p+2,\cP(x^{\vartheta_1}))=1}}\log(p+6),\\
			&\cC^*_3=\sum_{\substack{x^{\vartheta_1}\lqs p_1<x^{\vartheta_2}\lqs p_2<p_3\lqs x/(p_1p_2)\\(p_1p_2p_3-2,\cP(z(p_2)))=1}}\log(p_1p_2p_3+4),\\
			&\cC^*_4=\sum_{\substack{x^{\vartheta_2}\lqs p_1<p_2<p_3\lqs x/(p_1p_2)\\(p_1p_2p_3-2,\cP(z(p_1)))=1}}\log(p_1p_2p_3+4),
		\end{align*}
		where 
		\begin{align}\label{eq:Liusievechosen}
			(z(t))^3=\begin{cases}
				x^{1/2}t^{1/4}, &\text{for~}\vartheta_2\lqs \log t/\sL\lqs 2/5,\\
				xt^{-1}, & \text{for~}2/5\lqs \log t/\sL\lqs (1-\vartheta_1)/2.
			\end{cases}
		\end{align}
		Then we have
		\begin{align*}
			\sum_{\substack{p\lqs x\\p+2\in\sP_2}}\log(p+6)\gqs \cC^*_1-\frac{1}{2}\cC^*_2-\frac{1}{2}\cC^*_3-\cC^*_4\gqs 1.015\sC_0\pi(x).
		\end{align*}
	\end{corollary}

	\begin{remark}
		In \cite{Li89}, Liu chose the sieve parameter by replacing $(z(t))^3$ with $(z(t))^2$ in \eqref{eq:Liusievechosen}. Other choices of $z(t)$ with very small undulations do not affect Liu's constant $1.015$. However, the refined choice of $z(t)$ introduced in \eqref{eq:Liusievechosen} yields a better lower bound for the underlying $\cC_3$ and $\cC_4$.
	\end{remark}

	Motivated by the Chebyshev--Hooley method, we turn to consider the truncated sums:
	\begin{align}
		&\cC_1=\cC_1(\vartheta)=\sum_{p_0\lqs x^\vartheta}\log p_0\sum_{\substack{p\lqs x\\(p+2,\cP(x^{\vartheta_1}))=1\\p+6\equiv0\bmod {p_0}}}1,\label{eq:c1}\\
		&\cC_2=\cC_2(\vartheta)=\sum_{p_0\lqs x^\vartheta}\log p_0\sum_{x^{\vartheta_1}\lqs p_1<x^{\vartheta_2}}\sum_{\substack{p\lqs x\\p+2\equiv 0\bmod {p_1}\\(p+2,\cP(x^{\vartheta_1}))=1\\p+6\equiv 0\bmod {p_0}}}1,\label{eq:c2}\\
		&\cC_3=\cC_3(\vartheta)=\sum_{p_0\lqs x^\vartheta}\log p_0\sum_{\substack{x^{\vartheta_1}\lqs p_1<x^{\vartheta_2}\lqs p_2<p_3\lqs x/(p_1p_2)\\(p_1p_2p_3-2,\cP(z(p_2)))=1\\p_1p_2p_3+4\equiv 0\bmod {p_0}}}1,\label{eq:c3}\\
		&\cC_4=\cC_4(\vartheta)=\sum_{p_0\lqs x^\vartheta}\log p_0\sum_{\substack{x^{\vartheta_2}\lqs p_1<p_2<p_3\lqs x/(p_1p_2)\\(p_1p_2p_3-2,\cP(z(p_1)))=1\\p_1p_2p_3+4\equiv 0\bmod {p_0}}}1.\label{eq:c4}
	\end{align}
	Therefore, the problem is reduced to finding a $\vartheta$ as large as possible, while the following inequality holds:
	\begin{align*}
		\cC_1-\frac{1}{2}\cC_2-\frac{1}{2}\cC_3-\cC_4\lqs (1.015-\varepsilon)\sC_0\pi(x).
	\end{align*}
	
	In fact, for $\vartheta=0.217$ we have
	\begin{align}\label{eq:inequality1}
		\cC_1\lqs 1.239859\sC_0\pi(x),
	\end{align}
	\begin{align}\label{eq:inequality2}
		\cC_2\gqs 0.367562\sC_0\pi(x),
	\end{align}
	\begin{align}\label{eq:inequality3}
		\cC_3\gqs 0.078732\sC_0\pi(x),
	\end{align}
	\begin{align}\label{eq:inequality4}
		\cC_4\gqs 0.003035\sC_0\pi(x).
	\end{align}
	
	In the next four subsections we will prove the inequalities above. For sieve arguments, we will follow the notation in Section \ref{bilinearsieve} and use $q$ to denote the prime variable $p_0$ in $\cC_1,\cC_2,\cC_3,\cC_4.$

	\subsection{Estimation of $\cC_1$}\label{1}
	Recall that $\vartheta_1=0.1058$, $\vartheta_2=0.2981$, and the definition of functions $F$ and $f$. In particular, we only need 
	$$\begin{cases}
		F(s)=2e^{C_0}s^{-1}, & \text{for~}1\lqs s\lqs 3,\\
		F(s)=2e^{C_0}s^{-1}\big(1+\int_{2}^{s-1}t^{-1}\log(t-1)\ud t\big), & \text{for~}3\lqs s\lqs 5,\\
		f(s)=2e^{C_0}s^{-1}\log(s-1), & \text{for~}2\lqs s\lqs 4,\\
		f(s)=2e^{C_0}s^{-1}\big(\log(s-1)+\int_{3}^{s-1}t^{-1}\ud t\int_2^{t-1}u^{-1}\log(u-1)\ud u\big), & \text{for~}4\lqs s\lqs6.
	\end{cases}$$
	
	We write
	\begin{align*}
		\cC_1=\sum_{q\lqs x^{7/36}}+\sum_{x^{7/36}<q\lqs x^{\vartheta}}:=\cC_{11}+\cC_{12},
	\end{align*}
	and let
	\begin{align}
		\sA_{q}=\{p+2:p\lqs x,\,q\mid p+6\},\ \cP_{q}=\{p:p\nmid 2q\},\ X_{q}=\pi(x)/\varphi(q). 
	\end{align}
	For fixed $q\lqs x^{\vartheta},$ we employ the linear sieve to detect the condition $(p+2,\cP(x^{\vartheta_1}))=1$ for primes $p$ with $q\mid p+6.$
	For any square-free numbers $d$ we put
	$$g_q(d)=\begin{cases}
		\varphi^{-1}(d), & \text{if~}d\mid \cP(x^{\vartheta_1})\ \text{and~}(d,q)=1,\\
		0, & \text{otherwise}.
	\end{cases}$$
	Correspondingly, $r(\sA_q,d)=|\sA_q(d)|-g_q(d)X_q$. Applying Proposition \ref{linearsieve} with an upper bound sieve $(\lambda_d)$ of level $D,$
	we have
	\begin{align*}
		S(\sA_q,\cP_q;x^{\vartheta_1})\lqs X_qV(g_q,x^{\vartheta_1})\Big(F\Big(\frac{\log D}{\vartheta_1\sL}\Big)+O(\varepsilon)\Big)+\sum_{\substack{d\mid \cP_q(x^{\vartheta_1})\\d\lqs D}}\lambda_d r(\sA_q,d).
	\end{align*}
	Summing the above terms over $q\lqs x^{7/36}$, together with the prime number theorem, yields the following estimate, where the remainder is handled by Theorem \ref{thm:meanvalue-primes} (for $y\gqs 5$ we replace $F(y)$ by $F(5)$):
	\begin{align*}
		\cC_{11}&\lqs\sum_{q\lqs x^{7/36}}(\log q) X_{q}V(g_q,x^{\vartheta_1})\bigg(F\Big(\frac{\cL(\log q/\sL )}{\vartheta_1}\Big)+O(\varepsilon)\bigg)+O(x\sL^{-A})\\
		&=(1+O(\varepsilon))\sC_0\frac{\pi(x)}{e^{C_0}\vartheta_1}\int_{0}^{7/36}F(\cL(t)/\vartheta_1)\ud t\lqs 1.086489\sC_0\pi(x).
	\end{align*}
	Similarly, for $\cC_{12},$ we apply the Bombieri--Vinogradov theorem to handle the remainder, and to get
	\begin{align*}
		\cC_{12}\lqs (1+O(\varepsilon))\sC_0\frac{\pi(x)}{e^{C_0}\vartheta_1}\int_{7/36}^{\vartheta}F((1/2-t)/\vartheta_1)\ud t\lqs 0.153370\sC_0\pi(x).
	\end{align*}
	Concluding the arguments in this subsection, we have $(\ref{eq:inequality1})$.

	\subsection{Estimation of $\cC_2$}\label{2}
	Let
	\begin{align*}
		\cC_2\gqs &\sum_{q\lqs x^{7/36}}\sum_{x^{\vartheta_1}\lqs p_1<x^{1/4}/q^{1/2}}+\sum_{q\lqs x^{7/36}}\sum_{x^{1/4}/q^{1/2}\lqs p_1< x^{\vartheta_2}}:=\cC_{21}+\cC_{22},
	\end{align*}
	and differently
	\begin{align*}
		\sA_{q,p_1}=\{p+2:p\lqs x,\,q\mid p+6,\,p_1\mid p+2\},\ \cP_{q,p_1}=\{p:p\nmid 2qp_1\},\ X_{q,p_1}=\pi(x)/\varphi(qp_1).
	\end{align*}
	
	Notice that there is an extra summation (over $p_1$) in $\cC_2$ compared to $\cC_1$, forcing us to recombine the variables before applying Theorem \ref{thm:meanvalue-primes}. To this end, we introduce the following Corollary, which can be deduced from Theorem \ref{thm:meanvalue-primes} directly.
	
	\begin{corollary}\label{thm:meanvalue-primes-coro}
		Under the assumptions of Theorem $\ref{thm:meanvalue-primes}$, we further assume that $\brho=(\rho_t)$ is a divisor-bounded sequence supported on $ t \leq x^{\sT(\theta)} $ with $(t, da_1) = 1$ and $\mu^2(t)=1$. When $\sT(\theta)\lqs \cL(\theta)$ and
		$$\cL(\theta)+\sT(\theta)=\begin{cases}
			7/13-3\theta, & \text{for~}\theta\in [0,1/78], \\
			1/2, & \text{for~}\theta\in [1/78,1/40], \\
			7/13-20\theta/13, & \text{for~}\theta\in [1/40,1/35], \\
			19/36-20\theta/17, & \text{for~}\theta\in [1/35,17/192], \\
			295/576-\theta, & \text{for~}\theta\in [17/192,7/72],\\
			151/288-9\theta/8, & \text{for~}\theta\in [7/72,7/36],	
		\end{cases}$$
		we have
		\begin{align}\label{eq:coro102}
			\mathop{\sum_{q}\sum_{d}\sum_{t}}_{(q,dt)=1}\gamma_q\lambda_d\rho_t\bigg(\sum_{\substack{n\lqs x\\n\equiv a_1\bmod {dt}\\n\equiv a_2\bmod {q}}}\Lambda(n)-\frac{1}{\varphi(qdt)}\sum_{\substack{n\lqs x\\(n,qdt)=1}}\Lambda(n)\bigg)\ll x\sL^{-A}
		\end{align}
		\noindent for any $A>0$.
	\end{corollary}
	
	\proof By the arguments in \cite[Lemma 5]{FG89}, the Dirichlet convolution $\blambda*\brho$ is well-factorable of level $\sT(\theta) \cL(\theta)$ since $\sT(\theta)\lqs \cL(\theta)$. Now the proof can be deduced immediately by taking $\lambda_d$ in Theorem \ref{thm:meanvalue-primes} to be the Dirichlet convolution of $\lambda_d$ and $\rho_t$ in Corollary \ref{thm:meanvalue-primes-coro}. \endproof

	Now we apply Corollary \ref{thm:meanvalue-primes-coro} to treat the remainder from $\cC_{21},$ and the Bombieri--Vinogradov theorem to treat the remainder from $\cC_{22}$ and $\cC_{23}$. Notice that the condition $p_1<x^{1/4}/q^{1/2}$ is not crucial compared to $\sT(\theta)\lqs \cL(\theta)$ in Corollary \ref{thm:meanvalue-primes-coro}, but can simplify the calculation. We have (in Corollary \ref{thm:meanvalue-primes-coro}, $\cL(\theta)$ depends on $\sT(\theta)$; likewise, in the equations below, $\cL(t)$ depends on $u$)
	\begin{align*}
		\cC_{21}&\gqs \sum_{\substack{q\lqs x^{7/36}\\x^{\vartheta_1}\lqs p_1<x^{1/4}/{q}^{1/2}}}\frac{\sC_0\pi(x)\log q}{e^{C_0}\vartheta_1(q-1)(p_1-1)\sL}\bigg(f\Big(\frac{\cL(\log q/\sL)}{\vartheta_1}\Big)+O(\varepsilon)\bigg)+O\Big(\frac{x}{\sL^A}\Big)\\
		&\gqs(1+O(\varepsilon))\sC_0\frac{\pi(x)}{e^{C_0}\vartheta_1}\int_{0}^{7/36}\ud t\int_{\vartheta_1}^{1/4-t/2}\frac{f(\cL(t)/\vartheta_1)}{u}\ud u\\
		&\gqs 0.361841\pi(x),\\
		\cC_{22}&\gqs(1+O(\varepsilon))\sC_0\frac{\pi(x)}{e^{C_0}\vartheta_1}\int_{0}^{7/36}\ud t\int_{1/4-t/2}^{\vartheta_2}\frac{f((1/2-t-u)/\vartheta_1)}{u}\ud u\gqs 0.005721\pi(x).
	\end{align*}
	Concluding the arguments in this subsection, we have $(\ref{eq:inequality2})$.

	\subsection{Estimation of $\cC_3$}\label{3}
	Let
	\begin{align*}
		\cC_3=\sum_{q\lqs x^{1/60}}+\sum_{x^{1/60}<q\lqs x^{2/23}}+\sum_{x^{2/23}<q\lqs x^{\vartheta}}:=\cC_{31}+\cC_{32}+\cC_{33},
	\end{align*}
	and differently (more specifically, we need to split the sum over $p_2$ into bisection intervals like $p_2\sim P'_2$ with some $P'_2>0$)
	\begin{align*}
		&\sA_q=\{p_1p_2p_3-2:x^{\vartheta_1}\lqs p_1<x^{\vartheta_2}\lqs p_2<p_3\lqs x/(p_1p_2),\,q\mid p_1p_2p_3+4\},\\ 
		&X_q=\frac{\pi(x)}{\varphi(q)}\int_{\vartheta_1}^{\vartheta_2}\frac{\ud u}{u}\int_{\vartheta_2}^{(1-u)/2}\frac{\ud v}{v(1-u-v)}. 
	\end{align*}
	By Proposition \ref{linearsieve}, the Bombieri--Vinogradov theorem and the first part of Theorem \ref{thm:meanvalue-multilinearform}, we  obtain
	\begin{align*}
		\cC_{31}&\gqs (1+O(\varepsilon))\sum_{q\lqs x^{1/60}}\frac{\log q}{e^{C_0}(q-1)}\sC_0\pi(x)\int_{\vartheta_1}^{\vartheta_2}\frac{\ud u}{u}\int_{\vartheta_2}^{(1-u)/2}\frac{f\big(\cL(t,u)\sL/(\log z(v))\big)}{(\log z(v))v(1-u-v)}\ud v \\
		&\gqs(1+O(\varepsilon))\sC_0\frac{\pi(x)}{e^{C_0}}\Big\{\int_{0}^{1/60}\ud t\int_{\vartheta_1}^{1/5}\frac{\ud u}{u}\int_{\vartheta_2}^{2/5}\frac{f(12\cL(t,u)/(2+v))}{v(1-u-v)(1/6+v/12)}\ud v\\
		&\qquad\qquad\qquad\qquad\quad\ \ +\int_{0}^{1/60}\ud t\int_{\vartheta_1}^{1/5}\frac{\ud u}{u}\int_{2/5}^{(1-u)/2}\frac{f(3\cL(t,u)/(1-v))}{v(1-u-v)(1/3-v/3)}\ud v\\
		&\qquad\qquad\qquad\qquad\quad\ \ +\int_{0}^{1/60}\ud t\int_{1/5}^{\vartheta_2}\frac{\ud u}{u}\int_{\vartheta_2}^{(1-u)/2}\frac{f((6-12t)/(2+v))}{v(1-u-v)(1/6+v/12)}\ud v\Big\}\\
		&\gqs (0.013336+0.003157+0.006364)\sC_0\pi(x)=0.022857\sC_0\pi(x).
	\end{align*}
	Similarly, applying the Bombieri--Vinogradov theorem and the third part of Theorem \ref{thm:meanvalue-multilinearform}, we obtain
	\begin{align*}
		\cC_{32}&\gqs(1+O(\varepsilon))\sC_0\frac{\pi(x)}{e^{C_0}}\Big\{\int_{1/60}^{2/23}\ud t\int_{\vartheta_1}^{\vartheta_2}\frac{\ud u}{u}\int_{\vartheta_2}^{8/23}\frac{f((6-12t)/(2+v))}{v(1-u-v)(1/6+v/12)}\ud v\\
		&\qquad\qquad\qquad\qquad\quad\ \ +\int_{1/60}^{2/23}\ud t\int_{\vartheta_1}^{1/5}\frac{\ud u}{u}\int_{8/23}^{2/5}\frac{f(12\cL(t,u)/(2+v))}{v(1-u-v)(1/6+v/12)}\ud v\\
		&\qquad\qquad\qquad\qquad\quad\ \ +\int_{1/60}^{2/23}\ud t\int_{\vartheta_1}^{1/5}\frac{\ud u}{u}\int_{2/5}^{(1-u)/2}\frac{f(3\cL(t,u)/(1-v))}{v(1-u-v)(1/3-v/3)}\ud v\\
		&\qquad\qquad\qquad\qquad\quad\ \ +\int_{1/60}^{2/23}\ud t\int_{1/5}^{\vartheta_2}\frac{\ud u}{u}\int_{8/23}^{(1-u)/2}\frac{f(12\cL(t,u)/(2+v))}{v(1-u-v)(1/6+v/12)}\ud v\Big\}\\
		&\gqs (0.027684+0.009611+0.005196+0.005438)\sC_0\pi(x)=0.047929\sC_0\pi(x),
	\end{align*}
	\begin{align*}
		\cC_{33}&\gqs(1+O(\varepsilon))\sC_0\frac{\pi(x)}{e^{C_0}}\Big\{\int_{2/23}^{\vartheta}\ud t\int_{\vartheta_1}^{1/5}\frac{\ud u}{u}\int_{\vartheta_2}^{2/5}\frac{f((6-12t)/(2+v))}{v(1-u-v)(1/6+v/12)}\ud v\\
		&\qquad\qquad\qquad\qquad\quad\ \ +\int_{2/23}^{\vartheta}\ud t\int_{\vartheta_1}^{1/5}\frac{\ud u}{u}\int_{2/5}^{(1-u)/2}\frac{f((3/2-3t)/(1-v))}{v(1-u-v)(1/3-v/3)}\ud v\\
		&\qquad\qquad\qquad\qquad\quad\ \ +\int_{2/23}^{\vartheta}\ud t\int_{1/5}^{\vartheta_2}\frac{\ud u}{u}\int_{\vartheta_2}^{(1-u)/2}\frac{f((6-12t)/(2+v))}{v(1-u-v)(1/6+v/12)}\ud v\Big\}\\ 
		&\gqs (0.003995+0.001219+0.002732)\sC_0\pi(x)=0.007946\sC_0\pi(x).
	\end{align*}
	Concluding the arguments in this subsection we have $(\ref{eq:inequality3})$.

	\subsection{Estimation of $\cC_4$}\label{4}
	Very similar to Section \ref{3}, 
	let
	\begin{align*}
		\cC_4=\sum_{q\lqs x^{1/30}}+\sum_{x^{1/30}<q\lqs x^{\vartheta}}:=\cC_{41}+\cC_{42}.
	\end{align*}
	Then by Proposition \ref{linearsieve}, the Bombieri--Vinogradov theorem and the second part of Theorem \ref{thm:meanvalue-multilinearform}, we have
	\begin{align*}
		\cC_{41}&\gqs (1+O(\varepsilon))\sum_{q\lqs x^{1/30}}\frac{\log q}{e^{C_0}(q-1)}\sC_0\pi(x)\int_{\vartheta_1}^{\vartheta_2}\frac{\ud u}{u}\int_{\vartheta_2}^{(1-u)/2}\frac{f\big(\cL(t,u)\sL/(\log z(v))\big)}{\log(z(v))v(1-u-v)}\ud v \\
		&\gqs(1+O(\varepsilon))\sC_0\frac{\pi(x)}{e^{C_0}}\Big\{\int_{0}^{1/30}\ud t\int_{\vartheta_2}^{14/45}\frac{\ud u}{u}\int_{31/90}^{(1-u)/2}\frac{f(12\cL(t,v)/(2+v))}{v(1-u-v)(1/6+v/12)}\ud v\\
		&\qquad\qquad\qquad\qquad\quad\ \ +\int_{0}^{1/30}\ud t\int_{\vartheta_2}^{14/45}\frac{\ud u}{u}\int_{u}^{31/90}\frac{f((6-12t)/(2+v))}{v(1-u-v)(1/6+v/12)}\ud v\\
		&\qquad\qquad\qquad\qquad\quad\ \ +\int_{0}^{1/30}\ud t\int_{14/45}^{1/3}\frac{\ud u}{u}\int_{u}^{(1-u)/2}\frac{f((6-12t)/(2+v))}{v(1-u-v)(1/6+v/12)}\ud v\Big\}\\
		&\gqs (0.000063+0.000783+0.000550)\sC_0\pi(x)=0.001396\sC_0\pi(x),
	\end{align*}
	\begin{align*}
		\cC_{42}&\gqs(1+O(\varepsilon))\sC_0\frac{\pi(x)}{e^{C_0}}\int_{1/30}^{\vartheta}\ud t\int_{\vartheta_2}^{1/3}\frac{\ud u}{u}\int_{u}^{(1-u)/2}\frac{f((6-12t)/(2+v))}{v(1-u-v)(1/6+v/12)}\ud v\\
		&\gqs 0.001639\sC_0\pi(x).
	\end{align*}
	Concluding the arguments in this subsection we have $(\ref{eq:inequality4})$.
	\begin{remark}
		In Sections $\ref{3}$ and $\ref{4}$ we use Chen's switching principle, which succeed in switching the problem of sieving $p+2$ to be $p_1p_2p_3$ into the problem of sieving $p_1p_2p_3-2$ to be primes.
	\end{remark}

	\section{Comments and remarks}
	
	\subsection{Possible refinements of Theorems \ref{thm:meanvalue-primes} and \ref{thm:meanvalue-multilinearform}}
	
	Theorem \ref{thm:meanvalue-primes} deals with a well-factorable weight $\blambda$. If $\blambda$ comes from upper-bound linear sieves or satisfies triply-well-factorable properties, it is very possible that the new ideas in \cite{Ma25b,Li23,Pa24,Pa25a} can be adopted to enhance the levels of distributions in Theorems \ref{thm:meanvalue-primes} and \ref{thm:meanvalue-multilinearform}.
	
	Grimmelt and Merikoski \cite{GM24} developed a new framework that spectral theory of Poincar\'e series over congruence subgroups can be applied to study determinant equations with high uniformity in the modulus, instead of stepping via sums of Kloosterman sums as before. In an ongoing project, their new insights are applied to study distributions of primes in simultaneous arithmetic progressions as in Theorem \ref{thm:meanvalue-primes}. As noted in \cite{GM24}, their argument works in the slightly stronger form
	\begin{align*}
		\sum_{q}\gamma_q\max_{(a_2,q)=1}\sum_{(d,q)=1}\lambda_d\Big(\sum_{\substack{p\lqs x\\p\equiv a_1\bmod d\\p\equiv a_2\bmod q}}1-\frac{\pi(x)}{\varphi(qd)}\Big),
	\end{align*}  
	where $\blambda=(\lambda_d)$ is well-factorable. We believe that their approach allows one to go further with larger support of $\bgamma=(\gamma_q).$

	\subsection{Brun--Titchmarsh theorem in simultaneous arithmetic progressions}
	By the Chinese remainder theorem, Theorem \ref{thm:meanvalue-primes} provides an asymptotic formula for the counting function $\pi(x;dq,\fa)$ on average over $d,q$, where
	$\fa\equiv a_1q\overline{q}+a_2d\overline{d}\bmod{dq}$. Instead, one may also study the upper bound for $\pi(x;dq,\fa)$ on average over $d,q,$ in the spirit of many previous works on the Brun--Titchmarsh theorem on average. Moreover, the underlying ideas should be helpful to produce larger exponents in Theorem \ref{thm:p+6} for shifted Chen primes. We will revisit this issue in the forthcoming works.

	\subsection{Smooth numbers in simultaneous arithmetic progressions}
	In place of primes, it is also fundamental to study smooth numbers in arithmetic progressions. Define
	\begin{align*}
		S(x,y)=\{n\lqs x: p\mid n\Rightarrow p\lqs y\}.
	\end{align*}
	The problem of determining when the asymptotic formula
	\begin{align*}
		\sum_{\substack{n\in S(x,y)\\n\equiv a\bmod{q}}}1\sim\frac{1}{\varphi(q)}\sum_{\substack{n\in S(x,y)\\(n,q)=1}}1
	\end{align*}
	holds, depending on $q,x,y,$ has been studied in \cite{Fr81,FT91,BP92,Gr93a,Gr93b,So08,Ha12b}. In many applications, it suffices that this property holds for almost all $q$ within a prescribed range. Results in this direction $($see \cite{Wo73a,Wo73b,Gr93a,FT96,Ha12a}$),$ which provide an averaged estimate over $q,$ are of the Bombieri--Vinogradov type
	\begin{align}\label{eq:smoothaverage}
		\sum_{q\lqs x^{1/2-\varepsilon}}\max_{(a,q)=1}\Big|\sum_{\substack{n\in S(x,y)\\n\equiv a\bmod{q}}}1-\frac{1}{\varphi(q)}\sum_{\substack{n\in S(x,y)\\(n,q)=1}}1\Big|\ll_{\varepsilon,A}|S(x,y)|\sL^{-A}.
	\end{align}
	By sacrificing uniformity in the residue class $a,$ one aims to surpass the $x^{1/2}$-barrier in \eqref{eq:smoothaverage}. Fouvry--Tenenbaum \cite{FT96} enhanced the exponent to $3/5$ for a fixed residue $a$. This was recently improved by Drappeau \cite{Dr15} with much higher uniformity of $a$ compared to $q$. More recently, Pascadi \cite{Pa25b} obtained a larger admissible exponent $66/107$ for a fixed residue class $a.$ Analogous to Theorem \ref{thm:meanvalue-primes}, one may expect to have the estimate
	\begin{align*}
		\sum_{(q,d)=1}\gamma_q\lambda_d\Big(\sum_{\substack{n\in S(x,y)\\n\equiv a_1\bmod{d}\\n\equiv a_2\bmod{q}}}1-\frac{1}{\varphi(qd)}\sum_{\substack{n\in S(x,y)\\(n,q)=1}}1\Big)\ll_{\varepsilon,A}|S(x,y)|\sL^{-A}
	\end{align*}
	for fixed $a_1\neq a_2$, where the level of $\bgamma*\blambda$ reaches $x^{1/2+\delta}$ with some small $\delta>0$. This is the starting point of our further works.

	\bibliographystyle{plain}

\end{document}